\documentclass[11pt]{article}
\usepackage{amssymb}
\usepackage{amsmath}
\usepackage[all]{xy}
\usepackage{times}
\usepackage{graphicx}

\title{On the complexity of Borel equivalence relations\\ with some countability property\indent}
\author{Dominique LECOMTE}
\date{\today}

\def\ufootnote#1{\let\savedthfn\thefootnote\let\thefootnote\relax
\footnote{#1}\let\thefootnote\savedthfn\addtocounter{footnote}{-1}}

\newcommand{\Ana}{{\it\Sigma}^{1}_{1}}

\newcommand{\Ca}{{\it\Pi}^{1}_{1}}
\newcommand{\Boraone}{{\it\Sigma}^{0}_{1}}

\newcommand{\Borel}{{\it\Delta}^{1}_{1}}
\newcommand{\ana}{{\bf\Sigma}^{1}_{1}}
\newcommand{\borel}{{\bf\Delta}^{1}_{1}}

\newcommand{\boraone}{{\bf\Sigma}^{0}_{1}}
\newcommand{\boratwo}{{\bf\Sigma}^{0}_{2}}

\newcommand{\boraxi}{{\bf\Sigma}^{0}_{\xi}}

\newcommand{\borone}{{\bf\Delta}^{0}_{1}}
\newcommand{\bortwo}{{\bf\Delta}^{0}_{2}}

\newcommand{\bormone}{{\bf\Pi}^{0}_{1}}

\newcommand{\bormtwo}{{\bf\Pi}^{0}_{2}}

\newcommand{\bormep}{{\bf\Pi}^{0}_{\eta +1}}

\newcommand{\bormxi}{{\bf\Pi}^{0}_{\xi}}

\newcommand{\borxi}{{\bf\Delta}^{0}_{\xi}}
\newcommand{\borme}{{\bf\Pi}^{0}_{\eta}}

\newcommand{\boraep}{{\bf\Sigma}^{0}_{\eta +1}}

\newtheorem{thm} {Theorem} [section]

\newtheorem{cor} [thm] {Corollary}

\newtheorem{lem} [thm] {Lemma}
\newtheorem{prop} [thm] {Proposition}
\newtheorem{them} {Theorem} [subsection]
\newtheorem{defin} [them] {Definition}

\newtheorem{lemm} [them] {Lemma}

\begin{document}

\maketitle

\centerline{$\bullet$ Sorbonne Universit\' e, Institut de Math\'ematiques de Jussieu-Paris Rive Gauche,}

\centerline{CNRS, Universit\'e Paris Diderot, Projet Analyse Fonctionnelle,}

\centerline{Campus Pierre et Marie Curie, Case 247, 4, place Jussieu, 75 252 Paris cedex 5, France}

\centerline{dominique.lecomte@upmc.fr}\bigskip

\centerline{$\bullet$ Universit\'e de Picardie, I.U.T. de l'Oise, site de Creil,}

\centerline{13, all\'ee de la fa\"\i encerie, 60 107 Creil, France}\bigskip\bigskip\bigskip\bigskip\bigskip

\ufootnote{{\it 2010 Mathematics Subject Classification.}~Primary: 03E15, Secondary: 28A05, 54H05}

\ufootnote{{\it Keywords and phrases.}~Borel class, Borel equivalence relation, descriptive complexity, Borel reducibiity, continuous reducibility, dichotomy}

\noindent {\bf Abstract.} We study the class of Borel equivalence relations under continuous reducibility. In particular, we characterize when a Borel equivalence relation with countable equivalence classes is $\boraxi$ (or $\bormxi$). We characterize when all the equivalence classes of such a relation are $\boraxi$ (or 
$\bormxi$). We prove analogous results for the Borel equivalence relations with countably many equivalence classes. We also completely solve these two problems for the first two ranks. In order to do this, we prove some extensions of the Louveau-Saint Raymond theorem which itself generalized the Hurewicz theorem characterizing when a Borel subset of a Polish space is $G_\delta$.

\vfill\eject

\section{$\!\!\!\!\!\!$ Introduction}\indent

 The present paper is about descriptive set theory, which is the study of definable subsets of Polish spaces (recall that a topological space is {\bf Polish} if it is separable and completely metrizable). The reader should see [K1] for the standard descriptive set theoretic notions and notation. The most classical hierarchy of topological complexity in descriptive set theory is the one given by the Borel classes. If $\bf\Gamma$ is a class of subsets of the metrizable spaces, then $\check {\bf\Gamma}\! :=\!\{\neg S\mid S\!\in\! {\bf\Gamma}\}$ is its {\bf dual class}. Recall that the Borel hierarchy is the inclusion from left to right in the following picture:\bigskip
 
\scalebox{0.75}{$$\!\!\!\!\!\!\!\!\!\!\!\!\!\!\!\!\!\!\!\!\!\!\!\xymatrix@1{ 
& & \boraone\! =\!\mbox{open} & & \boratwo\! =\! F_\sigma & & & 
\boraxi\! =\! (\bigcup_{\eta <\xi}~\borme )_\sigma & \\ 
& \borone\! =\!\mbox{clopen} & & \bortwo\! =\!\boratwo\cap\bormtwo & & \cdots & \borxi\! =\!\boraxi\cap\bormxi & & \cdots\\ 
& & \bormone\! =\!\mbox{closed} & & \bormtwo\! =\! G_\delta & & & \bormxi\! =\!\check\boraxi & }$$}\bigskip
 
\noindent This hierarchy is strict in uncountable Polish spaces, in which the non self-dual classes are those of the form $\boraxi$ or $\bormxi$. In the sequel, by non self-dual Borel class, we mean exactly those classes.\bigskip

 The study of {\bf Borel equivalence relations} under Borel reducibilty is one of the major topics in descriptive set theory since more than three decades now. Several important dichotomy results concerning the Borel equivalence relations have been proved (see, for example, [S], [Ha-K-Lo], [H-K]). They are of the following form: a relation is either simple, or more complicated than a typical complicated relation. Several quasi-orders have been used to compare the Borel equivalence relations (recall that a {\bf quasi-order} is a reflexive and transitive relation). The most common is {\bf Borel reducibility}. Recall that if $X,Y$ are topological (or standard Borel) spaces and $E\!\subseteq\! X^2$, 
$F\!\subseteq\! Y^2$, 
$$(X,E)\leq_B(Y,F)\Leftrightarrow\exists f\! :\! X\!\rightarrow\! Y\mbox{ Borel with }E\! =\! (f\!\times\! f)^{-1}(F)$$
(we say that $f$ {\bf reduces} $E$ to $F$). However, very early in the theory, the quasi-order $\sqsubseteq_c$ of {\bf injective continuous reducibility} defined by 
$$(X,E)\sqsubseteq_c(Y,F)\Leftrightarrow\exists f\! :\! X\!\rightarrow\! Y\mbox{ injective continuous with }E\! =\! (f\!\times\! f)^{-1}(F)$$
has also been considered, for example in the main result of [S].

\begin{thm} (Silver) Let $E$ be a co-analytic equivalence relation on a Polish space $X$. Then exactly one of the following holds:\smallskip

(a) the relation $E$ has countably many equivalence classes,\smallskip

(b) $(2^\omega ,=)\sqsubseteq_c(X,E)$.\end{thm}

 The quasi-order $\leq_c$ of {\bf continuous reducibility} can also be mentioned. We are interested in the descriptive complexity of Borel equivalence relations on Polish spaces. In order to approach this problem, it is useful to consider invariants for the considered quasi-order. In the context of Borel relations on a Polish space, a natural invariant for Borel reducibility has been studied, the notion of potential complexity (see, for example, [L2], [L3], and [Lo2] for the definition). A Borel relation $R$ on a Polish space $X$ is {\bf potentially} in a Wadge class $\bf\Gamma$ if we can find a finer Polish topology $\tau$ on $X$ such that $R$ is in $\bf\Gamma$ in the product $(X,\tau )^2$.
 
\vfill\eject
 
 This is an invariant in the sense that any relation which is Borel reducible to a relation potentially in 
$\bf\Gamma$ has also to be potentially in $\bf\Gamma$. Along similar lines, any relation which is continuously reducible to a relation in $\bf\Gamma$ has also to be in $\bf\Gamma$. Moreover, the pre-image of an equivalence relation by a square map is an equivalence relation, which is not the case with arbitrary continuous maps. This motivates the  work in the present paper. We are looking for characterizations of the Borel equivalence relations either in a fixed Borel class $\bf\Gamma$, or whose equivalence classes are in $\bf\Gamma$. So we will consider the continuous and injective continuous reducibilities. In other words, we want to give answers to the following very simple questions.\bigskip

\noindent\bf Questions.\rm\ (1) When is a Borel equivalence relation ${\bf\Sigma}^0_\xi$ (or ${\bf\Pi}^0_\xi$)?\smallskip

\noindent (2) When are the equivalence classes of a Borel equivalence relation ${\bf\Sigma}^0_\xi$ (or ${\bf\Pi}^0_\xi$)?\bigskip

 Question (1) is the most natural one. Question (2) is also natural, in particular when we think about classical uniformization results for instance (see Section 18 in [K1]). As we will see, it turns out that the solution to Question (2) is an important step towards the solution to Question (1) (see Theorem \ref{contB}). There are several possible approaches to try to solve these problems. One can try an approach ``from above", which means finding a relation universal for  (i.e., above for the considered quasi-order) the relations in 
$\bf\Gamma$. For instance, it is known that there is a universal $K_\sigma$ equivalence relation for Borel reducibility (see [R]). It is an open and difficult problem to find a universal $F_\sigma$ equivalence relation for Borel reducibility, and thus for continuous reducibility also. In this paper, we will follow another approach, ``from below", which means that we will prove dichotomies of the form above. In particular, we will be able to characterize the $F_\sigma$ equivalence relations this way. We provide a complete solution for the Borel equivalence relations with some countability property, namely those with countably many equivalence classes and those with countable equivalence classes. In order to describe this, we now introduce, for some Borel classses $\bf\Gamma$ and some natural numbers $n$, useful examples of complex equivalence relations 
$\mathbb{E}^{\bf\Gamma}_n\!\notin\! {\bf\Gamma}$.\bigskip

\noindent\bf Notation.\rm\ Let $\bf\Gamma$ be a non self-dual Borel class, $\mathbb{K}$ be a metrizable compact space, and $\mathbb{C}\!\in\!\check {\bf\Gamma}(\mathbb{K})\!\setminus\! {\bf\Gamma}$.\bigskip 

 If the {\bf rank} of $\bf\Gamma$ is one (i.e., if ${\bf\Gamma}\!\in\!\{\boraone ,\bormone\}$), then we set 
$\mathbb{K}\! :=\!\{ 0\}\cup\{ 2^{-k}\mid k\in\omega\}\!\subseteq\!\mathbb{R}$, $\mathbb{C}\! :=\!\{ 0\}$ if ${\bf\Gamma}\! =\!\boraone$, and 
$\mathbb{C}\! :=\!\mathbb{K}\!\setminus\!\{ 0\}$ if ${\bf\Gamma}\! =\!\bormone$, since we want some injectivity results.\bigskip

 If the rank of $\bf\Gamma$ is at least two, then we set $\mathbb{K}\! :=\! 2^\omega$, and 
$\mathbb{C}\cap N_s\!\in\!\check {\bf\Gamma}(N_s)\!\setminus\! {\bf\Gamma}$ for each $s\!\in\! 2^{<\omega}$ (we will check that this is possible). In particular, $\mathbb{C}$ is dense and co-dense in $2^\omega$. We set 
$$\mathbb{C}\! :=\!\{\alpha\!\in\! 2^\omega\mid\exists^\infty n\!\in\!\omega ~~\alpha (n)\! =\! 1\}$$ 
if ${\bf\Gamma}\! =\!\boratwo$, and 
$\mathbb{C}\! :=\!\{\alpha\!\in\! 2^\omega\mid\forall^\infty n\!\in\!\omega ~~\alpha (n)\! =\! 0\}$ if  
${\bf\Gamma}\! =\!\bormtwo$, for injectivity reasons again. In the sequel, we will say that $\mathbb{K}$ or $\mathbb{C}$ is {\bf as above} if it satisfies all the properties mentioned here.\bigskip

\noindent\bf Examples.\rm\ We introduce two equivalence relations on 
$\mathbb{K}$. We set
$$\begin{array}{ll}
& x~\mathbb{E}^{\bf\Gamma}_0~y\Leftrightarrow (x,y\!\in\!\mathbb{C})\vee (x\! =\! y)\mbox{,}\cr\cr
& x~\mathbb{E}^{\bf\Gamma}_1~y\Leftrightarrow (x,y\!\in\! \mathbb{C})\vee (x,y\!\notin\!\mathbb{C}).
\end{array}$$
Note that $\mathbb{E}^{\bormone}_1\! =\!\mathbb{E}^{\bormone}_0$.
 
\vfill\eject
 
\setlength{\unitlength}{0.8cm}
\begin{picture}(6.0,5.5)(-5.0,-0.0)

   \put(6,2){\line(1,1){3}}
   \put(6,2){\line(1,0){3}}
   \put(6,2){\line(0,1){3}}
   \put(6,5){\line(1,0){3}}
   \put(9,2){\line(0,1){3}}
   \put(6,1){\rule[8mm]{1.6cm}{1.6cm}}
   \put(8,3){\rule[8mm]{0.8cm}{0.8cm}}
   \put(6.85,1.5){$\mathbb{C}$}
   \put(8.2,1.5){$\neg\mathbb{C}$}
   \put(5.5,2.85){$\mathbb{C}$}
   \put(5.25,4.35){$\neg\mathbb{C}$}
   \put(7.45,0.5){$\mathbb{E}^{\bf\Gamma}_1$}
   
   \put(-1,2){\line(1,1){3}}
   \put(-1,2){\line(1,0){3}}
   \put(-1,2){\line(0,1){3}}
   \put(-1,5){\line(1,0){3}}
   \put(2,2){\line(0,1){3}}
   \multiput(1.1,4)(0.2,0){5}{\circle*{0.05}}
   \multiput(1,4)(0,0.2){5}{\circle*{0.05}}
   \put(-1,1){\rule[8mm]{1.6cm}{1.6cm}}
   \put(-0.15,1.5){$\mathbb{C}$}
   \put(1.2,1.5){$\neg\mathbb{C}$}
   \put(-1.5,2.85){$\mathbb{C}$}
   \put(-1.75,4.35){$\neg\mathbb{C}$}
   \put(0.45,0.5){$\mathbb{E}^{\bf\Gamma}_0$}
   
\end{picture}

 We solve Question (2) for the Borel classes of rank at most two. Recall that if $(Q,\leq )$ is a quasi-ordered class, then a {\bf basis} is a subclass $B$ of $Q$ such that any element of $Q$ is $\leq$-above an element of $B$. We are looking for basis as small as possible, so in fact for antichains (an {\bf antichain} is a subclass of $Q$ made of pairwise $\leq$-incomparable elements). So we want antichain basis. We set 
$${\mathcal A}^{\bf\Gamma}\! :=\!\left\{\!\!\!\!\!\!\!
\begin{array}{ll}
& \{ (\mathbb{K},\mathbb{E}^{\bf\Gamma}_0)\}\mbox{ if }{\bf\Gamma}\! =\!\bormone\mbox{,}\cr\cr
& \{ (\mathbb{K},\mathbb{E}^{\bf\Gamma}_n)\mid n\!\leq\! 1\}\mbox{ if }{\bf\Gamma}\! =\!\boraone
\mbox{ or the rank of }{\bf\Gamma}\mbox{ is two.}
\end{array}
\right.$$
Most of our results will hold in analytic spaces and not only in Polish spaces. Recall that a separable metrizable space is an {\bf analytic space} if it is homeomorphic to an analytic subset of a Polish space.

\begin{thm} \label{main1} Let ${\bf\Gamma}$ be a non self-dual Borel class of rank at most two, 
$\mathbb{K},\mathbb{C}$ as above, $X$ be an analytic space, and $E$ be a Borel equivalence relation on $X$. Then exactly one of the following holds:\smallskip  

(a) the equivalence classes of $E$ are in $\bf\Gamma$,\smallskip  

(b) there is $(\mathbb{X},\mathbb{E})\!\in\! {\mathcal A}^{\bf\Gamma}$ such that 
$(\mathbb{X},\mathbb{E})\sqsubseteq_c(X,E)$.\smallskip

\noindent Moreover, ${\mathcal A}^{\bf\Gamma}$ is a $\leq_c$-antichain (and thus a 
$\sqsubseteq_c$ and a $\leq_c$-antichain basis).\end{thm}

 In order to state our results concerning Question (1), we introduce some other examples of complex equivalence relations.\bigskip

\noindent\bf Examples.\rm\ We define three equivalence relations on 
$\mathbb{H}\! :=\! 2\!\times\!\mathbb{K}$. We set  
$$\begin{array}{ll}
& (\varepsilon ,x)~\mathbb{E}^{\bf\Gamma}_3~(\eta ,y)\Leftrightarrow 
(\varepsilon ,x)\! =\! (\eta ,y)\vee (x\! =\! y\!\in\!\mathbb{C})\mbox{,}\cr\cr
& (\varepsilon ,x)~\mathbb{E}^{\bf\Gamma}_4~(\eta ,y)\Leftrightarrow (\varepsilon ,x)\! =\! (\eta ,y)\vee (x\! =\! y\!\in\!\mathbb{C})\vee 
(\eta\! =\!\varepsilon\! =\! 1\wedge x,y\!\notin\!\mathbb{C})\mbox{,}\cr\cr
& (\varepsilon ,x)~\mathbb{E}^{\bf\Gamma}_5~(\eta ,y)\Leftrightarrow (\varepsilon ,x)\! =\! (\eta ,y)\vee (x\! =\! y\!\in\!\mathbb{C})\vee 
(\eta\! =\!\varepsilon\wedge x,y\!\notin\!\mathbb{C}).
\end{array}$$
Note that $\mathbb{E}^{\bormone}_4\! =\!\mathbb{E}^{\bormone}_5\! =\!\mathbb{E}^{\bormone}_3$. We solve Question (1) for the Borel classes of rank at most two. We set 
$${\mathcal B}^{\bf\Gamma}\! :=\! {\mathcal A}^{\bf\Gamma}\cup
\left\{\!\!\!\!\!\!\!
\begin{array}{ll}
& \emptyset\mbox{ if }{\bf\Gamma}\! =\!\boraone\mbox{,}\cr\cr
& \{ (\mathbb{H},\mathbb{E}^{\bf\Gamma}_3)\}\mbox{ if }{\bf\Gamma}\! =\!\bormone\mbox{,}\cr\cr
& \{ (\mathbb{H},\mathbb{E}^{\bf\Gamma}_n)\mid 3\!\leq\! n\!\leq\! 5\}\mbox{ if the rank of }{\bf\Gamma}\mbox{ is two.}
\end{array}
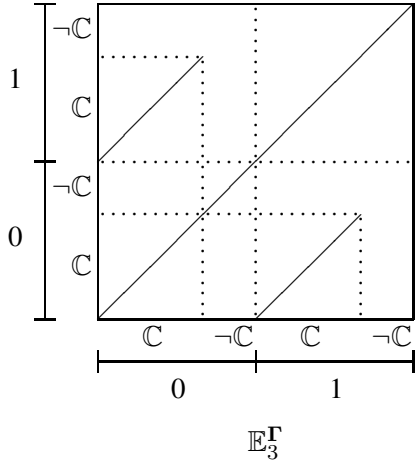
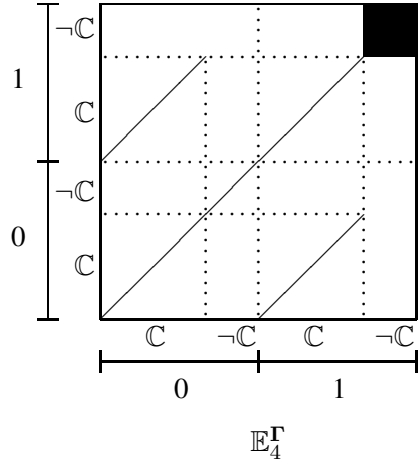
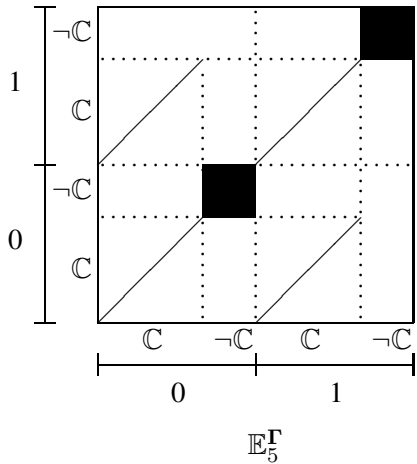
\right.$$

\setlength{\unitlength}{0.7cm}
\begin{picture}(6,6)(-2,-6)
   
   \multiput(-0.9,-4)(0.2,0){30}{\circle*{0.05}}
   \put(-1,-7){\line(1,1){6}}
   \put(-1,-7){\line(1,0){6}}
   \put(-1,-7){\line(0,1){6}}
   \put(-1,-1){\line(1,0){6}}
   \put(5,-7){\line(0,1){6}}
   \multiput(-0.9,-5)(0.2,0){25}{\circle*{0.05}}
   \multiput(1,-6.9)(0,0.2){25}{\circle*{0.05}}
   \multiput(2,-6.9)(0,0.2){30}{\circle*{0.05}}
   \multiput(-0.9,-2)(0.2,0){10}{\circle*{0.05}}
   \multiput(4,-6.9)(0,0.2){10}{\circle*{0.05}}
   \put(-1,-4){\line(1,1){2}}
   \put(2,-7){\line(1,1){2}}
   \put(-0.15,-7.5){$\mathbb{C}$}
   \put(1.2,-7.5){$\neg\mathbb{C}$}
   \put(2.85,-7.5){$\mathbb{C}$}
   \put(4.2,-7.5){$\neg\mathbb{C}$}
   \put(-1,-7.8){\line(1,0){6}}
   \put(-1,-8){\line(0,1){0.4}}
   \put(2,-8){\line(0,1){0.4}}
   \put(5,-8){\line(0,1){0.4}}
   \put(0.4,-8.5){0}
   \put(3.4,-8.5){1}
   \put(-1.5,-6.15){$\mathbb{C}$}
   \put(-1.87,-4.65){$\neg\mathbb{C}$}
   \put(-1.5,-3.15){$\mathbb{C}$}
   \put(-1.87,-1.65){$\neg\mathbb{C}$}
   \put(-2,-7){\line(0,1){6}}
   \put(-2.2,-7){\line(1,0){0.4}}
   \put(-2.2,-4){\line(1,0){0.4}}
   \put(-2.2,-1){\line(1,0){0.4}}
   \put(-2.7,-5.6){0}
   \put(-2.7,-2.6){1}
   \put(1.85,-9.5){$\mathbb{E}^{\bf\Gamma}_3$}
   
   \multiput(-0.9,-14)(0.2,0){10}{\circle*{0.05}}
   \multiput(2.1,-14)(0.2,0){15}{\circle*{0.05}}
   \put(-1,-17){\line(1,1){6}}
   \put(-1,-17){\line(1,0){6}}
   \put(-1,-17){\line(0,1){6}}
   \put(-1,-11){\line(1,0){6}}
   \put(5,-17){\line(0,1){6}}
   \multiput(-0.9,-15)(0.2,0){10}{\circle*{0.05}}
   \multiput(2.1,-15)(0.2,0){10}{\circle*{0.05}}
   \multiput(1,-16.9)(0,0.2){10}{\circle*{0.05}}
   \multiput(1,-14.1)(0,0.2){10}{\circle*{0.05}}
   \multiput(2,-16.9)(0,0.2){10}{\circle*{0.05}}
   \multiput(2,-13.9)(0,0.2){15}{\circle*{0.05}}
   \multiput(-0.9,-12)(0.2,0){25}{\circle*{0.05}}
   \multiput(4,-16.9)(0,0.2){25}{\circle*{0.05}}
   \put(-1,-14){\line(1,1){2}}
   \put(2,-17){\line(1,1){2}}
   \put(4,-13.15){\rule[8mm]{0.7cm}{0.7cm}}
   \put(1,-16.15){\rule[8mm]{0.7cm}{0.7cm}}
   \put(-0.15,-17.5){$\mathbb{C}$}
   \put(1.2,-17.5){$\neg\mathbb{C}$}
   \put(2.85,-17.5){$\mathbb{C}$}
   \put(4.2,-17.5){$\neg\mathbb{C}$}
   \put(-1,-17.8){\line(1,0){6}}
   \put(-1,-18){\line(0,1){0.4}}
   \put(2,-18){\line(0,1){0.4}}
   \put(5,-18){\line(0,1){0.4}}
   \put(0.4,-18.5){0}
   \put(3.4,-18.5){1}
   \put(-1.5,-16.15){$\mathbb{C}$}
   \put(-1.87,-14.65){$\neg\mathbb{C}$}
   \put(-1.5,-13.15){$\mathbb{C}$}
   \put(-1.87,-11.65){$\neg\mathbb{C}$}
   \put(-2,-17){\line(0,1){6}}
   \put(-2.2,-17){\line(1,0){0.4}}
   \put(-2.2,-14){\line(1,0){0.4}}
   \put(-2.2,-11){\line(1,0){0.4}}
   \put(-2.7,-15.6){0}
   \put(-2.7,-12.6){1}
   \put(1.85,-19.5){$\mathbb{E}^{\bf\Gamma}_5$}

  \multiput(11.1,-4)(0.2,0){30}{\circle*{0.05}}
   \put(11,-7){\line(1,1){6}}
   \put(11,-7){\line(1,0){6}}
   \put(11,-7){\line(0,1){6}}
   \put(11,-1){\line(1,0){6}}
   \put(17,-7){\line(0,1){6}}
   \multiput(11.1,-5)(0.2,0){25}{\circle*{0.05}}
   \multiput(13,-6.9)(0,0.2){25}{\circle*{0.05}}
   \multiput(14,-6.9)(0,0.2){30}{\circle*{0.05}}
   \multiput(11.1,-2)(0.2,0){25}{\circle*{0.05}}
   \multiput(16,-6.9)(0,0.2){25}{\circle*{0.05}}
   \put(11,-4){\line(1,1){2}}
   \put(14,-7){\line(1,1){2}}
   \put(16,-3.15){\rule[8mm]{0.7cm}{0.7cm}}
   \put(11.85,-7.5){$\mathbb{C}$}
   \put(13.2,-7.5){$\neg\mathbb{C}$}
   \put(14.85,-7.5){$\mathbb{C}$}
   \put(16.2,-7.5){$\neg\mathbb{C}$}
   \put(11,-7.8){\line(1,0){6}}
   \put(11,-8){\line(0,1){0.4}}
   \put(14,-8){\line(0,1){0.4}}
   \put(17,-8){\line(0,1){0.4}}
   \put(12.4,-8.5){0}
   \put(15.4,-8.5){1}
   \put(10.5,-6.15){$\mathbb{C}$}
   \put(10.13,-4.75){$\neg\mathbb{C}$}
   \put(10.5,-3.25){$\mathbb{C}$}
   \put(10.13,-1.65){$\neg\mathbb{C}$}
   \put(10,-7){\line(0,1){6}}
   \put(9.8,-7){\line(1,0){0.4}}
   \put(9.8,-4){\line(1,0){0.4}}
   \put(9.8,-1){\line(1,0){0.4}}
   \put(9.3,-5.6){0}
   \put(9.3,-2.6){1}
   \put(13.85,-9.5){$\mathbb{E}^{\bf\Gamma}_4$}
     
\end{picture}\bigskip\bigskip\bigskip\bigskip\bigskip\bigskip\bigskip\bigskip\bigskip\bigskip\bigskip\bigskip
\bigskip\bigskip\bigskip\bigskip\bigskip\bigskip\bigskip\bigskip\bigskip\bigskip\bigskip

\begin{thm} \label{main2} Let ${\bf\Gamma}$ be a non self-dual Borel class of rank at most two, 
$\mathbb{K},\mathbb{C}$ as above, $X$ be an analytic space, and $E$ be a Borel equivalence relation on $X$. Then exactly one of the following holds:\smallskip  

(a) the relation $E$ is a $\bf\Gamma$ subset of $X^2$,\smallskip  

(b) there is $(\mathbb{X},\mathbb{E})\!\in\! {\mathcal B}^{\bf\Gamma}$ such that $(\mathbb{X},\mathbb{E})\sqsubseteq_c(X,E)$.\smallskip

\noindent Moreover, ${\mathcal B}^{\bf\Gamma}$ is a $\leq_c$-antichain (and thus a $\sqsubseteq_c$ and a $\leq_c$-antichain basis).\end{thm}

 In particular, this characterizes the $F_\sigma$ equivalence relations, from below.\bigskip

\noindent\bf Remarks.\rm\ (1) This result contrasts with the case of potentially open equivalence relations. Indeed, by Silver's theorem, if $E$ is a Borel equivalence relation on a Polish space $X$, then either $E$ is potentially open, or $\big( 2^\omega ,\Delta (2^\omega )\big)\sqsubseteq_c (X,E)$ (where 
$\Delta (X)\! :=\!\{ (x,y)\!\in\! X^2\mid x\! =\! y\}$ is the {\bf diagonal} of $X$, see [S]). So there is a $\sqsubseteq_c$-minimum non-potentially open Borel equivalence relation, and no $\sqsubseteq_c$-minimum non-open Borel equivalence relation.\smallskip

\noindent (2) This result also contrasts with the case of potentially closed equivalence relations. Indeed, by the Harrington-Kechris-Louveau theorem, if $E$ is a Borel equivalence relation on a Polish space $X$, then either $E$ is potentially closed, or $(2^\omega ,\mathbb{E}_0)\sqsubseteq_c (X,E)$ (where 
$$\mathbb{E}_0\! :=\!\{ (\alpha ,\beta )\!\in\! (2^\omega)^2\mid\forall^\infty n\!\in\!\omega ~~\alpha (n)\! =\!\beta (n)\}\mbox{,}$$ 
see [Ha-K-Lo]). So there is a $\sqsubseteq_c$-minimum non-potentially closed Borel equivalence relation, and no 
$\sqsubseteq_c$-minimum non-closed Borel equivalence relation.\smallskip

\noindent (3) As mentioned in [C-L-M], there is no equivalence relation which is $\leq_B$-minimum among non-potentially in $\bf\Gamma$ Borel equivalence relations if ${\bf\Gamma}\!\supseteq\!\boratwo$ is a Borel class. Theorem \ref{main2} gives a $\sqsubseteq_c$-antichain basis among non-$\boratwo$ Borel equivalence relations. This leads  to the question of knowing whether there is a $\leq_B$-antichain basis among non-potentially $\boratwo$ Borel equivalence relations.\bigskip

 We now turn our attention to the class $\mathfrak C$ of countable Borel equivalence relations. Recall that an 
{\bf equivalence relation} is {\bf countable} if all its equivalence classes are countable. The class 
$\mathfrak C$ has been widely studied (see, for instance, [J-K-Lo], [K2]). This class is extremely big. For instance, Adams and Kechris proved in [A-K] that we can embed the quasi-order of inclusion on the Borel subsets of $\mathbb{R}$ into the quasi-order of Borel reducibility on $\mathfrak C$. Moreover, the Feldman-Moore theorem (see 18.16 in [K1]) says that such relations are induced by a Borel action of a countable group, and the study of Borel actions of Polish groups is currently a very active area of research. Theorem \ref{main2} solves Question (1) for the Borel classes of rank at most two. Our main result, which solves Question (1) for the other Borel classes in the case of countable equivalence relations (and in fact more), is as follows.
   
\begin{thm} \label{ctble} Let ${\bf\Gamma}$ be a non self-dual Borel class of rank at least three, 
$\mathbb{C}$ as above, $X$ be an analytic space, and $E$ be a Borel equivalence relation on $X$ with 
$F_\sigma$ classes. Then exactly one of the following holds:\smallskip  

(a) the relation $E$ is a $\bf\Gamma$ subset of $X^2$,\smallskip  

(b) $(\mathbb{H},\mathbb{E}^{\bf\Gamma}_3)\sqsubseteq_c(X,E)$.\end{thm}

 Theorem \ref{ctble} can be extended to the first ranks, using Theorem \ref{main2}. The set $\{ (\mathbb{H},\mathbb{E}^{\bf\Gamma}_3)\}$ has to be replaced with 
$$\left\{\!\!\!\!\!\!\!\!
\begin{array}{ll}
& \{ (\mathbb{K},\mathbb{E}^{\bf\Gamma}_0),(\mathbb{K},\mathbb{E}^{\bf\Gamma}_1)\}\mbox{ if }
{\bf\Gamma}\! =\!\boraone\mbox{,}\cr\cr
& \{ (\mathbb{K},\mathbb{E}^{\bf\Gamma}_0),(\mathbb{H},\mathbb{E}^{\bf\Gamma}_3)\}\mbox{ if }
{\bf\Gamma}\!\in\!\{\bormone ,\bormtwo\}\mbox{,}\cr\cr
& \{ (\mathbb{H},\mathbb{E}^{\bf\Gamma}_n)\mid 3\!\leq\! n\!\leq\! 5\}\mbox{ if }
{\bf\Gamma}\! =\!\boratwo\mbox{,}
\end{array}
\right.$$
since $E$ has $F_\sigma$ classes. So we completely solved Questions (1) and (2) for countable equivalence relations.\bigskip

 In the case of Borel reducibility, the Borel equivalence relations with countably many equivalence classes are trivial. We can pick a point in each equivalence class, which Borel reduces such a relation to $(\kappa ,=)$, where $\kappa\!\leq\!\omega$ is an ordinal, and the reduction works in both directions. The situation is much more complicated in the case of continuous reducibility. We solve Question (2) for these relations. 

\begin{thm} \label{C1denclPi} Let $\xi\!\geq\! 1$ be a countable ordinal, $\mathbb{K}$ as above, 
$\mathbb{C}\!\in\!\boraxi (\mathbb{K})$ not in $\bormxi$ (as above if $\xi\!\leq\! 2$), $X$ be an analytic space, and $E$ be a Borel equivalence relation on $X$ with countably many classes. Then exactly one of the following holds:\smallskip  

(a) the equivalence classes of $E$ are $\bormxi$,\smallskip  

(b) $(\mathbb{K},\mathbb{E}^{\bormxi}_1)\sqsubseteq_c(X,E)$.\end{thm}

 If ${\bf\Gamma}\! =\!\boraxi$ with $\xi\!\geq\! 3$, then we will have to consider another equivalence relation on $2^\omega$. We can write 
$\neg\mathbb{C}\! =\!\bigcup_{n\in\omega}~\mathbb{C}_n$, where $(\mathbb{C}_n)_{n\in\omega}$ is a sequence of pairwise disjoint $\borxi$ sets (which will not be arbitrary and be given by Theorem \ref{Hur+} to come, as we will see). We set
$$x~\mathbb{E}^{\boraxi}_2~y\Leftrightarrow 
(x,y\!\in\!\mathbb{C})\vee (\exists n\!\in\!\omega ~~x,y\!\in\!\mathbb{C}_n).$$
Note that $\mathbb{E}^{\boraxi}_0\!\subseteq\!\mathbb{E}^{\boraxi}_2\!\subseteq\!\mathbb{E}^{\boraxi}_1$. Also, we set $\mathbb{E}^{\boraxi}_2\! :=\!\mathbb{E}^{\boraxi}_0$ if $\xi\!\leq\! 2$ since the 
$\mathbb{C}_n$'s are singletons in this case.

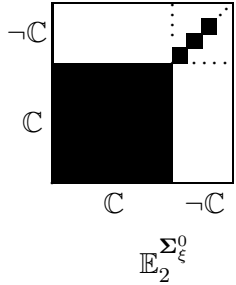
\setlength{\unitlength}{0.8cm}
\begin{picture}(6.0,6.0)(6.0,-0.0)
   
   \multiput(13,2)(0.1,0.1){30}{\circle*{0.05}}
   \put(13,2){\line(1,0){3}}
   \put(13,2){\line(0,1){3}}
   \put(13,5){\line(1,0){3}}
   \put(16,2){\line(0,1){3}}
   \put(13,1){\rule[8mm]{1.6cm}{1.6cm}}
   \put(15,3){\rule[8mm]{0.2cm}{0.2cm}}
   \put(15.24,3.24){\rule[8mm]{0.2cm}{0.2cm}}
   \put(15.48,3.48){\rule[8mm]{0.2cm}{0.2cm}}
   \multiput(15.36,4)(0.17,0){4}{\circle*{0.05}}
   \multiput(15,4.46)(0,0.17){4}{\circle*{0.05}}
   \put(13.85,1.5){$\mathbb{C}$}
   \put(15.2,1.5){$\neg\mathbb{C}$}
   \put(12.5,2.85){$\mathbb{C}$}
   \put(12.25,4.35){$\neg\mathbb{C}$}
   \put(14.45,0.5){$\mathbb{E}^{\boraxi}_2$}
   
\end{picture}

\begin{thm} \label{c1dencls} Let $\xi\!\geq\! 1$ be a countable ordinal, $\mathbb{K}$ as above, 
$\mathbb{C}\!\in\!\bormxi (\mathbb{K})$ not in $\boraxi$ (as above if $\xi\!\leq\! 2$), $X$ be an analytic space, and $E$ be a Borel equivalence relation on $X$ with countably many classes. Then exactly one of the following holds:\smallskip  

(a) the equivalence classes of $E$ are $\boraxi$ (exactly when $E$ is a $\boraxi$ subset of $X^2$),\smallskip  

(b) there is $n\!\in\!\{ 1,2\}$ such that $(\mathbb{K},\mathbb{E}^{\boraxi}_n)\sqsubseteq_c(X,E)$.\smallskip

\noindent Moreover, $\{ (\mathbb{K},\mathbb{E}^{\boraxi}_n)\mid 1\!\leq\! n\!\leq\! 2\}$ is a 
$\leq_c$-antichain (and thus a $\sqsubseteq_c$ and a $\leq_c$-antichain basis).\end{thm}

 Note that Theorem \ref{c1dencls} characterizes when a Borel equivalence relation with countably many classes is $\boraxi$. In order to finish the study of Borel equivalence relations with countably many classes, it remains to characterize those which are not $\bormxi$ if $\xi\!\geq\! 3$. The partition 
$(\mathbb{C}_n)_{n\in\omega}$ of $\mathbb{C}$ into $\borxi$ subsets of $2^\omega$ allows us to define an equivalence relation on $2\!\times\! 2^\omega$ by 
$$(\varepsilon ,\alpha )~\mathbb{E}^{\bormxi}_8~(\eta ,\beta )\Leftrightarrow (\exists n\!\in\!\omega ~~\alpha ,\beta\!\in\!\mathbb{C}_n)\vee 
(\eta\! =\!\varepsilon\wedge\alpha ,\beta\!\notin\!\mathbb{C})$$ 
(we use the number 8 here because we can consider some examples $\mathbb{E}^{\bf\Gamma}_n$ for
$n\!\in\!\{ 6,7,8\}$, in the spirit of those for $n\!\in\!\{ 3,4,5\}$ respectively, to state a general conjecture that we will not give here).

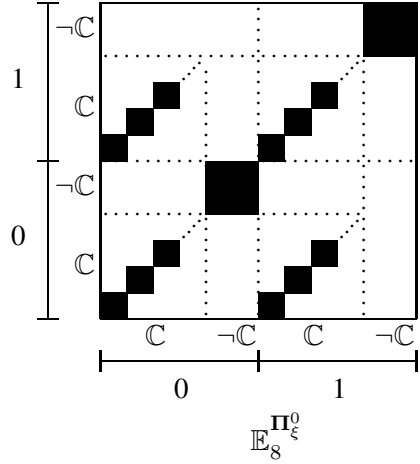
\vfill\eject

\setlength{\unitlength}{0.7cm}
\begin{picture}(1,1)(4.0,-22.0)

 \multiput(11.1,-24)(0.2,0){10}{\circle*{0.05}}
   \multiput(14.1,-24)(0.2,0){15}{\circle*{0.05}}
   \multiput(12.58,-25.44)(0.1,0.1){4}{\circle*{0.05}}
   \multiput(12.58,-22.44)(0.1,0.1){4}{\circle*{0.05}}
   \multiput(15.58,-22.44)(0.1,0.1){4}{\circle*{0.05}}
   \multiput(15.58,-25.44)(0.1,0.1){4}{\circle*{0.05}}
   \put(11,-27){\line(1,0){6}}
   \put(11,-27){\line(0,1){6}}
   \put(11,-21){\line(1,0){6}}
   \put(17,-27){\line(0,1){6}}
   \multiput(11.1,-25)(0.2,0){10}{\circle*{0.05}}
   \multiput(14.1,-25)(0.2,0){10}{\circle*{0.05}}
   \multiput(13,-26.9)(0,0.2){10}{\circle*{0.05}}
   \multiput(13,-24.1)(0,0.2){10}{\circle*{0.05}}
   \multiput(14,-26.9)(0,0.2){10}{\circle*{0.05}}
   \multiput(14,-23.9)(0,0.2){15}{\circle*{0.05}}
   \multiput(11.1,-22)(0.2,0){25}{\circle*{0.05}}
   \multiput(16,-26.9)(0,0.2){25}{\circle*{0.05}}
   \put(16,-23.15){\rule[8mm]{0.7cm}{0.7cm}}
   \put(13,-26.15){\rule[8mm]{0.7cm}{0.7cm}}
   \put(11.85,-27.5){$\mathbb{C}$}
   \put(13.2,-27.5){$\neg\mathbb{C}$}
   \put(14.85,-27.5){$\mathbb{C}$}
   \put(16.2,-27.5){$\neg\mathbb{C}$}
   \put(11,-27.8){\line(1,0){6}}
   \put(11,-28){\line(0,1){0.4}}
   \put(14,-28){\line(0,1){0.4}}
   \put(17,-28){\line(0,1){0.4}}
   \put(12.4,-28.5){0}
   \put(15.4,-28.5){1}
   \put(10.5,-26.15){$\mathbb{C}$}
   \put(10.13,-24.65){$\neg\mathbb{C}$}
   \put(10.5,-23.15){$\mathbb{C}$}
   \put(10.13,-21.65){$\neg\mathbb{C}$}
   \put(10,-27){\line(0,1){6}}
   \put(9.8,-27){\line(1,0){0.4}}
   \put(9.8,-24){\line(1,0){0.4}}
   \put(9.8,-21){\line(1,0){0.4}}
   \put(9.3,-25.6){0}
   \put(9.3,-22.6){1}
   \put(14,-25.15){\rule[8mm]{0.35cm}{0.35cm}}
   \put(14.5,-24.65){\rule[8mm]{0.35cm}{0.35cm}}
   \put(15,-24.15){\rule[8mm]{0.35cm}{0.35cm}}
   \put(11,-25.15){\rule[8mm]{0.35cm}{0.35cm}}
   \put(11.5,-24.65){\rule[8mm]{0.35cm}{0.35cm}}
   \put(12,-24.15){\rule[8mm]{0.35cm}{0.35cm}}
   \put(14,-28.15){\rule[8mm]{0.35cm}{0.35cm}}
   \put(14.5,-27.65){\rule[8mm]{0.35cm}{0.35cm}}
   \put(15,-27.15){\rule[8mm]{0.35cm}{0.35cm}}
   \put(11,-28.15){\rule[8mm]{0.35cm}{0.35cm}}
   \put(11.5,-27.65){\rule[8mm]{0.35cm}{0.35cm}}
   \put(12,-27.15){\rule[8mm]{0.35cm}{0.35cm}}
   \put(13.85,-29.5){$\mathbb{E}^{\bormxi}_8$}
     
\end{picture}
\bigskip\bigskip\bigskip\bigskip\bigskip\bigskip\bigskip\bigskip\bigskip\bigskip\bigskip\bigskip\bigskip

\begin{thm} \label{pixidencl} Let $\xi\!\geq\! 1$ be a countable ordinal, $\mathbb{K},\mathbb{C}$ as above, $X$ be an analytic space, and $E$ be a Borel equivalence relation on $X$ with countably many classes. Then exactly one of the following holds:\smallskip  

(a) the relation $E$ is a $\bormxi$ subset of $X^2$,\smallskip  

(b) there is $(\mathbb{X},\mathbb{E})\!\in\!\{ (\mathbb{K},\mathbb{E}^{\bormxi}_1),
(\mathbb{H},\mathbb{E}^{\bormxi}_8)\}$ such that $(\mathbb{X},\mathbb{E})\sqsubseteq_c(X,E)$.\smallskip

\noindent Moreover, $\{ (\mathbb{K},\mathbb{E}^{\bormxi}_1),
(\mathbb{H},\mathbb{E}^{\bormxi}_8)\}$ is a $\leq_c$-antichain (and thus a 
$\sqsubseteq_c$ and a $\leq_c$-antichain basis).\end{thm}

 We now say a few words about some of the methods used in this paper, and state some general results interesting for themselves. First, we make a strong use of the representation theorem for Borel sets due to Debs and Saint Raymond. In particular, it provides the sequence $(\mathbb{C}_n)_{n\in\omega}$ mentioned before. Secondly, our work is partly based on the Louveau-Saint Raymond theorem (see page 433 in [Lo-SR]) generalizing the Hurewicz theorem (see 21.22 in [K1]). 
 
\begin{thm} \label{Lo-SR} (Louveau-Saint Raymond) Let $\xi\!\geq\! 1$ be a countable ordinal, $\mathbb{K}$ as above, 
$\mathbb{C}\!\in\!\bormxi (\mathbb{K})$ not in $\boraxi$ (as above if $\xi\!\leq\! 2$), $X$ be a Polish space, and $A,B$ be disjoint analytic subsets of $X$. Then exactly one of the following holds:\smallskip

(a) the set $A$ is separable from $B$ by a $\boraxi$ set,\smallskip

(b) we can find $f\! :\!\mathbb{K}\!\rightarrow\! X$ injective continuous such that 
$\mathbb{C}\!\subseteq\! f^{-1}(A)$ and $\neg\mathbb{C}\!\subseteq\! f^{-1}(B)$.\end{thm}

 We will prove and use the following extension of Theorem \ref{Lo-SR}. 

\begin{thm} \label{coord} Let ${\bf\Gamma}$ be a non self-dual Borel class, $\mathbb{K}$ as above, 
$\mathbb{C}\!\in\!\check {\bf\Gamma}(\mathbb{K})$ not in $\bf\Gamma$ (as above if the rank of 
$\bf\Gamma$ is at most two), $X$ be an analytic space, and $A,B$ be disjoint analytic relations on $X$, the sections of $A$ being in $\bf\Gamma$. Then exactly one of the following holds:\smallskip  

(a)  the set $A$ is separable from $B$ by a $\bf\Gamma$ set,\smallskip  

(b) we can find $f\! :\!\mathbb{K}\!\rightarrow\! X^2$ continuous with injective coordinate functions such that 
$\mathbb{C}\!\subseteq\! f^{-1}(A)$ and $\neg\mathbb{C}\!\subseteq\! f^{-1}(B)$.\end{thm}

 The proofs of Theorems \ref{c1dencls} and \ref{pixidencl} use our following other extension of Theorem \ref{Lo-SR}, which provides more control than Theorem \ref{Lo-SR} on where the $\borxi$ sets decomposing $\neg\mathbb{C}$ are sent by $f$. 

\begin{thm} \label{Hur+} Let $\xi\!\geq\! 1$ be a countable ordinal, $\mathbb{K}$ as above, and $\mathbb{C}\!\in\!\bormxi (\mathbb{K})$ (as above if 
$\xi\!\leq\! 2$). Then we can find $\mathbb{I}\!\subseteq\!\omega$ and a partition $(\mathbb{C}_n)_{n\in\mathbb{I}}$ of $\neg\mathbb{C}$ into 
$\borxi$ subsets of $\mathbb{K}$ such that for any analytic space $X$, for any analytic subset $A$ of $X$, and for any sequence 
$(D_n)_{n\in\omega}$ of pairwise disjoint analytic subsets of $X$ such that $A$ is both disjoint from $\bigcup_{n\in\omega}~D_n$ and separable from any of the $D_n$'s by a $\boraxi$ set, one of the following holds:\smallskip

(a)  the set $A$ is separable from $\bigcup_{n\in\omega}~D_n$ by a $\boraxi$ set,\smallskip

(b) we can find $\phi\! :\!\mathbb{I}\!\rightarrow\!\omega$ and $f\! :\!\mathbb{K}\!\rightarrow\! X$ injective continuous such that 
$\mathbb{C}\subseteq\! f^{-1}(A)$ and $\mathbb{C}_n\!\subseteq\! f^{-1}(D_{\phi (n)})$ for each $n\!\in\!\mathbb{I}$. Moreover, we can ensure that if 
$(F_n)_{n\in\omega}$ is a sequence of finite subsets of $\omega$, then $\phi (n)\!\notin\! F_{\phi (p)}$ whenever $p\! <\! n$ are in $\mathbb{I}$. In particular, 
$\phi$ can be injective.\smallskip

 If moreover $\mathbb{C}\!\notin\!\boraxi$, then this is a dichotomy.\end{thm}
  
 The organization of the paper is as follows. In Section 2, we first recall the material about representation of Borel sets and effective topologies needed here. Then we prove our extensions of the Louveau-Saint Raymond Theorem. In Section 3, we prove the main lemma used in the proof of our main result. Essentially, it provides a Cantor set $C$ as in the Mycielski-Kuratowski theorem (see 19.1 in [K1]), with the additional property that the set $\mathbb{C}$ as above remains complex on $C$. In Section 4, we prove some general facts giving additional motivation for introducing our examples, and prove that ${\cal A}^{\bf\Gamma}$ and ${\cal B}^{\bf\Gamma}$ are $\leq_c$-antichains. In Sections 5 and 6, we solve our main questions for the Borel classes of rank one. In Sections 7, 8 and 9, we solve our main questions for the Borel classes of rank two. In Section 10, we prove our results about relations with countably many equivalence classes. In Section 11, we prove our main result.

\section{$\!\!\!\!\!\!$ Extensions of the Louveau-Saint Raymond theorem}\indent

 As in [L2] and [L3], the main results in this section are based on the representation theorem for Borel sets due to Debs and Saint Raymond. We first recall the material related to that needed here.

\subsection{$\!\!\!\!\!\!$ Representation of Borel sets}\indent

 The following definition can be found in [D-SR].

\begin{defin} (Debs-Saint Raymond) A partial order relation $R$ on $2^{<\omega}$ is a {\bf tree relation} if, for $s\!\in\! 2^{<\omega}$,\smallskip

(a) $\emptyset ~R~s$,\smallskip

(b) the set $P_{R}(s)\! :=\!\{ t\!\in\! 2^{<\omega}\mid t~R~s\}$ is finite and linearly ordered by $R$ ($h_R(s)$ will denote the number of strict $R$-predecessors of $s$, so that $h_R(s)\! =\!\mbox{Card}\big( P_R(s)\big)\! -\! 1$).\smallskip

\noindent $\bullet$ Let $R$ be a tree relation. An $R$-{\bf branch} is a $\subseteq$-maximal subset of $2^{<\omega}$ linearly ordered by $R$. We denote by $[R]$ the set of all {\rm infinite} $R$-branches.\smallskip

 We equip $(2^{<\omega})^\omega$ with the product of the discrete topology on $2^{<\omega}$. If $R$ is a tree relation, then the space 
$[R]\!\subseteq\! (2^{<\omega})^\omega$ is equipped with the topology induced by that of $(2^{<\omega})^\omega$, and is a Polish space. A basic clopen set is of the form $N^R_s\! :=\!\big\{\gamma\!\in\! [R]\mid\gamma\big(h_R(s)\big)\! =\! s\big\}$, where $s\!\in\! 2^{<\omega}$.

\vfill\eject

\noindent $\bullet$ Let $R$, $S$ be tree relations with $R\!\subseteq\! S$. The 
{\bf canonical map} $\Pi\! :\! [R]\!\rightarrow\! [S]$ is defined by
$$\Pi (\gamma )\! :=\!\mbox{ the unique }S\mbox{-branch containing }\gamma .$$
The canonical map is continuous.\smallskip

\noindent $\bullet$ Let $S$ be a tree relation. We say that $R\!\subseteq\! S$ is {\bf distinguished} in $S$ if
$$\forall s,t,u\!\in\! 2^{<\omega}\ \ \left.
\begin{array}{ll}
& s~S~t~S~u\cr & \cr
& \ \ s~R~u
\end{array}\right\} ~\Rightarrow ~s~R~t.$$
$\bullet$ Let $\eta\! <\!\omega_1$. A family $(R^\rho )_{\rho\leq\eta}$ of tree relations is a {\bf resolution family} if\smallskip

(a) $R^{\rho +1}$ is a distinguished subtree of $R^\rho$, for each $\rho\! <\!\eta$.\smallskip

(b) $R^\lambda\! =\!\bigcap_{\rho <\lambda}~R^\rho$, for each limit ordinal $\lambda\!\leq\!\eta$.
\end{defin}

 The representation theorem of Borel sets is as follows in the successor case (see Theorems I-6.6 and I-3.8 in [D-SR]).

\begin{them} \label{rep} (Debs-Saint Raymond) Let $\eta$ be a countable ordinal, and $P\!\in\! {\bf\Pi}^0_{\eta +1}([\subseteq ])$. Then there is a resolution family $(R^\rho )_{\rho\leq\eta}$ such that\smallskip

(a) $R^0\! =\subseteq$,\smallskip

(b) the canonical map $\Pi\! :\! [R^\eta ]\!\rightarrow\! [R^0]$ is a continuous bijection with $\boraep$-measurable inverse,\smallskip

(c) the set $\Pi^{-1}(P)$ is a closed subset of $[R^\eta ]$.\end{them}

 For the limit case, we need some more definition that can be found in [D-SR].

\begin{defin} (Debs-Saint Raymond) Let $\xi$ be an infinite limit countable ordinal. We say that a resolution family $(R^\rho )_{\rho\leq\xi}$ with $R^0\! =\subseteq$ is {\bf uniform} if
$$\forall k\!\in\!\omega ~~\exists\xi_k\! <\!\xi ~~\forall s,t\!\in\! 2^{<\omega}~~
\Big(\mbox{min}\big( h_{R^\xi}(s),h_{R^\xi}(t)\big)\!\leq\! k\wedge s\ R^{\xi_k}\ t\Big)\Rightarrow s\ R^\xi\ t.$$
We may (and will) assume that $\xi_k\!\geq\! 1$.\end{defin}

 The representation theorem of Borel sets is as follows in the limit case (see Theorems I-6.6 and I-4.1 in [D-SR]).

\begin{them} \label{replim} (Debs-Saint Raymond) Let $\xi$ be an infinite limit countable ordinal, and $P\!\in\!\bormxi ([\subseteq ])$. Then there is a uniform resolution family 
$(R^\rho )_{\rho\leq\xi}$ such that\smallskip

(a) $R^0\! =\subseteq$,\smallskip

(b) the canonical map $\Pi\! :\! [R^\xi ]\!\rightarrow\! [R^0]$ is a continuous bijection with $\boraxi$-measurable inverse,\smallskip

(c) the set $\Pi^{-1}(P)$ is a closed subset of $[R^\xi ]$.\end{them}

 We will use the following extension of the property of distinction (see Lemma 2.3.2 in [L2]):

\begin{lemm} \label{extdist} Let $\eta\! <\!\omega_1$, $(R^\rho )_{\rho\leq\eta}$ be a resolution family, and 
$\rho\! <\!\eta$. Assume that $s,t,u\!\in\! 2^{<\omega}$, $s~R^0~t~R^\rho ~u$ and $s~R^{\rho +1}~u$. Then 
$s~R^{\rho +1}~t$.\end{lemm}

\vfill\eject

\noindent\bf Notation.\rm ~Let $\eta\! <\!\omega_1$, $(R^\rho )_{\rho\leq\eta}$ be a resolution family with $R^0\! =\subseteq$, $s\!\in\! 2^{<\omega}$, and 
$\rho\!\leq\!\eta$. We define
$$s^\rho\! :=\!\left\{\!\!\!\!\!\!\!\!
\begin{array}{ll}
 & \emptyset\mbox{ if }s\! =\!\emptyset\mbox{,}\cr
 & s\vert\mbox{max}\{ l\! <\!\vert s\vert\mid s\vert l~R^\rho ~s\}\mbox{ if }s\!\not=\!\emptyset .
\end{array}
\right.$$ 
The map $h\! :\! 2^\omega\!\rightarrow\! [\subseteq]$, for which $h (\alpha )$ is the strictly $\subseteq$-increasing sequence of initial segments of $\alpha$, is a homeomorphism.

\subsection{$\!\!\!\!\!\!$ Topologies}

\bf Notation.\rm\ Let $S$ be a recursively presented Polish space.\bigskip

\noindent (1) The {\bf Gandy-Harrington topology} on $S$ is generated by $\Ana (S)$ and denoted 
${\it\Sigma}_S$. Recall the following facts about ${\it\Sigma}_S$ (see [L1]).\smallskip

\noindent - ${\it\Sigma}_S$ is finer than the initial topology of $S$.\smallskip

\noindent - We set $\Omega_S :=\{ s\!\in\! S\mid\omega_1^s\! =\!\omega_1^{\mbox{CK}}\}$. Then 
$\Omega_S$ is $\Ana (S)$ and dense in $(S,{\it\Sigma}_S)$.\smallskip

\noindent - $W\cap\Omega_S$ is a clopen subset of $(\Omega_S,{\it\Sigma}_S)$ for each 
$W\!\in\!\Ana (S)$.\smallskip

\noindent - $(\Omega_S,{\it\Sigma}_S)$ is a zero-dimensional Polish space. So we fix a complete compatible metric on $(\Omega_S,{\it\Sigma}_S)$.\bigskip

\noindent (2) We call $T_1$ the usual topology on $S$, and $T_\eta$ is the topology generated by the 
$\Ana\cap {\bf\Pi}^0_{<\eta}$ subsets of $S$ if $2\!\leq\!\eta\! <\!\omega_1^{\mbox{CK}}$ (see Definition 1.5 in [Lo1]).\bigskip

 The next result is essentially Lemma 2.2.2 and the claim in the proof of Theorem 2.4.1 in [L2].

\begin{lemm} \label{top} Let $S$ be a recursively presented Polish space, and 
$1\!\leq\!\eta\! <\!\omega^{\mbox{CK}}_1$.\smallskip

(a) (Louveau) Fix $A\!\in\!\Ana (S)$. Then $\overline{A}^{T_\eta}$ is $\borme$ and $\Ana$.\smallskip

(b) (Louveau) Fix $A,D\!\in\!\Ana (S)$ disjoint. Then $D$ is separable from $A$ by a $\borme$ set exactly when 
$A\cap\overline{D}^{T_\eta}\! =\!\emptyset$.\smallskip

(c) Let $p\!\geq\! 1$ be a natural number, $1\!\leq\!\eta_{1}\! <\!\eta_{2}\! <\!\ldots\! <\!\eta_{p}\!\leq\!\eta$, $S_1$, 
$\ldots$, $S_p\!\in\!\Ana (S)$, and $O\!\in\!\Boraone (S)$. Assume that 
$S_i\!\subseteq\!\overline{S_{i+1}}^{T_{\eta_i+1}}$ if $1\!\leq\! i\! <\! p$. Then 
${S_p\cap\bigcap_{1\leq i<p}~\overline{S_i}^{T_{\eta_i}}}\cap O$ is $T_1$-dense in $\overline{S_1}^{T_1}\cap O$.\smallskip

(d) Let $(R^\rho )_{\rho\leq\eta}$ be a resolution family with $R^0\! =\subseteq$, 
$s\!\in\! 2^{<\omega}\!\setminus\!\{\emptyset\}$, $S_{s^\rho}\!\in\!\Ana (S)$ (for 
$1\!\leq\!\rho\!\leq\!\eta$),  $E\!\in\!\Ana (S)$, and $O\!\in\!\Boraone (S)$. We assume that 
$S_{s^\eta}\!\subseteq\!\overline{S}^{T_{\eta +1}}$ and 
$S_t\!\subseteq\!\overline{S_u}^{T_\rho}\mbox{ if }u~R^\rho ~t\!\subsetneqq\! s$ and 
$1\!\leq\!\rho\!\leq\!\eta$. Then 
$S_{s^\eta}\cap\bigcap_{1\leq\rho <\eta}~\overline{S_{s^\rho}}^{T_\rho}\cap O$ and 
$E\cap\bigcap_{1\leq\rho\leq\eta}~\overline{S_{s^\rho}}^{T_\rho}\cap O$ are $T_1$-dense in 
$\overline{S_{s^1}}^{T_1}\cap O$.\end{lemm}

\noindent\bf Proof.\rm ~(a) and (b) See Lemmas 1.7 and 1.9 in [Lo1].\bigskip

\noindent (c) and (d) See Lemma 2.2.1 in [L4].\hfill{$\square$}

\begin{lemm} \label{top2} Let $S$ be a recursively presented Polish space.\smallskip

(a) The set $\Borel\cap S$ is countable, $\Ca$, and $T_\eta$-clopen if $3\!\leq\!\eta\! <\!\omega^{\mbox{CK}}_1$.\smallskip

(b) If $A,B$ are disjoint subsets of $S$ and $2\!\leq\!\eta\! <\!\omega^{\mbox{CK}}_1$, then 
$A\cap\overline{B}^{T_\eta}$ does not meet $\Borel\cap S$.\end{lemm}

\noindent\bf Proof.\rm ~(a) By 4D.2 and 4D.14 in [M], $\Borel\cap S$ is countable and $\Ca$, so that its complement is $\Ana\cap\bormtwo$ and thus $T_\eta$-open if $\eta\!\geq\! 3$. Moreover, $\Borel\cap S$ is the union of its singletons, which are closed and $\Borel$ and thus $T_\eta$-open if $\eta\!\geq\! 2$. This shows that 
$\Borel\cap S$ is $T_\eta$-open if $\eta\!\geq\! 2$.\bigskip

\noindent (b) We argue by contradiction, which gives $s$ in the intersection of $A\cap\overline{B}^{T_\eta}$ and 
$\Borel\cap S$. By (a), $\{ s\}$ is $T_\eta$-open, so that $s\!\in\! B\cap A$, contradicting the disjointness of $A$ and $B$.\hfill{$\square$}

\subsection{$\!\!\!\!\!\!$ Proof of Theorem \ref{Hur+}}

\label{SectionHur+}

\bf (A) The successor case\rm\bigskip

 Assume that $\xi\! =\!\eta\!+\! 1\!\geq\! 3$. As $h$ is a homeomorphism, $P\! :=\! h[\mathbb{C}]\!\in\!\bormep ([\subseteq ])$. Theorem \ref{rep} gives a resolution family $(R^\rho )_{\rho\leq\eta}$ such that ${\cal C}\! :=\!\Pi^{-1}(P)$ is a closed subset of $[R^\eta ]$. If $\gamma\!\in\! [R^\eta ]\!\setminus\! {\cal C}$, then there is $k_\gamma\!\in\!\omega$ minimal such that $N^{R^\eta}_{\gamma (k_\gamma )}\cap {\cal C}\! =\!\emptyset$. We set 
$s_\gamma\! :=\!\gamma (k_\gamma )$, so that $\neg {\cal C}\! =\!\bigcup_{\gamma\notin {\cal C}}~N^{R^\eta}_{s_\gamma}$. We enumerate injectively 
${\cal S}\! :=\!\{ s_\gamma\mid\gamma\!\notin\! {\cal C}\}$, which gives $\mathbb{I}\!\subseteq\!\omega$ such that ${\cal S}\! =\!\{ s_n\mid n\!\in\!\mathbb{I}\}$. As ${\cal S}\!\subseteq\! 2^{<\omega}$, we can ensure that the enumaration is made in the increasing order of the lengths, and in the lexicographical order inside each length. We set ${\cal C}_n\! :=\! N^{R^\eta}_{s_n}$, so that $\neg {\cal C}\! =\!\bigcup_{n\in\mathbb{I}}~{\cal C}_n$. By minimality, this union is disjoint. We then set $\mathbb{C}_n\! :=\! h^{-1}(\Pi[{\cal C}_n])$, so that $(\mathbb{C}_n)_{n\in\mathbb{I}}$ is a partition of $\neg\mathbb{C}$ into $\borxi$ subsets of $2^\omega$. Note that 
$\mbox{diam}(\mathbb{C}_n)$ tends to $0$ as $n$ goes to infinity if $\mathbb{I}$ is infinite since 
$\mathbb{C}_n\!\subseteq\! N_{s_n}$.\bigskip

 Assume that (a) does not hold. As $X$ is analytic, we may assume that $X$ is an analytic subset of $S\! :=\! [0,1]^\omega$, as well as $A$ and the $D_n$'s. As our separation assumptions also hold in $S$, we may assume that $X\! =\! S$. In order to simplify the notation, we will assume that 
$\xi\! <\!\omega_1^{\mbox{CK}}$, $A\!\in\!\Ana$ and the relation defined by 
$R(\alpha ,n)\Leftrightarrow\alpha\!\in\! D_n$ is $\Ana$. As $\bigcup_{n\in\omega}~D_n$ is not separable from $A$ by a $\bormxi$ set, 
$N\! :=\! A\cap\overline{\bigcup_{n\in\omega}~D_n}^{T_\xi}$ is a nonempty $\Ana$ subset of $S$, by Lemma \ref{top}.(a). We set $U\! :=\!\Borel\cap S$. By Lemma \ref{top2}.(b), $N\cap U\! =\!\emptyset$. By Lemma \ref{top2}.(a), $U$ is $T_\xi$-clopen since $\xi\!\geq\! 3$. This shows that 
$N\! =\! A\cap\overline{\bigcup_{n\in\omega}~D_n}^{T_\xi}\!\setminus\! U\! =\! 
A\cap\overline{\bigcup_{n\in\omega}~D_n\!\setminus\! U}^{T_\xi}$. By Lemma \ref{top}.(b), 
$A\cap\overline{D_n}^{T_\xi}\! =\!\emptyset$ for each $n\!\in\!\omega$, so that 
$N\! =\! A\cap\overline{\bigcup_{n\in\omega\setminus F}~D_n\!\setminus\! U}^{T_\xi}$ for each finite subset $F$ of $\omega$.\bigskip
  
 We set, for $n\!\in\!\mathbb{I}$, ${\cal O}_n\! :=\!\{ s\!\in\! 2^{<\omega}\mid s_n~R^\eta ~s\}$, and 
${\cal I}\! :=\! 2^{<\omega}\!\setminus(\bigcup_{n\in\mathbb{I}}~{\cal O}_n)$. By definition of the $s_n$'s, 
$2^{<\omega}$ is the disjoint union of $\cal I$ and the ${\cal O}_n$'s. If $\emptyset\!\notin\! {\cal I}$, then there is 
$n\!\in\!\mathbb{I}$ such that $\emptyset\!\in\! {\cal O}_n$, so that $s_n\! =\!\emptyset$, ${\cal O}_n\! =\! 2^{<\omega}$, 
$\mathbb{I}\! =\!\{ n\}$, $\mathbb{C}_n\! =\! 2^\omega$ and $\mathbb{C}\! =\!\emptyset$. There is $p\!\in\!\omega$ such that $D_p$ is uncountable since 
$\xi\!\geq\! 3$, so that we can set $\phi (n)\! :=\! p$ by the perfect set theorem. So in the sequel we will assume that $\emptyset\!\in\! {\cal I}$. We construct\bigskip

\noindent - a sequence $(x_s)_{s\in 2^{<\omega}}$ of points of $S$,\smallskip

\noindent - a sequence $(X_s)_{s\in 2^{<\omega}}$ of $\Boraone$ subsets of $S$,\smallskip

\noindent - a sequence $(S_s)_{s\in 2^{<\omega}}$ of $\Ana$ subsets of $S$,\smallskip

\noindent - $\phi\! :\!\mathbb{I}\!\rightarrow\!\omega$,\bigskip

\noindent satisfying the following conditions.

\vfill\eject

$$\begin{array}{ll}
& (1)~\left\{\!\!\!\!\!\!\!\!
\begin{array}{ll}
& \overline{X_t}\!\subseteq\! X_s\mbox{ if }~s~R^0~t\wedge s\!\not=\! t\cr
& S_t\!\subseteq\! S_s\mbox{ if }~s~R^\eta ~t\wedge (s,t\!\in\! {\cal I}\vee
\exists n\!\in\!  I~~s,t\!\in\! {\cal O}_n)
\end{array}
\right.\cr\cr
& (2)~x_s\!\in\! S_s\!\subseteq\! X_s\cap\Omega_S\!\setminus\! U\cr\cr
& (3)~\mbox{diam}(X_s),\mbox{diam}_{\mbox{GH}}(S_s)\!\leq\! 2^{-\vert s\vert}\cr\cr
& (4)~X_{s0}\cap X_{s1}\! =\!\emptyset\cr\cr
& (5)~S_s\!\subseteq\!\left\{\!\!\!\!\!\!\!\!
\begin{array}{ll}
& N\mbox{ if }~s\!\in\! {\cal I}\cr
& D_{\phi (n)}\mbox{ if }s\!\in\! {\cal O}_n
\end{array}
\right.\cr\cr
& (6)~S_t\!\subseteq\!\overline{S_s}^{T_\rho}\mbox{ if }s~R^\rho ~t\wedge 1\!\leq\!\rho\!\leq\!\eta
\end{array}$$
Assume that this is done. Let $\alpha\!\in\! 2^\omega$. Then 
$(\overline{X_{\alpha\vert l}})_{l\in\omega}$ is a decreasing sequence of nonempty closed subsets of $S$ with vanishing diameters, which defines $f\! :\! 2^\omega\!\rightarrow\! S$ injective continuous. If $\alpha\!\in\!\mathbb{C}$, then 
$\Pi^{-1}\big( h(\alpha )\big)(k)\!\in\! {\cal I}$ for each $k\!\in\!\omega$. Note that 
$(S_{\Pi^{-1}(h(\alpha ))(k)})_{k\in\omega}$ is a decreasing sequence of nonempty clopen subsets of 
$N\cap\Omega_S\!\subseteq\! A$ with vanishing GH-diameters, which defines $G(\alpha )\!\in\! A$. As 
$S_s\!\subseteq\! X_s$, $G(\alpha )\! =\! f(\alpha )$, so that $f(\alpha )\!\in\! A$. If now $\alpha\!\in\!\mathbb{C}_n$, then 
$\Pi^{-1}\big( h(\alpha )\big)\!\in\! {\cal C}_n$ and $\Pi^{-1}\big( h(\alpha )\big)(k)$ is in ${\cal O}_n$ if 
$k\!\geq\! k_0$. Note that $(S_{\Pi^{-1}(h(\alpha ))(k)})_{k\geq k_0}$ is a decreasing sequence of nonempty clopen subsets of $D_{\phi (n)}\cap\Omega_S$ with vanishing GH-diameters, which defines $H(\alpha )\!\in\! D_{\phi (n)}$. As $S_s\!\subseteq\! X_s$, $H(\alpha )\! =\! f(\alpha )$, so that $f(\alpha )\!\in\! D_{\phi (n)}$.\bigskip

 Let us prove that the construction is possible. We first choose 
$x_\emptyset\!\in\! N\cap\Omega_S$, $X_\emptyset$ semi-recursive with diameter at most $1$ containing $x_\emptyset$, and $S_\emptyset\!\in\!\Ana$ with GH-diameter at most $1$ containing $x_\emptyset$ and contained in $X_\emptyset\cap N\cap\Omega_S$. Assume that our objects satisfying (1)-(6) have been contructed up to the length $l$ (which means that $\phi (n)$ is constructed if $\vert s_n\vert\!\leq\! l$), which is the case for $l\! =\! 0$.\bigskip

 Let $t\!\in\! 2^l$, and $s\! :=\! t0$. We first define $x_s$, and $X_s$ and $S_s$ later, after the definition of $x_{t1}$. Our definitions are in the lexicographical order of the $t$'s.\bigskip

\noindent\bf Claim.\it ~(a) The set 
$S_{s^\eta}\cap\bigcap_{1\leq\rho <\eta}~\overline{S_{s^\rho}}^{T_\rho}\cap X_{s^0}$ is $\Ana$ and uncountable.\smallskip

(b) If $s^\eta\!\in\! {\cal I}$ and $F\!\subseteq\!\omega$ is finite, then 
$(\bigcup_{n\in\omega\setminus F}~D_n\!\setminus\! U)\cap
\bigcap_{1\leq\rho\leq\eta}~\overline{S_{s^\rho}}^{T_\rho}\cap X_{s^0}\cap\Omega_S$ is $\Ana$ and uncountable.\bigskip\rm  

 Indeed, by Lemma \ref{top}.(d) applied to $E\! :=\! \bigcup_{n\in\omega\setminus F}~D_n\!\setminus\! U$ and 
$O\! :=\! X_{s^0}$, the sets 
$$S_{s^\eta}\cap\bigcap_{1\leq\rho <\eta}~\overline{S_{s^\rho}}^{T_\rho}\cap X_{s^0}$$ 
and $(\bigcup_{n\in\omega\setminus F}~D_n\!\setminus\! U)\cap
\bigcap_{1\leq\rho\leq\eta}~\overline{S_{s^\rho}}^{T_\rho}\cap X_{s^0}$ are $T_1$-dense in 
$\overline{S_{s^1}}^{T_1}\cap X_{s^0}$. As $s^1~R^1~s^0$, 
$$S_{s^0}\!\subseteq\!\overline{S_{s^1}}^{T_1}\cap X_{s^0}.$$ 

 This proves that the intersections in the statement are not empty since they are $\Ana$ by Lemma \ref{top}.(a). We argue by contradiction to see that they are uncountable. By 4F.1 in [Mos], they are contained in $U$, which contradicts the induction assumption.\hfill{$\diamond$}\bigskip

\noindent\bf Case 1\rm\ $s\!\in\! {\cal I}$, which implies that $s^\eta\!\in\! {\cal I}$.\bigskip

 We choose 
$x_s\!\in\! S_{s^\eta}\cap\bigcap_{1\leq\rho <\eta}~\overline{S_{s^\rho}}^{T_\rho}\cap X_{s^0}$.\bigskip
 
\noindent\bf Case 2\rm\ $s^\eta\!\in\! {\cal O}_n$, which implies that $s\!\in\! {\cal O}_n$.\bigskip

 We proceed as in Case 1.\bigskip

\noindent\bf Case 3\rm\ $s^\eta\!\in\! {\cal I}$ and $s\!\in\! {\cal O}_n$.\bigskip

 In this case, $s\! =\! s_n$ and $\phi (n)$ has to be defined. We choose it outside 
$F\! :=\!\bigcup_{p<n,p\in\mathbb{I}}~F_{\phi (p)}$ in such a way that 
$(D_{\phi (n)}\!\setminus\! U)\cap\bigcap_{1\leq\rho\leq\eta}~
\overline{S_{s^\rho}}^{T_\rho}\cap X_{s^0}\cap\Omega_S$ is uncountable. Then we choose $x_s$ in $(D_{\phi (n)}\!\setminus\! U)\cap\bigcap_{1\leq\rho\leq\eta}~
\overline{S_{s^\rho}}^{T_\rho}\cap X_{s^0}\cap\Omega_S$.\bigskip

 This finishes the construction of $x_{t0}$, which is in the right uncountable $\Ana$ set. The construction of $x_{t1}$ is similar, the difference being that we ensure moreover that 
$x_{t1}\!\not=\! x_{t0}$, which is possible since the right $\Ana$ set is uncountable. Then we choose disjoint $\Boraone$ sets $X_{t0}$ and $X_{t1}$ with diameter at most $2^{-l-1}$ such that $x_{t\varepsilon}\!\in\! X_{t\varepsilon}\!\subseteq\!\overline{X_{t\varepsilon}}\!\subseteq\! X_t$, and $S_{t\varepsilon}\!\in\!\Ana$ with GH-diameter at most $2^{-l-1}$ containing $x_{t\varepsilon}$ and contained in $X_{t\varepsilon}$ and the right $\Ana$ set. Note that we ensured (1) and (6) for the immediate predecessors of $t$, and not for an arbitrary $s$ at this point. These conditions are transitive enough to ensure the general case. For example, for (6), assume that $s~R^\rho ~t$. We may assume that $s\!\not=\! t$, so that $s~R^0~s^\rho ~R^\rho ~t$. By Lemma \ref{extdist}, $s~R^\rho ~s^\rho$. By induction assumption, 
$S_{s^\rho}\!\subseteq\!\overline{S_s}^{T_\rho}$. We ensured that 
$S_t\!\subseteq\!\overline{S_{s^\rho}}^{T_\rho}$, so that $S_t\!\subseteq\!\overline{S_s}^{T_\rho}$.
\hfill{$\square$}\bigskip

 We now study the case $\xi\! =\! 2$.\bigskip
 
\noindent\bf Notation.\rm\ In this case, $\mathbb{C}\! =\!\mathbb{P}_\infty\! :=\!
\{\alpha\!\in\! 2^\omega\mid\exists^\infty n\!\in\!\omega ~~\alpha (n)\! =\! 1\}$.  By 23.A in [K1], 
$\mathbb{P}_\infty\!\in\!\bormtwo (2^\omega )\!\setminus\!\boratwo$. Note that $\mathbb{P}_\infty$ is dense and co-dense in $2^\omega$. We set $\mathbb{P}_f\! :=\!\neg\mathbb{P}_\infty$ and enumerate $\mathbb{P}_f\! :=\!\{\alpha_n\mid n\!\in\!\omega\}$. We also set 
$O\! :=\!\{\emptyset\}\cup\{ u1\mid u\!\in\! 2^{<\omega}\}$, so that 
$\mathbb{P}_f\! =\!\{ t0^\infty\mid t\!\in\! O\}$.\bigskip

 We set $\mathbb{I}\! :=\!\omega$ and $\mathbb{C}_n\! :=\!\{\alpha_n\}$, which defines a partition of $\neg\mathbb{C}$ into $\bortwo$ subsets of $2^\omega$. We also set $R^0\! :=\subseteq$, and  
$$s~R^1~t\Leftrightarrow s~R^0~t\wedge\big( s\!\in\! O\vee\forall s~R^0~u~R^0~t~~u\!\notin\! O\big) .$$
Note that $(R^\rho )_{\rho\leq 1}$ is a resolution family such that\bigskip

(b) the canonical map $\Pi\! :\! [R^1]\!\rightarrow\! [R^0]$ is a continuous bijection with $\boratwo$-measurable inverse,\smallskip

(c) the sets ${\cal C}_n\! :=\!\Pi^{-1}(h[\mathbb{C}_n])$ are clopen subsets of $[R^1]$, so that 
${\cal C}\! :=\!\Pi^{-1}(h[\mathbb{C}])$ is a closed subset of $[R^1]$.\bigskip

 In fact, as $\mathbb{C}_n$ is a singleton, ${\cal C}_n$ too, and ${\cal C}_n\! =\! N^{R^1}_{s_n}$ for some $s_n\!\in\! 2^{<\omega}$ of the form $t_n0$ with $t_n\!\in\! O$ (so that $\mathbb{C}_n\! =\!\{ t_n0^\infty\}$), and the ${\cal C}_n$'s define a partition of $\neg {\cal C}$ as above.\bigskip

 We then argue as in the case $\xi\!\geq\! 3$, with the following differences. This time, we only write 
$$N\! :=\! A\cap\overline{\bigcup_{n\in\omega}~D_n}^{T_2}\! =\! 
A\cap\overline{\bigcup_{n\in\omega\setminus F}~D_n}^{T_2}$$ 
for each finite subset $F$ of $\omega$. Note that ${\cal O}_n\cap 2^l$ has cardinality at most $1$ for each 
$l\!\in\!\omega$. Condition (2) becomes
$$(2')~x_s\!\in\! S_s\!\subseteq\! X_s\cap\Omega_S$$
\bf Claim.\it ~(a) $S_{s^1}\cap X_{s^0}$ is $\Ana$, nonempty, and uncountable if $s^1\!\in\! O$.\smallskip

(b) If $s^1\!\in\! O$ and $F$ is a finite subset of $\omega$, then the set 
$\big(\bigcup_{n\in\omega\setminus F}~D_n\big)\cap
\overline{S_{s^1}}^{T_1}\cap X_{s^0}\cap\Omega_S$ is $\Ana$ and nonempty.\bigskip\rm  

 In Case 3, we choose $\phi (n)$ in such a way that 
$D_{\phi (n)}\cap\overline{S_{s^1}}^{T_1}\cap X_{s^0}\cap\Omega_S$ is nonempty. Then we choose $x_s\!\in\! D_{\phi (n)}\cap\overline{S_{s^1}}^{T_1}\cap X_{s^0}\cap\Omega_S$.\bigskip

 This finishes the construction of $x_{t0}$, which is in the right $\Ana$ set. The construction of 
$x_{t1}$ is similar. Note that $x_{t1}\!\not=\! x_{t0}$ since $t0\!\notin\! O$ and $t1\!\in\! O$, so that $x_{t0}\!\in\!\bigcup_{n\in\omega}~D_n$ and $x_{t1}\!\in\! A$.\bigskip

 Assume finally that $\xi\! =\! 1$, so that $\mathbb{C}\! =\!\{ 0\}$. We set $\mathbb{I}\! :=\!\omega$ and $\mathbb{C}_n\! :=\!\{ 2^{-n}\}$, which defines a partition of $\neg\mathbb{C}$ into clopen subsets of $\mathbb{K}$. Fix $x\!\in\! N$, $\phi (0)\!\in\!\omega$ with $B(x,1)\cap D_{\phi (0)}\!\not=\!\emptyset$, and 
$x_0\!\in\! B(x,1)\cap D_{\phi (0)}$. As $D_n$ is separable from $A$ by a closed set for each $n$, $A\cap\overline{D_n}^{T_1}\! =\!\emptyset$ and 
$N\! =\! A\cap\overline{\bigcup_{n\in\omega\setminus F}~D_n}^{T_1}$ for each finite subset $F$ of 
$\omega$. So we can choose $\phi (1)\!\in\!\omega\!\setminus\!\big( F_{\phi (0)}\cup\{\phi (0)\}\big)$ with 
$B(x,2^{-1})\cap D_{\phi (1)}\!\not=\!\emptyset$, and $x_1\!\in\! B(x,2^{-1})\cap D_{\phi (1)}$. And so on. It remains to set $f(0)\! :=\! x$ and $f(2^{-n})\! :=\! x_n$.\bigskip

\noindent\bf (B) The limit case\rm\bigskip

 Assume that $\xi$ is an infinite limit ordinal. We indicate the differences with the successor case. Theorem \ref{replim} gives a uniform resolution family 
$(R^\rho )_{\rho\leq\xi}$ such that $\cal C$ is a closed subset of $[R^\xi ]$. This time, 
${\cal O}_n\! :=\!\{ s\!\in\! 2^{<\omega}\mid s_n~R^\xi ~s\}$. If $s\!\in\! 2^{<\omega}$, then we set, as in the proof of Theorem 2.4.4 in [L2],
$$\xi (s)\! :=\!\mbox{max}\{\xi_{h_{R^\xi}(t)+1}\mid t\!\subseteq\! s\} .$$ 
Note that $\xi (t)\!\leq\!\xi (s)$ if $t\!\subseteq\! s$.\bigskip

 Conditions (1) and (6) are changed as follows:
$$\begin{array}{ll}
& (1')~\left\{\!\!\!\!\!\!\!\!
\begin{array}{ll}
& \overline{X_t}\!\subseteq\! X_s\mbox{ if }~s~R^0~t\wedge s\!\not=\! t\cr
& S_t\!\subseteq\! S_s\mbox{ if }~s~R^\xi ~t\wedge (s,t\!\in\! {\cal I}\vee
\exists n\!\in\! I~~s,t\!\in\! {\cal O}_n)
\end{array}
\right.\cr\cr
& (6')~S_t\!\subseteq\!\overline{S_s}^{T_\rho}\mbox{ if }s~R^\rho ~t\wedge 1\!\leq\!\rho\!\leq\!\xi (s)
\end{array}$$

 The next claim and the remark after it were already present in the proof of Theorem 2.4.4 in [L2].\bigskip

\noindent\bf Claim 1\rm ~Assume that $s^\rho\!\not=\! s^\xi$. Then $\rho\! +\! 1\!\leq\!\xi (s^{\rho +1})$.\bigskip 

 We argue by contradiction. We get
$$\rho\! +\! 1\! >\!\rho\!\geq\!\xi (s^{\rho +1})\!\geq\!\xi_{h_{R^\xi}(s^{\xi})+1}\! =\!\xi_{h_{R^\xi}(s)}.$$
As $s^\rho\ R^\rho\ s$, $s^\rho\ R^\xi\ s$ and $s^\rho\! =\! s^\xi$, which is absurd.\hfill{$\diamond$}\bigskip

 Note that $\xi_{n-1}\! <\!\xi_{n-1}\! +\! 1\!\leq\!\xi (s^{\xi_{n-1}+1})\!\leq\!\xi (s)$. Thus 
$s^{\xi (s)}\! =\! s^{\xi}$.\bigskip

\noindent\bf Claim 2\rm ~(a) The set 
$S_{s^\xi}\cap\bigcap_{1\leq\rho <\xi (s)}~\overline{S_{s^\rho}}^{T_\rho}\cap X_{s^0}$ 
is $\Ana$ and uncountable.\smallskip

(b) If $s^\xi\!\in\! {\cal I}$ and $F\!\subseteq\!\omega$ is finite, then 
$(\bigcup_{n\in\omega\setminus F}~D_n\!\setminus\! U)\cap
\bigcap_{1\leq\rho\leq\xi (s)}~
\overline{S_{s^\rho}}^{T_\rho}\cap X_{s^0}\cap\Omega_S$ is $\Ana$ and uncountable.\rm\bigskip

 We conclude as in the successor case, using the facts that $\xi_k\!\geq\! 1$ and $\xi (.)$ is increasing.
\hfill{$\square$}

\subsection{$\!\!\!\!\!\!$ Proof of Theorem \ref{coord}}\indent

 We prove Theorem \ref{coord} for ${\bf\Gamma}\! =\!\boraxi$, the other case being similar. Note that (a) and (b) cannot hold simultaneously. We indicate the differences with the proof of Theorem \ref{Hur+}.\bigskip

\noindent\bf (A) The successor case\rm\bigskip

 Assume that (a) does not hold. As $X$ is analytic, we may assume that $X$ is an analytic subset of $[0,1]^\omega$, and that $A$ and $B$ are analytic subsets of 
$S\! :=\! [0,1]^\omega\!\times\! [0,1]^\omega$. Note that $A$ is not separable from $B$ by a 
$\boraxi$ subset of $S$. In order to simplify the notation, we will assume that 
$\xi\! <\!\omega_1^{\mbox{CK}}$ and $X,A,B\!\in\!\Ana$. As $B$ is not separable from $A$ by a 
$\bormxi$ set, $M\! :=\! A\cap\overline{B}^{T_\xi}$ is a nonempty $\Ana$ subset of $S$, by Lemma \ref{top}.\bigskip

 Let us prove that $M$ is not contained in the $T_2$-open set 
$$U\ :=\!\bigcup_{x\in\Borel\cap [0,1]^\omega}~
(\{ x\}\!\times\! [0,1]^\omega\cup [0,1]^\omega\!\times\!\{ x\} ).$$ 
We argue by contradiction to see that. Note that $A\!\setminus\!\overline{B}^{T_\xi}$ is separable from $B$ by the $\boraxi$ set $\neg\overline{B}^{T_\xi}$. As $A$ is not separable from $B$ by a 
$\boraxi$ subset of $S$, this implies that $M$ is not separable from $B$ by a $\boraxi$ subset of $S$. This gives $x\!\in\!\Borel\cap [0,1]^\omega$ such that, for example, 
$A\cap (\{ x\}\!\times\! [0,1]^\omega )$ is not separable from $B\cap (\{ x\}\!\times\! [0,1]^\omega )$ by a $\boraxi$ subset of $S$ since $\Borel\cap [0,1]^\omega$ is countable (see Lemma \ref{top2}). Therefore $A\cap (\{ x\}\!\times\! [0,1]^\omega )$ is not a $\boraxi$ subset of $X^2$, which means that $x\!\in\! X$ and the vertical section $A_x\! :=\!\{ y\!\in\! X\mid (x,y)\!\in\! A\}$ is not a $\boraxi$ subset of $X$, which is absurd.\bigskip

 Note that any nonempty $\Ana$ subset of $S$ which is disjoint from $U$ has uncountable projections, by 4D.14 in [M]. By 4D.14 in [M], the set $U$ is $\Ca\cap\boratwo$, so that its complement is $\Ana\cap\bormtwo$ and thus $T_\xi$-open since $\xi\!\geq\! 3$. This shows that 
$N\! :=\! A\cap\overline{B}^{T_\xi}\!\setminus\! U\! =\! 
A\cap\overline{B\!\setminus\! U}^{T_\xi}$ is a nonempty $\Ana$ subset of $S$.

\vfill\eject
 
 We set ${\cal I}\! :=\!\{ s\!\in\! 2^{<\omega}\mid N^{R^\eta}_s\cap {\cal C}\!\not=\!\emptyset\}$. As 
$\mathbb{C}\!\not=\!\emptyset$, $\emptyset\!\in\! {\cal I}$. We construct a sequence 
$(p_s)_{s\in 2^{<\omega}}$ of points of $S$ (where $p_s\! :=\! (x_s,y_s)$), a sequence 
$(X_s)_{s\in 2^{<\omega}}$ of $\Boraone$ subsets of $S$, and a sequence 
$(S_s)_{s\in 2^{<\omega}}$ of $\Ana$ subsets of $S$ satisfying the following conditions.
$$\begin{array}{ll}
& (1)~\left\{\!\!\!\!\!\!\!\!
\begin{array}{ll}
& \overline{X_t}\!\subseteq\! X_s\mbox{ if }~s~R^0~t\wedge s\!\not=\! t\cr
& S_t\!\subseteq\! S_s\mbox{ if }~s~R^\eta ~t\wedge (s,t\!\in\! {\cal I}\vee s,t\!\notin\! {\cal I})
\end{array}
\right.\cr\cr
& (2)~p_s\!\in\! S_s\!\subseteq\! X_s\cap\Omega_S\!\setminus\! U\cr\cr
& (3)~\mbox{diam}(X_s),\mbox{diam}_{\mbox{GH}}(S_s)\!\leq\! 2^{-\vert s\vert}\cr\cr
& (4)~\forall\varepsilon\!\in\! 2~~\mbox{proj}_\varepsilon [X_{s0}]\cap\mbox{proj}_\varepsilon [X_{s1}]\! =\!\emptyset\cr\cr
& (5)~S_s\!\subseteq\!\left\{\!\!\!\!\!\!\!\!
\begin{array}{ll}
& N\mbox{ if }~s\!\in\! {\cal I}\cr
& B\mbox{ if }s\!\notin\! {\cal I}
\end{array}
\right.\cr\cr
& (6)~S_t\!\subseteq\!\overline{S_s}^{T_\rho}\mbox{ if }s~R^\rho ~t\wedge 1\!\leq\!\rho\!\leq\!\eta
\end{array}$$
Assume that this is done. Let $\alpha\!\in\! 2^\omega$. Then 
$(\overline{X_{\alpha\vert l}})_{l\in\omega}$ is a decreasing sequence of nonempty closed subsets of $S$ with vanishing diameters, which defines $f\! :\! 2^\omega\!\rightarrow\! S$ continuous with injective coordinates. If $\alpha\!\in\!\mathbb{C}$, then $f(\alpha )\!\in\! A$ again. If now $\alpha\!\notin\!\mathbb{C}$, then $\Pi^{-1}\big( h(\alpha )\big)(k)$ is not in ${\cal I}$ if $k\!\geq\! k_0$. Note that $(S_{\Pi^{-1}(h(\alpha ))(k)})_{k\geq k_0}$ is a decreasing sequence of nonempty clopen subsets of $B\cap\Omega_S$ with vanishing GH-diameters, which defines $H(\alpha )\!\in\! B$. As $S_s\!\subseteq\! X_s$, 
$H(\alpha )\! =\! f(\alpha )$, so that $f(\alpha )\!\in\! B$.\bigskip

 Let us prove that the construction is possible.\bigskip

\noindent\bf Claim.\it ~(a) The set 
$S_{s^\eta}\cap\bigcap_{1\leq\rho <\eta}~\overline{S_{s^\rho}}^{T_\rho}\cap X_{s^0}$ is $\Ana$ and nonempty.\smallskip

(b) If $s^\eta\!\in\! {\cal I}$, then the set $(B\!\setminus\! U)\cap
\bigcap_{1\leq\rho\leq\eta}~\overline{S_{s^\rho}}^{T_\rho}\cap X_{s^0}\cap\Omega_S$ is $\Ana$ and nonempty.\bigskip\rm  

\noindent\bf Case 1\rm\ $s\!\in\! {\cal I}$, which implies that $s^\eta\!\in\! {\cal I}$.\bigskip

 We choose 
$p_s\!\in\! S_{s^\eta}\cap\bigcap_{1\leq\rho <\eta}~\overline{S_{s^\rho}}^{T_\rho}\cap X_{s^0}$.\bigskip
 
\noindent\bf Case 2\rm\ $s^\eta\!\notin\! {\cal I}$, which implies that $s\!\notin\! {\cal I}$.\bigskip

 We proceed as in Case 1.\bigskip

\noindent\bf Case 3\rm\ $s^\eta\!\in\! {\cal I}$ and $s\!\notin\! {\cal I}$.\bigskip

 We choose $p_s\!\in\! (B\!\setminus\! U)\cap
\bigcap_{1\leq\rho\leq\eta}~\overline{S_{s^\rho}}^{T_\rho}\cap X_{s^0}\cap\Omega_S$. This finishes the construction of $p_{t0}$, which is in the right uncountable $\Ana$ set. The construction of $p_{t1}$ is similar, the difference being that we have to ensure that moreover $x_{t1}\!\not=\! x_{t0}$ 
and $y_{t1}\!\not=\! y_{t0}$. We first choose $\tilde p_{t1}$ in the right $\Ana$ set $\cal A$ as above, ensuring that $\tilde x_{t1}\!\not=\! x_{t0}$, which is possible since $\cal A$ is disjoint from $U$ and therefore has an uncountable first projection. We then choose $O\!\in\!\Boraone$ with 
$\tilde x_{t1}\!\in\! O$ and $x_{t0}\!\notin\! O$, so that ${\cal A}\cap (O\!\times\! [0,1]^\omega)$ is again a nonempty $\Ana$ set disjoint from $U$. We now choose $p_{t1}$ in 
${\cal A}\cap (O\!\times\! [0,1]^\omega)$, ensuring that $y_{t1}\!\not=\! y_{t0}$, which is possible since ${\cal A}\cap (O\!\times\! [0,1]^\omega)$ has an uncountable second projection.

\vfill\eject

 Then we choose $\Boraone$ sets $X_{t0}$ and $X_{t1}$ with disjoint projections and diameter at most $2^{-l-1}$ such that 
$p_{t\varepsilon}\!\in\! X_{t\varepsilon}\!\subseteq\!\overline{X_{t\varepsilon}}\!\subseteq\! X_t$, 
and $S_{t\varepsilon}\!\in\!\Ana$ with GH-diameter at most $2^{-l-1}$ containing 
$p_{t\varepsilon}$ and contained in $X_{t\varepsilon}$ and the right $\Ana$ set.\hfill{$\square$}\bigskip

 We now study the case $\xi\! =\! 2$. The following lemma is a variant of the Mycielski-Kuratowski theorem (see 19.1 in [K1]). Recall the notation after the proof of Theorem \ref{Hur+} in the successor case.\bigskip

\noindent \bf Lemma\it\ \label{closedndsect} Let $F$ be a symmetric closed relation on 
$2^\omega$ with nowhere dense sections. Then there is 
$f\! :\! 2^\omega\!\rightarrow\! 2^\omega$ injective continuous such that $\mathbb{P}_f\! =\! f^{-1}(\mathbb{P}_f)$ and $\big( f(\alpha ),f(\beta )\big)\!\notin\! F$ if $\alpha\!\not=\!\beta$.\rm\bigskip
 
\noindent\bf Proof.\rm\ We inductively construct a sequence $(n_t)_{t\in 2^{<\omega}}$ of natural numbers, and a sequence $(U_t)_{t\in 2^{<\omega}}$ of clopen subsets of $2^\omega$ satisfying the following conditions:
$$\begin{array}{ll}
& (1)~U_{t\varepsilon}\!\subseteq\! U_t\cr
& (2)~\alpha_{n_t}\!\in\! U_t\cr
& (3)~\mbox{diam}(U_t)\!\leq\! 2^{-\vert t\vert}\cr
& (4)~U_{t0}\cap U_{t1}\! =\!\emptyset\cr
& (5)~n_{t0}\! =\! n_t\cr
& (6)~U_{t1}\cap\{\alpha_n\mid n\!\leq\!\vert t\vert\}\! =\!\emptyset\cr
& (7)~(U_{t0}\!\times\! U_{t1})\cap F\! =\!\emptyset
\end{array}$$
Assume that this is done. Using (1)-(3), we define $f\! :\! 2^\omega\!\rightarrow\! 2^\omega$ by 
$\{ f(\beta )\}\! :=\!\bigcap_{n\in\omega}~U_{\beta\vert n}$, and $f$ is injective continuous by (4). If 
$t\!\in\! O$ and $\alpha\! =\! t0^\infty$, then $f(\alpha )\! =\!\alpha_{n_t}$ by (5). If 
$\beta\!\in\!\mathbb{P}_\infty$, then there is an infinite strictly increasing sequence 
$(l_k)_{k\in\omega}$ of natural numbers with $\beta\vert l_k\!\in\! O$. Condition (6) implies that 
$f(\beta )\!\in\!\mathbb{P}_\infty$. Condition (7) implies that $\big( f(\alpha ),f(\beta )\big)\!\notin\! F$ if $\alpha\!\not=\!\beta$, by symmetry. So we are done.\bigskip

 Let us prove that the construction is possible. For the first step of the induction, we set  
$n_\emptyset\! :=\! 0$ and $U_\emptyset\! :=\! 2^\omega$. Assume that $(n_t)_{\vert t\vert\leq l}$ and 
$(U_t)_{\vert t\vert\leq l}$ satisfying (1)-(7) have been constructed, which is the case for $l\! =\! 0$.\bigskip

  Let $t\!\in\! 2^l$. Condition (5) defines $n_{t0}$. As $F$ has nowhere dense vertical sections, we can choose $n_{t1}$ in such a way that $\alpha_{n_{t1}}\!\in\! U_t\!\setminus\! 
(\{\alpha_{n_t}\}\cup\{\alpha_n\mid n\!\leq\! l\}\cup F_{\alpha_{n_t}})$. Then we choose disjoint clopen sets $U_{t0},U_{t1}$ with diameter at most $2^{-l-1}$ such that 
$\alpha_{n_{t\varepsilon}}\!\in\! U_{t\varepsilon}\!\subseteq\! U_t$ and satisfying (1)-(7).
\hfill{$\square$}\bigskip

 If $\Delta (X)\cap A$ is not separable from $\Delta (X)\cap B$ by a $\bf\Gamma$ set, then  Theorem \ref{Lo-SR} gives $f\! :\! 2^\omega\!\rightarrow\! X^2$ injective continuous with $\mathbb{C}\!\subseteq\! f^{-1}\big(\Delta (X)\cap A\big)$ and $\neg\mathbb{C}\!\subseteq\! f^{-1}\big(\Delta (X)\cap B\big)$. If 
$\alpha\!\not=\!\beta\!\in\! 2^\omega$, then $f(\alpha )\!\not=\! f(\beta )$. As $f(\alpha ),f(\beta )\!\in\!\Delta (X)$, 
$f_\varepsilon (\alpha )\!\not=\! f_\varepsilon (\beta )$ for each $\varepsilon\!\in\! 2$ $\big( f\! =\! (f_0,f_1)\big)$.\bigskip

 If $\Delta (X)\cap A$ is separable from $\Delta (X)\cap B$ by a $\bf\Gamma$ set, then 
$A\!\setminus\!\Delta (X)$ is not separable from $B\!\setminus\!\Delta (X)$ by a $\bf\Gamma$ set. Theorem \ref{Lo-SR} gives $f\! :=\! (f_0,f_1)\! :\! 2^\omega\!\rightarrow\! X^2$ injective continuous such that $\mathbb{C}\!\subseteq\! f^{-1}\big( A\!\setminus\!\Delta (X)\big)$ and 
$\neg\mathbb{C}\!\subseteq\! f^{-1}\big( B\!\setminus\!\Delta (X)\big)$. By the choice of $\mathbb{C}$, we may assume that $f_0$ and $f_1$ have disjoint ranges.\bigskip

 Note that the $f_\varepsilon$'s are nowhere dense-to-one. Indeed, we argue by contradiction, which gives $\varepsilon\!\in\! 2$ and $s\!\in\! 2^{<\omega}$ such that $f_\varepsilon$ is constant on $N_s$ and equal to $x\!\in\! X$. Assume for example that $\varepsilon\! =\! 0$. As $f$ is injective continuous, 
${f_1}_{\vert N_s}$ is also injective continuous. Note also that $f_1[\mathbb{C}\cap N_s]\! =\! f_1[N_s]\cap A_x$. As $A_x$ is in $\bf\Gamma$, so are 
 $f_1[\mathbb{C}\cap N_s]$ and $\mathbb{C}\cap N_s$, which contradicts the choice of $\mathbb{C}$.
 
\vfill\eject

 We next prove that we may assume that $f_\varepsilon$ is injective for each $\varepsilon\!\in\! 2$. In order to do this, we set 
$$F\! :=\!\{ (\alpha ,\beta )\!\in\! 2^\omega\!\times\! 2^\omega\mid\exists\varepsilon\!\in\! 2~~f_\varepsilon (\alpha )\! =\! f_\varepsilon (\beta )\} .$$ 
Note that $F$ is a symmetric closed relation on $2^\omega$ with nowhere dense sections. We apply Lemma \ref{closedndsect} to $F$, which gives $\psi\! :\! 2^\omega\!\rightarrow\! 2^\omega$ injective continuous such that $\mathbb{P}_f\! =\!\psi^{-1}(\mathbb{P}_f)$ and 
$\big(\psi (\alpha ),\psi (\beta )\big)\!\notin\! F$ if $\alpha\!\not=\!\beta$. This proves that we may assume that the $f_\varepsilon$'s are injective.\bigskip

 Assume finally that $\xi\! =\! 1$ and ${\bf\Gamma}\! =\!\boraone$, the other case being similar. As $A$ is not separable from $B$ by an open set, we can find 
$(\alpha ,\beta )\!\in\! A$ and $\big( (\alpha_n,\beta_n)\big)_{n\in\omega}\!\in\! B^\omega$ converging to $(\alpha ,\beta )$. If $\alpha_n\! =\!\alpha$ for all but finitely many $n$'s, then 
$A_\alpha$ is not open, which is absurd. So, extracting a subsequence if necessary, we may assume that the sequence $(\alpha_n)_{n\in\omega}$ is made of pairwise distinct elements different from $\alpha$. Similarly, we may assume that the sequence $(\beta_n)_{n\in\omega}$ is made of pairwise distinct elements different from $\beta$. It remains to set 
$f(0)\! :=\! (\alpha ,\beta )$ and $f(2^{-n})\! :=\! (\alpha_n,\beta_n)$.\bigskip

\noindent\bf (B) The limit case\rm\bigskip

 Condition (1) is changed as follows:
$$(1')~\left\{\!\!\!\!\!\!\!\!
\begin{array}{ll}
& \overline{X_t}\!\subseteq\! X_s\mbox{ if }~s~R^0~t\wedge s\!\not=\! t\cr
& S_t\!\subseteq\! S_s\mbox{ if }~s~R^\xi ~t\wedge\big( s,t\!\in\! {\cal I}\vee s,t\!\notin\! {\cal I}\big)
\end{array}
\right.$$
\bf Claim 2\rm ~(a) The set 
$S_{s^\xi}\cap\bigcap_{1\leq\rho <\xi (s)}~\overline{S_{s^\rho}}^{T_\rho}\cap X_{s^0}$ is $\Ana$ and nonempty.\smallskip

(b) If $s^\xi\!\in\! {\cal I}$, then the set $(B\!\setminus\! U)\cap
\bigcap_{1\leq\rho\leq\xi (s)}~\overline{S_{s^\rho}}^{T_\rho}\cap X_{s^0}\cap\Omega_S$ is $\Ana$ and nonempty.\rm\bigskip

 We conclude as in the proof of Theorem \ref{Hur+}.\hfill{$\square$}

\section{$\!\!\!\!\!\!$ Preserving $\mathbb{C}$ and avoiding countably many Borel graphs of functions}\indent

 The next lemma is essentially due to Louveau, even if it is not formally written like this in [Lo1] and [Lo2].

\begin{lem} \label{potclo} (Louveau) Let $1\!\leq\!\xi\! <\!\omega_1^{\mbox{CK}}$, $X,Y$ be recursively presented Polish spaces, and $B$ be $\Borel$ subset of $X\!\times\! Y$ whose vertical sections are $\bormxi$. Then 
$B\!\in\!\bormxi\big( (X,<\Borel >)\!\times\! Y\big)$.\end{lem}

\noindent\bf Proof.\rm\ Theorem 3.6 in [Lo1] provides a Polish topology $\tau$ on $X$ finer than $T_1$ such that 
$B$ is in $\bormxi\big( (X,\tau )\!\times\! (Y,T_1)\big)$. We then argue as in the proof of Theorem 3.4 in [Lo2] (we use the notation and material in this paper). Note first that $B\!\in\!\big(\borel\!\times\!\boraone\big)_\xi$ if 
$\xi\! <\!\omega$, $\big(\borel\!\times\!\boraone\big)_{\xi +1}$ otherwise. By Example 2 of Chapter 3 in [Lo2], the family $\big( N(n,Y)\big)_{n\in\omega}$ is regular. By Corollary 2.10 in [Lo2], $\bormxi (Y)$, as well as 
$\boraxi (Y)\! =\!\big(\bigcup_{\eta <\xi}~\borme (Y)\big)_\sigma$, are regular.\bigskip

 By Theorem 2.12 in [Lo2], $\borel\!\times\!\boraone$ is also regular. By Corollary 2.10.(v) in [Lo2], 
$\big(\borel\!\times\!\boraone\big)_\xi$ is also regular. The equality $W^\xi_\Phi\! =\! W_{\Phi_\xi}$ of this corollary, applied to $\Phi\! :=\!\borel\!\times\!\boraone$, shows that $B\!\in\!\big(\Borel\!\times\! (\boraone\cap\Borel )\big)_\xi$ if $\xi\! <\!\omega$, $\big(\Borel\!\times\! (\boraone\cap\Borel )\big)_{\xi +1}$ otherwise, and 
$B\!\in\!\bormxi\big( (X,<\Borel >)\!\times\! Y\big)$.\hfill{$\square$}\bigskip

 In order to prove Theorem \ref{ctble}, the main lemma is as follows.
 
\begin{thm} \label{presC-R}  Let ${\bf\Gamma}$ be a non self-dual Borel class of rank 
$3\!\leq\!\xi\! <\!\omega_1^{\mbox{CK}}$, $\mathbb{C}\!\in\!\Borel\cap\check {\bf\Gamma}(2^\omega )$, and $R$ be a $\Borel$ relation on $2^\omega$ with $F_\sigma$ vertical sections. We assume that there is a 
$\Ana$ subset $V$ of $2^\omega$ disjoint from $\Borel\cap 2^\omega$ such that $R\cap V^2$ is 
$\Sigma_{2^\omega}^2$-meager in $V^2$, and $V\cap\mathbb{C}$ is not separable from 
$V\!\setminus\!\mathbb{C}$ by a set in $\bf\Gamma$. Then there is $f\! :\! 2^\omega\!\rightarrow\! 2^\omega$ injective continuous such that $\mathbb{C}\! =\! f^{-1}(\mathbb{C})$ and 
$\big( f(\alpha ),f(\beta )\big)\!\notin\! R$ if $\alpha\!\not=\!\beta$.\end{thm}

\noindent\bf Proof.\rm\ By Theorem 3.5.(ii) in [Lo1], there is an increasing sequence 
$(F_l)_{l\in\omega}$ of $\Borel$ sets with closed vertical sections whose union is $R$.\bigskip

 Assume first that ${\bf\Gamma}\! =\!\boraxi$, so that 
$\mathbb{C}\!\in\!\bormxi (2^\omega )\!\setminus\!\boraxi$. We set 
$N\! :=\! V\cap\mathbb{C}\cap\overline{V\!\setminus\!\mathbb{C}}^{T_\xi}$. By Lemma \ref{top}, 
$N$ is $\Ana$ and nonempty.\bigskip

\noindent\bf (A) The successor case\rm\bigskip

 As in Section 2, we represent $h[\mathbb{C}]$ and set 
${\cal I}\! :=\!\{ s\!\in\! 2^{<\omega}\mid N^{R^\eta}_s\cap {\cal C}\!\not=\!\emptyset\}$, so that 
$\emptyset\!\in\! {\cal I}$. We construct a sequence $(X_s)_{s\in 2^{<\omega}}$ of nonempty 
$\Boraone$ subsets of $2^\omega$, and a sequence $(S_s)_{s\in 2^{<\omega}}$ of nonempty 
$\Ana$ subsets of $2^\omega$ satisfying the following conditions.
$$\begin{array}{ll}
& (1)~\left\{\!\!\!\!\!\!\!\!
\begin{array}{ll}
& \overline{X_t}\!\subseteq\! X_s\mbox{ if }~s~R^0~t\wedge s\!\not=\! t\cr
& S_t\!\subseteq\! S_s\mbox{ if }~s~R^\eta ~t\wedge\big( s,t\!\in\! {\cal I}\vee s,t\!\notin\! {\cal I}\big)
\end{array}
\right.\cr\cr
& (2)~S_s\!\subseteq\! X_s\cap\Omega_{2^\omega}\cap V\cr\cr
& (3)~\mbox{diam}(X_s),\mbox{diam}_{\mbox{GH}}(S_s)\! <\! 2^{-\vert s\vert}\cr\cr
& (4)~X_{s0}\cap X_{s1}\! =\!\emptyset\cr\cr
& (5)~S_s\!\subseteq\!\left\{\!\!\!\!\!\!\!\!
\begin{array}{ll}
& N\mbox{ if }~s\!\in\! {\cal I}\cr
& \neg\mathbb{C}\mbox{ if }s\!\notin\! {\cal I}
\end{array}
\right.\cr\cr
& (6)~S_t\!\subseteq\!\overline{S_s}^{T_\rho}\mbox{ if }s~R^\rho ~t\wedge 1\!\leq\!\rho\!\leq\!\eta
\cr\cr
& (7)~(S_s\!\times\! X_t)\cap F_l\! =\!\emptyset\mbox{ if }s\!\not=\! t\!\in\! 2^l
\end{array}$$
As in Section 2, it is enough to prove that the construction is possible. Indeed, fix 
$\alpha\!\not=\!\beta$. Then the inequality $\alpha\vert l\!\not=\!\beta\vert l$ holds if $l\!\geq\! L_0$. We set $l_k\! :=\!\vert\Pi^{-1}(h(\alpha ))(k)\vert$, so that $\Pi^{-1}(h(\alpha ))(k)\! =\!\alpha\vert l_k$, for each $k\!\in\!\omega$. As in Section 2, there is $k_0\!\in\!\omega$ such that $l_k\!\geq\! L_0$ if $k\!\geq\! k_0$ and $\big( f(\alpha ),f(\beta )\big)$ is in 
$(\bigcap_{k\geq k_0}~S_{\alpha\vert l_k})\!\times\! (\bigcap_{l\in\omega}~X_{\beta\vert l})$. Thus 
$\big( f(\alpha ),f(\beta )\big)\!\in\!\bigcap_{k\geq k_0}~
(S_{\alpha\vert l_k}\!\times\! X_{\beta\vert l_k})$. By (7), $\big( f(\alpha ),f(\beta )\big)$ is not in 
$\bigcup_{k\geq k_0}~F_{l_k}$. Therefore $\big( f(\alpha ),f(\beta )\big)\!\notin\! R$.\bigskip

 We first choose $\alpha_\emptyset\!\in\! N\cap\Omega_{2^\omega}$, $X_\emptyset$ semi-recursive with diameter at most $1$ containing $\alpha_\emptyset$, and $S_\emptyset\!\in\!\Ana$ with GH-diameter at most $1$ containing $\alpha_\emptyset$ and contained in 
 $X_\emptyset\cap N\cap\Omega_{2^\omega}$. Assume that our objects satisfying (1)-(7) have been contructed up to the length $l$, which is the case for $l\! =\! 0$.\bigskip

\noindent\bf Claim.\it ~(a) 
$S_{s^\eta}\cap\bigcap_{1\leq\rho <\eta}~\overline{S_{s^\rho}}^{T_\rho}\cap X_{s^0}$ is $\Ana$ and uncountable.\smallskip

(b) If $s^\eta\!\in\! {\cal I}$, then $\bigcap_{1\leq\rho\leq\eta}~\overline{S_{s^\rho}}^{T_\rho}\cap 
X_{s^0}\cap\Omega_{2^\omega}\cap V\!\setminus\!\mathbb{C}$ is $\Ana$ and uncountable.\rm\bigskip

\noindent\bf Subase 2.1\rm\ $s\!\in\! {\cal I}$, which implies that $s^\eta\!\in\! {\cal I}$.\bigskip

 We choose $\alpha_s\!\in\! S_{s^\eta}\cap
 \bigcap_{1\leq\rho <\eta}~\overline{S_{s^\rho}}^{T_\rho}\cap X_{s^0}$.\bigskip
 
\noindent\bf Subcase 2.2\rm\ $s^\eta\!\notin\! {\cal I}$, which implies that $s\!\notin\! {\cal I}$.\bigskip

 We argue as in Case 1.\bigskip

\noindent\bf Subcase 2.3\rm\ $s^\eta\!\in\! {\cal I}$ and $s\!\notin\! {\cal I}$.\bigskip

 We choose $\alpha_s\!\in\!\bigcap_{1\leq\rho\leq\eta}~\overline{S_{s^\rho}}^{T_\rho}\cap X_{s^0}\cap\Omega_{2^\omega}\cap V\!\setminus\!\mathbb{C}$.\bigskip

 This finishes the construction of $\alpha_{t0}$, which is in the right uncountable $\Ana$ set. The construction of $\alpha_{t1}$ is similar, the difference being that we ensure morover that 
$\alpha_{t1}\!\not=\!\alpha_{t0}$, which is possible since the right $\Ana$ set is uncountable. Then we choose disjoint $\Boraone$ sets $\tilde X_{t0}$ and $\tilde X_{t1}$ with diameter at most 
$2^{-l-1}$ such that $\alpha_{t\varepsilon}\!\in\!\tilde  X_{t\varepsilon}\!\subseteq\!
\overline{\tilde  X_{t\varepsilon}}\!\subseteq\! X_t$. We choose, for each $u\!\in\! 2^{l+1}$, a 
$\Ana$ set $\tilde S_u$ with GH-diameter at most $2^{-l-1}$, containing $\alpha_u$ and contained in $\tilde X_u$ and the right $\Ana$ set. Doing this, we ensured (1)-(6).\bigskip

 It remains to ensure (7). Assume that $s\!\not=\! t\!\in\! 2^{l+1}$. We first note that 
$F_{l+1}\cap V^2$ is ${\it\Sigma}_{2^\omega}^2$-meager in $V^2$. By Theorem 4.2 in [Ha-K-L] and 8.16 in [K1], $(2^\omega ,{\it\Sigma}_{2^\omega})$, $(V,{\it\Sigma}_{2^\omega})$ and 
$(V,{\it\Sigma}_{2^\omega})^2$ are strong Choquet. By 8.15 and 8.11 in [K1], 
$(V,{\it\Sigma}_{2^\omega})^2$ is a Baire space. This implies that 
$\tilde S_s\!\times\!\tilde S_t$ is not contained in $F_{l+1}$. By Lemma \ref{potclo}, $F_{l+1}$ is closed for $<\Borel >\!\times T_1$, and thus for ${\it\Sigma}_{2^\omega}^2$. So we can choose nonempty $\Ana$ sets $S'_s,S'_t$ contained in $\tilde S_s,\tilde S_t$ respectively with 
$(S'_s\!\times\! S'_t)\cap F_{l+1}\! =\!\emptyset$. After finitely many steps, we can ensure that this holds for any $s,t$. We pick $\beta_u\!\in\! S'_u$ for each $u\!\in\! 2^{l+1}$, so that 
$(\beta_s,\beta_t)\!\notin\! F_{l+1}$ for any $s,t$. As $F_{l+1}$ is closed for $<\Borel >\!\times T_1$, it is also closed for  ${\it\Sigma}_{2^\omega}\!\times\! T_1$. This gives, for each $s,t$, 
$S''_{s,t}\!\in\!\Ana$ and $X''_{s,t}\!\in\!\Boraone$ with 
$(\beta_s,\beta_t)\!\in\! S''_{s,t}\!\times\! X''_{s,t}\!\subseteq\! 
(S'_s\!\times\!\tilde X_t)\!\setminus\! F_{l+1}$. It remains to set 
$X_t\! :=\!\bigcap_{s\in 2^{l+1}}~X''_{s,t}$ and $S_s\! :=\!\bigcap_{t\in 2^{l+1}}~S''_{s,t}\cap X_s$.\bigskip

\noindent\bf (B) The limit case\rm\bigskip

 We indicate the differences with the successor case and the proof of Theorem \ref{Hur+}. We set 
$${\cal I}\! :=\!\{ s\!\in\! 2^{<\omega}\mid N^{R^\xi}_s\cap {\cal C}\!\not=\!\emptyset\} .$$
Conditions (1) and (6) are changed as follows:
$$\begin{array}{ll}
& (1')~\left\{\!\!\!\!\!\!\!\!
\begin{array}{ll}
& \overline{X_t}\!\subseteq\! X_s\mbox{ if }~s~R^0~t\wedge s\!\not=\! t\cr
& S_t\!\subseteq\! S_s\mbox{ if }~s~R^\xi ~t\wedge\big( s,t\!\in\! {\cal I}\vee s,t\!\notin\! {\cal I}\big)
\end{array}
\right.\cr\cr
& (6')~S_t\!\subseteq\!\overline{S_s}^{T_\rho}\mbox{ if }s~R^\rho ~t\wedge 1\!\leq\!\rho\!\leq\!\xi (s)
\end{array}$$
\bf Claim 2.\it ~(a) $S_{s^\xi}\cap\bigcap_{1\leq\rho <\xi (s)}~\overline{S_{s^\rho}}^{T_\rho}\cap X_{s^0}$ is $\Ana$ and uncountable.\smallskip

(b) If $s^\xi\!\in\! {\cal I}$, then 
$\bigcap_{1\leq\rho\leq\xi (s)}~\overline{S_{s^\rho}}^{T_\rho}\cap X_{s^0}\cap\Omega_{2^\omega}\cap V\!\setminus\!\mathbb{C}$ is $\Ana$ and uncountable.\rm\bigskip

 We conclude as above. Assume now that ${\bf\Gamma}\! =\!\bormxi$, so that 
$\mathbb{C}\!\in\!\boraxi (2^\omega )\!\setminus\!\bormxi$. We indicate the differences with the case ${\bf\Gamma}\! =\!\boraxi$. We set 
$N\! :=\!\overline{V\cap\mathbb{C}}^{T_\xi}\cap V\!\setminus\!\mathbb{C}$. By Lemma \ref{top}, 
$N$ is $\Ana$ and nonempty.

\vfill\eject

\noindent\bf (A) The successor case\rm\bigskip

 We represent $h[\neg\mathbb{C}]$ and set 
${\cal I}\! :=\!\{ s\!\in\! 2^{<\omega}\mid N^{R^\eta}_s\cap\Pi^{-1}(h[\neg \mathbb{C}])\!\not=\!\emptyset\}$, so that 
$\emptyset\!\in\! {\cal I}$. We ensure 
$$(5)~S_s\!\subseteq\!\left\{\!\!\!\!\!\!\!\!
\begin{array}{ll}
& N\mbox{ if }~s\!\in\! {\cal I}\cr
& \mathbb{C}\mbox{ if }s\!\notin\! {\cal I}
\end{array}
\right.$$
\bf Claim.\it ~(a) $S_{s^\eta}\cap\bigcap_{1\leq\rho <\eta}~\overline{S_{s^\rho}}^{T_\rho}\cap X_{s^0}$ is $\Ana$ and uncountable.\smallskip

(b) If $s^\eta\!\in\! {\cal I}$, then 
$(V\cap\mathbb{C})\cap\bigcap_{1\leq\rho\leq\eta}~\overline{S_{s^\rho}}^{T_\rho}\cap X_{s^0}\cap\Omega_{2^\omega}$ is $\Ana$ and uncountable.\rm\bigskip   

\noindent\bf (B) The limit case\rm\bigskip

 We set ${\cal I}\! :=\!\{ s\!\in\! 2^{<\omega}\mid N^{R^\xi}_s\cap\Pi^{-1}(h[\neg\mathbb{C}])\!\not=\!\emptyset\}$.\bigskip

\noindent\bf Claim 2.\it ~(a) $S_{s^\xi}\cap\bigcap_{1\leq\rho <\xi (s)}~\overline{S_{s^\rho}}^{T_\rho}\cap X_{s^0}$ is $\Ana$ and uncountable.\smallskip

(b) If $s^\xi\!\in\! {\cal I}$, then the set $(V\cap\mathbb{C})\cap\bigcap_{1\leq\rho\leq\xi (s)}~
\overline{S_{s^\rho}}^{T_\rho}\cap X_{s^0}\cap\Omega_{2^\omega}$ is $\Ana$ and uncountable.\bigskip\rm

 We conclude as above.\hfill{$\square$}
 
\begin{cor} \label{avoidctble} Let ${\bf\Gamma}$ be a non self-dual Borel class of rank at least three, 
$\mathbb{C}\!\in\!\check {\bf\Gamma}(2^\omega )$ not in $\bf\Gamma$, and $R$ be a Borel relation on 
$2^\omega$ with countable vertical sections. Then there is $f\! :\! 2^\omega\!\rightarrow\! 2^\omega$ injective continuous such that $\mathbb{C}\! =\! f^{-1}(\mathbb{C})$, and $\big( f(\alpha ),f(\beta )\big)\!\notin\! R$ if 
$\alpha\!\not=\!\beta$.\end{cor}

\noindent\bf Proof.\rm\ The Lusin-Novikov theorem gives a sequence $(f_n)_{n\in\omega}$ of partial Borel maps from $2^\omega$ into itself with 
$R\! =\!\bigcup_{n\in\omega}~\mbox{Graph}(f_n)$ (see 18.10 in [K1]). Let $\xi$ be the rank of 
$\bf\Gamma$. In order to simplify the notation, we assume that $\xi\! <\!\omega_1^{\mbox{CK}}$ and $\mathbb{C},(f_n)_{n\in\omega}$ are $\Borel$. We set $U\! :=\!\Borel\cap 2^\omega$. Lemma \ref{top2}.(a) shows that $U$ is countable and $\Ca$. We will apply Theorem \ref{presC-R} to the 
$\Ana$ set $V\! :=\!\neg U$. Note that $\mbox{Graph}(f_n)\cap V^2$ is 
${\it\Sigma}_{2^\omega}^2$-closed in $V^2$ with nowhere dense vertical sections by definition of $U$. By the Kuratowski-Ulam theorem (see 8.41 in [K1]), $R\cap V^2$ is 
${\it\Sigma}_{2^\omega}^2$-meager in $V^2$. It remains to note that $V\cap\mathbb{C}$ is not separable from $V\!\setminus\!\mathbb{C}$ by a set in $\bf\Gamma$ since $U$ is countable and therefore in $F_\sigma\!\subseteq\! {\bf\Gamma}\cap\check {\bf\Gamma}$.\hfill{$\square$}\bigskip

\noindent\bf Remark.\rm\ This corollary cannot be extended to lower levels. Indeed, for the rank one, as $\mathbb{K}$ is countable, $R$ can be $\mathbb{K}^2$. For ${\bf\Gamma}\! =\!\boratwo$, $R$ can be $(\neg\mathbb{C})^2$ since $\neg\mathbb{C}$ is countable. Similarly, if 
${\bf\Gamma}\! =\!\bormtwo$, then $R$ can be $\mathbb{C}^2$.

\section{$\!\!\!\!\!\!$ Some general facts}\indent 

\label{SectionGeneral}

 We first note the following topological properties.

\begin{lem} \label{basic} Let $\bf\Gamma$ be a class of sets closed under continuous pre-images, 
$Y$ be a topological space, and $F$ be an equivalence relation on $Y$.\smallskip

(a) if $F$ is in $\bf\Gamma$, then the equivalence classes of $F$ are also in $\bf\Gamma$, \smallskip

(b) if $Z$ is a topological space, $G$ is an equivalence relation on $Z$ whose classes are in 
$\bf\Gamma$, and $(Y,F)\leq_c(Z,G)$, then the equivalence classes of $F$ are also in 
$\bf\Gamma$.\end{lem}

\noindent\bf Proof.\rm\ (a) comes from the fact that if $y\!\in\! Y$, then the map 
$i_y\! :\! y'\!\mapsto\! (y,y')$ is continuous and satisfies $[y]_F\! =\! i_y^{-1}(F)$. The statement (b) comes from the fact that $[y]_F\! =\! f^{-1}\big( [f(y)]_G\big)$.\hfill{$\square$}\bigskip

 The introduction of $\mathbb{E}^{\bf\Gamma}_0,\mathbb{E}^{\bf\Gamma}_1$ and $\mathbb{E}^{\boraxi}_2$ is motivated by the following fact.

\begin{prop} \label{injcontA} Let ${\bf\Gamma}$ be a non self-dual Borel class, $\mathbb{K}$ as above, 
$\mathbb{C}\!\in\!\check {\bf\Gamma}(\mathbb{K})\!\setminus\! {\bf\Gamma}$ (as above if the rank of 
$\bf\Gamma$ is at most two), $X$ be an analytic space, and $E$ be a Borel equivalence relation on $X$. Then exactly one of the following holds:\smallskip  

(a) the equivalence classes of $E$ are in $\bf\Gamma$,\smallskip  

(b) there is a Borel equivalence relation $\mathbb{E}$ on $\mathbb{K}$ such that 
$\mathbb{E}^{\bf\Gamma}_0\!\subseteq\!\mathbb{E}\!\subseteq\!\mathbb{E}^{\bf\Gamma}_1$ and 
$(\mathbb{K},\mathbb{E})\sqsubseteq_c(X,E)$.\end{prop}

\noindent\bf Proof.\rm\ Note that (a) and (b) cannot hold simultaneously by Lemma \ref{basic} since (b) implies that $\mathbb{C}$ is an $\mathbb{E}$-class. Assume that (a) does not hold, which gives $x\!\in\! X$ such that $[x]_E\!\notin\! {\bf\Gamma}$. Theorem \ref{Lo-SR} gives 
$i\! :\!\mathbb{K}\!\rightarrow\! X$ injective continuous such that $\mathbb{C}\! =\! i^{-1}([x]_E)$. It remains to set $\mathbb{E}\! :=\! (i\!\times\! i)^{-1}(E)$ to finish the proof.\hfill{$\square$}\bigskip
   
 The introduction of our equivalence relations on $\mathbb{H}$ is motivated by the following facts.

\begin{thm} \label{contB} Let ${\bf\Gamma}$ be a non self-dual Borel class, $\mathbb{K},\mathbb{C}$ as above, $X$ be an analytic space, and $E$ be a Borel equivalence relation on $X$ whose sections are in 
$\bf\Gamma$. Then exactly one of the following holds:\smallskip  

(a) the relation $E$ is a $\bf\Gamma$ subset of $X^2$,\smallskip  

(b) there is a Borel equivalence relation $\mathbb{E}$ on $\mathbb{H}\! :=\! 2\!\times\!\mathbb{K}$ such that 
$\big\{\big( (0,\alpha ),(1,\alpha )\big)\mid\alpha\!\in\!\mathbb{C}\big\}\!\subseteq\!\mathbb{E}$, 
$\big\{\big( (0,\alpha ),(1,\alpha )\big)\mid\alpha\!\notin\!\mathbb{C}\big\}\!\subseteq\!\neg\mathbb{E}$ and 
$(\mathbb{H},\mathbb{E})\sqsubseteq_c(X,E)$.\end{thm}

\noindent\bf Proof.\rm\ We first note that (a) and (b) cannot hold simultaneously. Indeed, we argue by contradiction, so that $\mathbb{E}\!\in\! {\bf\Gamma}(\mathbb{H}^2)$, and $\mathbb{E}\cap\big\{\big( (0,\alpha ),(1,\alpha )\big)\mid\alpha\!\in\!\mathbb{K}\big\}\!\in\! 
{\bf\Gamma}\Big(\big\{\big( (0,\alpha ),(1,\alpha )\big)\mid\alpha\!\in\!\mathbb{K}\big\}\Big)$. This implies that 
$\mathbb{C}\!\in\!\ {\bf\Gamma}(\mathbb{K})$, which is absurd. Assume now that (a) does not hold, so that ${\bf\Gamma}\!\not=\!\boraone$.\bigskip

 Theorem \ref{coord} gives $f\! :=\! (f_0,f_1)\! :\!\mathbb{K}\!\rightarrow\! X^2$ continuous with injective coordinates with $\mathbb{C}\! =\! f^{-1}(E)$. If the rank of $\bf\Gamma$ is at least two, then we may assume that $f_0$ and $f_1$ have disjoint ranges, by the choice of $\mathbb{C}$. We define 
$g\! :\!\mathbb{H}\!\rightarrow\! X$ by $g(\varepsilon ,\alpha )\! :=\! f_\varepsilon (\alpha )$. Note that $g$ is continuous, 
$$\big\{\big( (0,\alpha ),(1,\alpha )\big)\mid\alpha\!\in\!\mathbb{C}\big\}\!\subseteq\! (g\!\times\! g)^{-1}(E)$$ 
and $\big\{\big( (0,\alpha ),(1,\alpha )\big)\mid\alpha\!\notin\!\mathbb{C}\big\}\!\subseteq\! (g\!\times\! g)^{-1}(\neg E)$. It remains to set 
$\mathbb{E}\! :=\! (g\!\times\! g)^{-1}(E)$.\bigskip

 If ${\bf\Gamma}\! =\!\bormone$, then $f(0)\!\notin\! E$, so that $f_0(0)\!\not=\! f_1(0)$ and $f_0(2^{-k})\!\not=\! f_1(2^{-k})$ if $k\!\geq\! k_0$. So here again we may assume that $f_0$ and $f_1$ have disjoint ranges, and we conclude as above.\hfill{$\square$}
 
\begin{prop} \label{foursame} Let ${\bf\Gamma}$ be a non self-dual Borel class, $\mathbb{K}$ as above, 
$\mathbb{C}\!\subseteq\!\mathbb{K}$, $\mathbb{E}$ be an equivalence relation on $\mathbb{H}$ containing 
$\big\{\big( (0,\alpha ),(1,\alpha )\big)\mid\alpha\!\in\!\mathbb{C}\big\}$, $\varepsilon,\eta\!\in\! 2$ and 
$\alpha ,\beta\!\in\!\mathbb{C}$. Then $\big( (\varepsilon,\alpha ),(\eta ,\beta )\big)\!\in\!\mathbb{E}$ is equivalent to $\big( (0,\alpha ),(0,\beta )\big)\!\in\!\mathbb{E}$.\end{prop}

\noindent\bf Proof.\rm\ We may assume that $\eta\! =\! 1$. Assume first that 
$\big( (\varepsilon,\alpha ),(1,\beta )\big)\!\in\!\mathbb{E}$. As 
$$\big( (0,\alpha ),(1,\alpha )\big) ,\big( (0,\beta ),(1,\beta )\big)\!\in\!\mathbb{E}\mbox{,}$$ 
$\big( (0,\alpha ),(0,\beta )\big)\!\in\!\mathbb{E}$. Conversely, assume that $\big( (0,\alpha ),(0,\beta )\big)\!\in\!\mathbb{E}$. Similarly, 
$\big( (0,\alpha ),(1,\beta )\big)\!\in\!\mathbb{E}$ and $\big( (1,\alpha ),(1,\beta )\big)\!\in\!\mathbb{E}$.\hfill{$\square$}\bigskip
 
 We now check a useful fact announced in the introduction.

\begin{lem} \label{true} Let ${\bf\Gamma}$ be a non self-dual Borel class of rank at least two. Then there is 
$\mathbb{C}\!\in\!\check {\bf\Gamma}(2^\omega)\!\setminus\! {\bf\Gamma}$ such that 
$\mathbb{C}\cap N_s\!\in\!\check {\bf\Gamma}(N_s)\!\setminus\! {\bf\Gamma}$ for each 
$s\!\in\! 2^{<\omega}$. In particular, $\mathbb{C}$ is dense and co-dense in $2^\omega$.\end{lem}

\noindent\bf Proof.\rm\ We may assume that ${\bf\Gamma}\! =\!\bormxi$ with $\xi\!\geq\! 2$, passing to complements if ${\bf\Gamma}\! =\!\boraxi$. We will inductively construct $\mathbb{C}_\xi\!\in\!\boraxi$ as required. As required in the introduction, we set 
$$\mathbb{C}_2\! :=\!\{\alpha\!\in\! 2^\omega\mid\forall^\infty n\!\in\!\omega ~~\alpha (n)\! =\! 0\} .$$ 
Note that $\mathbb{C}_2$ is dense and co-dense in $2^\omega$, and we are done for $\xi\! =\! 2$, by Baire's theorem. Let $3\!\leq\!\xi\! =\!\mbox{sup}_{n\in\omega}~(\xi_n\! +\! 1)$, with $2\!\leq\!\xi_n\! <\!\xi$. We set 
$$\mathbb{C}_\xi\! :=\!\{\alpha\!\in\! 2^\omega\mid\exists n\!\in\!\omega ~~(\alpha )_n\!\notin\!\mathbb{C}_{\xi_{(n)_0}}\} .$$ 
By 22.10 in [K1], it is enough to check that $\mathbb{C}_\xi\cap N_s$ reduces any $\boraxi$ subset $S$ of 
$2^\omega$. Assume first that $s\! =\!\emptyset$. Write $S\! =\!\bigcup_{n\in\omega}~\neg S_n$, where 
$S_n\!\in\! {\bf\Sigma}^0_{\xi_n}$. The induction assumption gives 
$f_n\! :\! 2^\omega\!\rightarrow\! 2^\omega$ continuous with $S_n\! =\! f_n^{-1}(\mathbb{C}_{\xi_n})$. We define $f\! :\! 2^\omega\!\rightarrow\! 2^\omega$ by 
$\big( f(\alpha )\big)_n\! :=\! f_{(n)_0}(\alpha )$, so that $f$ is continuous. Then 
$$\begin{array}{ll}
\alpha\!\in\! S\!\!\!
& \Leftrightarrow\exists n\!\in\!\omega ~~\alpha\!\notin\! S_n
\Leftrightarrow\exists n\!\in\!\omega ~~f_n(\alpha )\!\notin\!\mathbb{C}_{\xi_n}
\Leftrightarrow\exists n\!\in\!\omega ~~f_{(n)_0}(\alpha )\!\notin\!\mathbb{C}_{\xi_{(n)_0}}\cr
& \Leftrightarrow\exists n\!\in\!\omega ~~\big( f(\alpha )\big)_n\!\notin\!\mathbb{C}_{\xi_{(n)_0}}
\Leftrightarrow  f(\alpha )\!\in\!\mathbb{C}_\xi .
\end{array}$$
If now $s$ is arbitrary, then we define $g\! :\! 2^\omega\!\rightarrow\! N_s$ by
$$\big( g(\alpha )\big)_n\! :=\!\left\{\!\!\!\!\!\!\!\!
\begin{array}{ll}
& (s)_n0^\infty\mbox{ if }(n)_1\!\leq\!\vert s\vert\mbox{,}\cr
& (\alpha )_{<(n)_0,(n)_1-\vert s\vert -1>}\mbox{ if }(n)_1\! >\!\vert s\vert\mbox{,}
\end{array}
\right.$$
so that $g$ is continuous and reduces $\mathbb{C}_\xi$ to $\mathbb{C}_\xi\cap N_s$ since 
$$\begin{array}{ll}
\alpha\!\in\!\mathbb{C}_\xi\!\!\!
& \Leftrightarrow\exists n\!\in\!\omega ~~(\alpha )_n\!\notin\!\mathbb{C}_{\xi_{(n)_0}}
\Leftrightarrow\exists n,p\!\in\!\omega ~~(\alpha )_{<n,p>}\!\notin\!\mathbb{C}_{\xi_n}\cr
& \Leftrightarrow\exists n\!\in\!\omega ~~\exists p\! >\!\vert s\vert 
~~(\alpha )_{<n,p-\vert s\vert -1>}\!\notin\!\mathbb{C}_{\xi_n}\cr
& \Leftrightarrow\exists n,p\!\in\!\omega ~~\big( g(\alpha )\big)_{<n,p>}\!\notin\!\mathbb{C}_{\xi_n}\Leftrightarrow\exists n\!\in\!\omega ~~
\big( g(\alpha )\big)_n\!\notin\!\mathbb{C}_{\xi_{(n)_0}}\cr 
& \Leftrightarrow g(\alpha )\!\in\!\mathbb{C}_\xi\cap N_s.
\end{array}$$
This finishes the proof.\hfill{$\square$}\bigskip

\noindent\bf Notation.\rm\ If ${\bf\Gamma}$ is a non self-dual Borel class, then 
$D_2({\bf\Gamma})\! =\!\{ A\!\setminus\! B\mid A,B\!\in\! {\bf\Gamma}\}$, and 
$${\bf\Gamma}^+\! :=\!\{ (A\cap C)\cup (B\!\setminus\! C)\mid 
A\!\in\! {\bf\Gamma}\wedge B\!\in\!\check {\bf\Gamma}\wedge C\!\in\!\borone\}$$ 
is the successor of $\bf\Gamma$ in the Wadge quasi-order.

\vfill\eject

 In order to state the next result, we extend our sets ${\mathcal A}^{\bf\Gamma}$ and 
${\mathcal B}^{\bf\Gamma}$. We set 
$${\mathcal A}^{\bf\Gamma}\! :=\!\left\{\!\!\!\!\!\!\!
\begin{array}{ll}
& \{ (\mathbb{K},\mathbb{E}^{\bf\Gamma}_0)\}\mbox{ if }{\bf\Gamma}\! =\!\bormone\mbox{,}\cr
& \{ (\mathbb{K},\mathbb{E}^{\bf\Gamma}_n)\mid n\!\leq\! 1\}\mbox{ if }{\bf\Gamma}\! =\!\boraone
\mbox{ or the rank of }{\bf\Gamma}\mbox{ is two,}\cr
& \{ (\mathbb{K},\mathbb{E}^{\bf\Gamma}_n)\mid 1\!\leq\! n\!\leq\! 2\}\mbox{ if }
{\bf\Gamma}\!\in\!\{\boraxi\mid\xi\!\geq\! 3\}\mbox{,}\cr
& \{ (\mathbb{K},\mathbb{E}^{\bf\Gamma}_1)\}\mbox{ if }
{\bf\Gamma}\!\in\!\{\bormxi\mid\xi\!\geq\! 3\}\mbox{,}
\end{array}
\right.$$
$${\mathcal B}^{\bf\Gamma}\! :=\! {\mathcal A}^{\bf\Gamma}\cup
\left\{\!\!\!\!\!\!\!
\begin{array}{ll}
& \emptyset\mbox{ if }{\bf\Gamma}\! =\!\boraone\mbox{,}\cr
& \{ (\mathbb{H},\mathbb{E}^{\bf\Gamma}_3)\}\mbox{ if }{\bf\Gamma}\! =\!\bormone\mbox{,}\cr
& \{ (\mathbb{H},\mathbb{E}^{\bf\Gamma}_n)\mid 3\!\leq\! n\!\leq\! 5\}\mbox{ if the rank of }{\bf\Gamma}\mbox{ is two,}\cr
& \{ (\mathbb{H},\mathbb{E}^{\bf\Gamma}_8)\}\mbox{ if }{\bf\Gamma}\!\in\!\{\bormxi\mid\xi\!\geq\! 3\} .
\end{array}
\right.$$

\begin{thm} \label{prope} Let ${\bf\Gamma}$ be a non self-dual Borel class, 
$\mathbb{K},\mathbb{C}$ as above.\smallskip

\noindent (a) The following properties of 
$\mathbb{E}^{\bf\Gamma}_n\!\in\! {\mathcal B}^{\bf\Gamma}$ hold:
$$\begin{tabular}{|c|c|c|c|c|}
\hline
$n$ & $\bf\Gamma$ & Number of classes & Complexity of the classes & 
Complexity of the relation\\
\hline
0 & 
& $\begin{array}{ll}
& \cr
& 2\mbox{ if }{\bf\Gamma}\! =\!\bormone\cr 
& \omega\mbox{ if }{\bf\Gamma}\!\in\!\{\boraone ,\boratwo\}\cr 
& 2^\omega\mbox{ if }{\bf\Gamma}\!\supseteq\!\bormtwo\cr
& \end{array}$
& $\begin{array}{ll}
& {\bf\Gamma}^+\mbox{ if }{\bf\Gamma}\! =\!\bormone\cr 
& \check {\bf\Gamma}\mbox{ if }{\bf\Gamma}\!\not=\!\bormone\end{array}$
& $\begin{array}{ll}
& \check D_2({\bf\Gamma})\!\setminus\! D_2({\bf\Gamma})\mbox{ if }{\bf\Gamma}\! =\!\bormone\cr 
& \check {\bf\Gamma}\!\setminus\! {\bf\Gamma}\mbox{ if }{\bf\Gamma}\!\not=\!\bormone\end{array}$
\\
\hline
1 & 
& 2
& ${\bf\Gamma}^+$
& $\begin{array}{ll}
& \cr
& \check D_2({\bf\Gamma})\!\setminus\! D_2({\bf\Gamma})\mbox{ if }{\bf\Gamma}\! =\!\boraone\cr 
& \check D_2({\bf\Gamma})\!\setminus\!  ({\bf\Gamma}\cup\check {\bf\Gamma})
\mbox{ if rk}({\bf\Gamma})\!\geq\! 2\cr
& 
\end{array}$
\\
\hline
2 & $\boraxi$
& $\omega$
& $\bormxi$
& $\begin{array}{ll}
& \cr
& \check D_2(\boraxi )\!\setminus\!\boraxi\cr
& 
\end{array}$
\\
\hline
3 & 
& $\begin{array}{ll}
& \cr
& \omega\mbox{ if }{\bf\Gamma}\! =\!\bormone\cr 
& 2^\omega\mbox{ if rk}({\bf\Gamma})\!\geq\! 2\cr
& \end{array}$
& $\bormone$
& $\begin{array}{ll}
& {\bf\Gamma}^+\!\setminus\! ({\bf\Gamma\cup\check {\bf\Gamma}})\mbox{ if }
{\bf\Gamma}\! =\!\bormone\cr 
& \check {\bf\Gamma}\!\setminus\! {\bf\Gamma}\mbox{ if rk}({\bf\Gamma})\!\geq\! 2\end{array}$
\\
\hline
4 & 
& $2^\omega$
& $\bf\Gamma$
& $\begin{array}{ll}
& \cr 
& {\bf\Gamma}^+\!\setminus\! ({\bf\Gamma}\cup\check {\bf\Gamma})\cr
& \end{array}$\\
\hline
5 & 
& $\begin{array}{ll}
& \cr
& \omega\mbox{ if }{\bf\Gamma}\! =\!\bormtwo\cr 
& 2^\omega\mbox{ if }{\bf\Gamma}\!\supseteq\!\boratwo\cr
& \end{array}$
& $\bf\Gamma$
& ${\bf\Gamma}^+\!\setminus\! ({\bf\Gamma}\cup\check {\bf\Gamma})$\\
\hline
8 & $\bormxi$
& $\omega$
& $\bormxi$
& $\begin{array}{ll}
& \cr
& \check D_2(\boraxi )\!\setminus\! (\boraxi\cup\bormxi )\cr
& \end{array}$\\
\hline
\end{tabular}$$
(b) ${\mathcal A}^{\bf\Gamma}$ and ${\mathcal B}^{\bf\Gamma}$ are $\leq_c$-antichains.
\end{thm}

\noindent\bf Proof.\rm\ (a).(0) Note that the equivalence classes of $\mathbb{E}^{\bf\Gamma}_0$ are 
$\mathbb{C}$, and $\{ x\}$ for $x\!\notin\!\mathbb{C}$. Note also that 
${\mathbb{E}^{\boraone}_0\! =\!\Delta (\mathbb{K})}$. If ${\bf\Gamma}\!\not=\!\bormone$, then 
$\mathbb{E}^{\bf\Gamma}_0$ is in $\check {\bf\Gamma}\!\setminus\! {\bf\Gamma}$, and its equivalence classes are not all in $\bf\Gamma$, and all in $\check {\bf\Gamma}$.\bigskip

\noindent (1) Note that the equivalence classes of $\mathbb{E}^{\bf\Gamma}_1$ are $\mathbb{C}$ and 
$\mathbb{K}\!\setminus\!\mathbb{C}$. In particular, $\mathbb{E}^{\bf\Gamma}_1$ is not in $\bf\Gamma$, not in $\check {\bf\Gamma}$, and its equivalence classes are all in ${\bf\Gamma}^+$, not all in $\bf\Gamma$, and not all in $\check {\bf\Gamma}$. However, it is in $\check D_2({\bf\Gamma})$. If 
${\bf\Gamma}\! =\!\boraone$, as $\mathbb{K}\!\setminus\!\{ 0\}$ is dense in $\mathbb{K}$, 
$\mathbb{E}^{\boraone}_1$ is not in $D_2(\boraone )$. For ${\bf\Gamma}\! =\!\bormone$, note that 
$\mathbb{E}^{\bormone}_1\! =\!\mathbb{E}^{\boraone}_1$.\bigskip

\noindent (2) The equivalence classes of $\mathbb{E}^{\boraxi}_2$ are $\mathbb{C}$ and the 
$\mathbb{C}_n$'s.\bigskip

\noindent (3) The equivalence classes of $\mathbb{E}^{\bf\Gamma}_3$ are $2\!\times\!\{ x\}$ for 
$x\!\in\!\mathbb{C}$, and $\{ (\varepsilon ,x)\}$ for $\varepsilon\!\in\! 2$ and $x\!\notin\!\mathbb{C}$, and thus closed. If ${\bf\Gamma}\! =\!\bormone$, then $\mathbb{E}^{\bf\Gamma}_3$ is $\check D_2({\bf\Gamma})$. It is not closed since 
$$\overline{\mathbb{E}^{\bf\Gamma}_3}\!\setminus\!\mathbb{E}^{\bf\Gamma}_3\! =\!
\big\{\big( (0,0),(1,0)\big) ,\big\{\big( (1,0),(0,0)\big)\big\} .$$ 
In particular, $\mathbb{E}^{\bf\Gamma}_3$ is $D_2({\bf\Gamma})$. It is not open since 
$\big( (0,0),(0,0)\big)\!\in\!\mathbb{E}^{\bf\Gamma}_3\cap
\overline{\neg\mathbb{E}^{\bf\Gamma}_3}$. So the exact complexity of 
$\mathbb{E}^{\bf\Gamma}_3$ is ${\bf\Gamma}^+$.\bigskip

\noindent (4) The equivalence classes of $\mathbb{E}^{\bf\Gamma}_4$ are $2\!\times\!\{ x\}$ for 
$x\!\in\!\mathbb{C}$, $\{ (0,x)\}$ for $x\!\notin\!\mathbb{C}$, and $\{1\}\!\times\! (\neg\mathbb{C})$, and thus in $\bf\Gamma$ if ${\bf\Gamma}\!\not=\!\boraone$, ${\bf\Gamma}^+$ otherwise.\bigskip

\noindent (5) Note that the equivalence classes of $\mathbb{E}^{\bf\Gamma}_5$ are $2\!\times\!\{ x\}$ for 
$x\!\in\!\mathbb{C}$, and $\{\varepsilon\}\!\times\! (\neg\mathbb{C})$ for $\varepsilon\!\in\! 2$, and thus in 
$\bf\Gamma$ if ${\bf\Gamma}\!\not=\!\boraone$, ${\bf\Gamma}^+$ otherwise.\bigskip

\noindent (8) The equivalence classes of $\mathbb{E}^{\bormxi}_8$ are $2\!\times\!\mathbb{C}_n$ for 
$n\!\in\!\omega$, and $\{\varepsilon\}\!\times\! (\neg\mathbb{C})$ for $\varepsilon\!\in\! 2$.\bigskip

\noindent (b).(1) Assume that ${\bf\Gamma}\!\not=\!\bormone$. Note that 
$(\mathbb{K},\mathbb{E}^{\bf\Gamma}_1)$ is not $\leq_c$-below 
$(\mathbb{K},\mathbb{E}^{\bf\Gamma}_0)$ since $\mathbb{E}^{\bf\Gamma}_0$ is in 
$\check {\bf\Gamma}$ and $\mathbb{E}^{\bf\Gamma}_1$ is not. Moreover, 
$(\mathbb{K},\mathbb{E}^{\bf\Gamma}_0)$ is not $\leq_c$-below 
$(\mathbb{K},\mathbb{E}^{\bf\Gamma}_1)$ since $\mathbb{E}^{\bf\Gamma}_0$ has infinitely many classes and $\mathbb{E}^{\bf\Gamma}_1$ has only two classes.\bigskip

 Assume now that $n\! >\! 1$. Similarly, $\mathbb{E}^{\bf\Gamma}_n$ is not below 
$\mathbb{E}^{\bf\Gamma}_1$. Conversely, as the classes of 
$\mathbb{E}^{\bf\Gamma}_n\!\in\! {\mathcal B}^{\bf\Gamma}$ are all in $\bf\Gamma$ or all in 
$\check {\bf\Gamma}$, $\mathbb{E}^{\bf\Gamma}_1$ is not below $\mathbb{E}^{\bf\Gamma}_n$. Thus 
$\mathbb{E}^{\bf\Gamma}_1$ is incomparable with the other relations in ${\mathcal B}^{\bf\Gamma}$. In particular, 
${\mathcal A}^{\bf\Gamma}$ is a $\leq_c$-antichain.\bigskip

\noindent (0) $\mathbb{E}^{\bf\Gamma}_0$ is not below the other relations in ${\mathcal B}^{\bf\Gamma}$, because of the complexity of the classes.\bigskip

 Let us prove that $\mathbb{E}^{\bf\Gamma}_3$ is not below $\mathbb{E}^{\bf\Gamma}_0$ if 
${\bf\Gamma}\!\not=\!\boraone$. We argue by contradiction, which gives 
$f\! :\!\mathbb{H}\!\rightarrow\!\mathbb{K}$. As ${\bf\Gamma}\!\not=\!\boraone$, $\mathbb{C}$ is dense in 
$\mathbb{K}$. This gives $\alpha\!\in\!\mathbb{C}$ with $f(0,\alpha )\!\not=\! f(1,\alpha )$, since otherwise 
$$f(0,\beta )\! =\! f(1,\beta )$$ 
for each $\beta\!\in\! 2^\omega$, and thus 
$\big( (0,\beta ),(1,\beta )\big)\!\in\!\mathbb{E}^{\bf\Gamma}_3$ for some $\beta\!\notin\!\mathbb{C}$, which cannot be. In particular, $f(0,\alpha ),f(1,\alpha )\!\in\!\mathbb{C}$. Similarly, working in $\mathbb{C}\cap N_{1-\alpha (0)}$ if necessary, we can find $\beta\!\in\!\mathbb{C}\!\setminus\!\{\alpha\}$ with $f(0,\beta )\!\not=\! f(1,\beta )$ and $f(0,\beta ),f(1,\beta )\!\in\!\mathbb{C}$. As $\big( f(0,\alpha ),f(1,\beta )\big)\!\in\!\mathbb{E}^{\bf\Gamma}_0$, 
$\big( (0,\alpha ),(1,\beta )\big)\!\in\!\mathbb{E}^{\bf\Gamma}_3$, which is absurd.\bigskip

 Appealing to the number of classes or the complexity of the relations, we see that 
$\mathbb{E}^{\bf\Gamma}_0$ is above neither $\mathbb{E}^{\bf\Gamma}_4$, nor 
$\mathbb{E}^{\bf\Gamma}_5$.

\vfill\eject

\noindent (2) $\mathbb{E}^{\bf\Gamma}_2$ is not above the other relations in 
${\mathcal B}^{\bf\Gamma}$, because of the number of classes. Appealing to the complexity of the classes, we see that $\mathbb{E}^{\bf\Gamma}_2$ is not below the other relations in ${\mathcal B}^{\bf\Gamma}$.\bigskip

\noindent (3) $\mathbb{E}^{\bf\Gamma}_3$ is not above the other relations, because of the complexity of the classes. Let us prove that $\mathbb{E}^{\bf\Gamma}_3$ is not below $\mathbb{E}^{\bf\Gamma}_5$ if the rank of $\bf\Gamma$ is at least two. We argue by contradiction, which gives 
$$g\! =\! (g_0,g_1)\! :\!\mathbb{H}\!\rightarrow\!\mathbb{H}.$$ 
Pick $(\varepsilon ,\alpha )\!\in\! 2\!\times\!\mathbb{C}$. If 
$g(\varepsilon ,\alpha )\! =\! (\varepsilon_0,\gamma )$ with $\gamma\!\notin\!\mathbb{C}$, then 
$g(1\! -\!\varepsilon ,\alpha )\! =\! (\varepsilon_0,\delta )$ with $\delta\!\notin\!\mathbb{C}$. The continuity of $g$ gives $l\!\in\!\omega$ such that $g_0(\varepsilon',\beta )\! =\!\varepsilon_0$ if 
$(\varepsilon',\beta )\!\in\! 2\!\times\! N_{\alpha\vert l}$. Note that there is $s$ in $2^{<\omega}$ such that 
$g\big(0,(\alpha\vert l)\beta\big)\!\not=\! g\big(1,(\alpha\vert l)\beta\big)$ if $\beta\!\in\! N_s$, since otherwise there is $(\alpha\vert l)\beta\!\notin\!\mathbb{C}$ with 
$g\big(0,(\alpha\vert l)\beta\big)\! =\! g\big(1,(\alpha\vert l)\beta\big)$, which is absurd. Then the map 
$\delta\!\mapsto\!\big( g_1(0,\delta ) ,g_1(1,\delta )\big)$ reduces $\mathbb{C}$ to $(\neg\mathbb{C})^2$ on $N_{(\alpha\vert l)s}$, which contradicts Lemma \ref{true}. This shows that 
$g_1(\varepsilon ,\alpha )\!\in\!\mathbb{C}$. As the rank of $\bf\Gamma$ is at least two, $\mathbb{C}$ is dense, so that we may assume that there are $\alpha\!\in\!\mathbb{C}$, $\varepsilon_0\!\in\! 2$ and 
$\gamma\!\in\!\mathbb{C}$ with $g(0,\alpha )\! =\! (\varepsilon_0,\gamma )$ and 
$g(1,\alpha )\! =\! (1\! -\!\varepsilon_0,\gamma )$. The continuity of $g$ gives $l\!\in\!\omega$ and 
$G\! :\! N_{\alpha\vert l}\!\rightarrow\! 2^\omega$ continuous with 
${g(0,\beta )\! =\!\big(\varepsilon_0,G(\beta )\big)}$ and 
$g(1,\beta )\! =\!\big( 1\! -\!\varepsilon_0,G(\beta )\big)$ if $\beta\!\in\! N_{\alpha\vert l}$. Note that $G$ reduces $\mathbb{C}\cap N_{\alpha\vert l}$ to $\mathbb{C}$. As the set $\mathbb{C}\cap N_{\alpha\vert l}$ is not open, there are $\beta ,\beta'\!\in\! N_{\alpha\vert l}\!\setminus\!\mathbb{C}$ with 
$G(\beta )\!\not=\! G(\beta')$. Note that 
$\Big(\big(\varepsilon_0,G(\beta )\big) ,\big(\varepsilon_0,G(\beta')\big)\Big)\!\in\!\mathbb{E}^{\bf\Gamma}_5$ and $\big( (0,\beta ),(0,\beta')\big)\!\in\!\mathbb{E}^{\bf\Gamma}_3$, which is absurd.\bigskip

 This argument also shows that $\mathbb{E}^{\bf\Gamma}_3$ is not below $\mathbb{E}^{\bf\Gamma}_4$ if the rank of $\bf\Gamma$ is at least two.\bigskip 

\noindent (4)-(8) As in (3), $\mathbb{E}^{\bf\Gamma}_4$ is not below $\mathbb{E}^{\bf\Gamma}_5$, and $\mathbb{E}^{\bf\Gamma}_5$ is not below 
$\mathbb{E}^{\bf\Gamma}_4$ since $\mathbb{E}^{\check {\bf\Gamma}}_0$ is not below $\Delta (2^\omega )$.
\hfill{$\square$}

\section{$\!\!\!\!\!\!$ Non-$\boraone$ equivalence relations}\indent 

 A strong form of Theorem \ref{main1} holds.\bigskip

\noindent\bf Theorem\it\label{open}\ Let $X$ be a metrizable topological space, and $E$ be an equivalence  relation on $X$. Then exactly one of the following holds:\smallskip  

(a) the equivalence classes of $E$ are $\boraone$ (exactly when $E$ is a $\boraone$ subset of $X^2$),\smallskip  

(b) there is $(\mathbb{X},\mathbb{E})\!\in\! {\mathcal A}^{\boraone}$ such that $(\mathbb{X},\mathbb{E})\sqsubseteq_c(X,E)$.\smallskip

\noindent Moreover, ${\mathcal A}^{\boraone}$ is a $\leq_c$-antichain (and thus a 
$\sqsubseteq_c$ and a $\leq_c$-antichain basis).\rm\bigskip 

\noindent\bf Proof.\rm\ By Lemma \ref{basic}.(a), the equivalence classes of $E$ are $\boraone$ if $E$ is an open subset of $X^2$. The converse comes from the fact that $E$ is the union of the square of its equivalence classes. By Theorem \ref{prope}.(a), (a) and (b) cannot hold simultaneously. So assume that (a) does not hold, which gives $x\!\in\! X$ such that 
$x\!\in\!\overline{\neg [x]_E}$.\bigskip

\noindent\bf Case 1\rm\ $x\!\notin\!\overline{\cal C}$ if $\cal C$ is an $E$-class which does not contain $x$.\bigskip

 We inductively construct an injective sequence $(x_k)_{k\in\omega}$ of points of 
$X\!\setminus\! [x]_E$ as follows. We first choose $x_0\!\in\! X\!\setminus\! [x]_E$. As 
$x\!\notin\!\overline{[x_0]_E}$, we choose 
$x_1\!\in\! B(x,2^{-1})\!\setminus\! ([x]_E\cup\overline{[x_0]_E})$. Then we choose 
$x_2\!\in\! B(x,2^{-2})\!\setminus\! ([x]_E\cup\overline{[x_0]_E}\cup\overline{[x_1]_E})$, and so on. Note that $(x_k)_{k\in\omega}$ converges to $x$. We define $f\! :\!\mathbb{K}\!\rightarrow\! X$ by setting $f(0)\! :=\! x$ and $f(2^{-k})\! :=\! x_k$. Note that $f$ is injective continuous and reduces 
$\mathbb{E}^{\boraone}_0$ to $E$.

\vfill\eject
 
\noindent\bf Case 2\rm\ There is an $E$-class $\cal C$ with $x\!\in\!\overline{\cal C}\!\setminus\! {\cal C}$.\bigskip

 As $X$ is metrizable, there is $(x_k)_{k\in\omega}$ injective in $\cal C$ converging to 
$x$.  We define $f\! :\!\mathbb{K}\!\rightarrow\! X$ by setting $f(0)\! :=\! x$ and $f(2^{-k})\! :=\! x_k$. Note that $f$ is injective continuous and reduces $\mathbb{E}^{\boraone}_1$ to $E$.\hfill{$\square$}\bigskip

\noindent\bf Remark.\rm\ This result does not hold for arbitrary relations, not even for linear quasi-orders. Indeed, assume that $(\mathbb{K},\mathbb{E}^{\boraone}_n)\leq_c(\mathbb{K},Q)$, where $Q$ is a non-$\boraone$ linear quasi-order on $\mathbb{K}$ like 
$$\{ (x,y)\!\in\!\mathbb{K}^2\mid x\!\leq\! y\} .$$
Pick $(x,y)\!\in\!\mathbb{K}^2\!\setminus\!\mathbb{E}^{\boraone}_n$. Then 
$\big( f(x),f(y)\big)\!\notin\! Q$, so that $\big( f(y),f(x)\big)\!\in\! Q$ and 
$(y,x)\!\in\!\mathbb{E}^{\boraone}_n$, which contradicts the symmetry of 
$\mathbb{E}^{\boraone}_n$.

\section{$\!\!\!\!\!\!$ Non-$\bormone$ equivalence relations}\indent

 A strong form of Theorem \ref{main1} holds.

\begin{thm} \label{closed} Let $X$ be a metrizable topological space, and $E$ be an equivalence  relation on $X$. Then exactly one of the following holds:\smallskip  

(a) the equivalence classes of $E$ are $\bormone$,\smallskip  

(b) $(\mathbb{K},\mathbb{E}^{\bormone}_0)\sqsubseteq_c(X,E)$.\end{thm}

\noindent\bf Proof.\rm\ By Lemma \ref{basic}.(b), (a) and (b) cannot hold simultaneously since 
$[1]_{\mathbb{E}^{\bormone}_0}\! =\!\mathbb{K}\!\setminus\!\{ 0\}$ is not closed. So assume that (a) does not hold, which gives $x\!\in\! X$ such that $[x]_E$ is not closed. Pick 
$y\!\in\!\overline{[x]_E}\!\setminus\! [x]_E$. As $X$ is metrizable, there is an injective sequence 
$(x_k)_{k\in\omega}$ in $[x]_E$ converging to $y$. We define $f\! :\!\mathbb{K}\!\rightarrow\! X$ by setting $f(0)\! :=\! y$ and $f(2^{-k})\! :=\! x_k$. Note that $f$ is injective continuous and reduces 
$\mathbb{E}^{\bormone}_0$ to $E$.\hfill{$\square$}\bigskip

 A strong form of Theorem \ref{main2} holds.
 
\begin{thm} \label{closed2} Let $X$ be a metrizable topological space, and $E$ be an equivalence  relation on $X$. Then exactly one of the following holds:\smallskip  

(a) $E$ is a $\bormone$ subset of $X^2$,\smallskip  

(b) there is $(\mathbb{X},\mathbb{E})\!\in\!\{ (\mathbb{K},\mathbb{E}^{\bormone}_0),
(\mathbb{H},\mathbb{E}^{\bormone}_3)\}$ such that $(\mathbb{X},\mathbb{E})\sqsubseteq_c(X,E)$.\smallskip

\noindent Moreover, $\{ (\mathbb{K},\mathbb{E}^{\bormone}_0),
(\mathbb{H},\mathbb{E}^{\bormone}_3)\}$ is a $\leq_c$-antichain (and thus a 
$\sqsubseteq_c$ and a $\leq_c$-antichain basis).\end{thm} 

\noindent\bf Proof.\rm\ By Theorem \ref{prope}.(a), (a) and (b) cannot hold simultaneously. So assume that (a) does not hold, which gives 
$(x,y)\!\in\!\overline{E}\!\setminus\! E$, and $\big( (x_k,y_k)\big)_{k\in\omega}\!\in\! E^\omega$ converging to $(x,y)$. Note that 
$x\!\not=\! y$, so that we may assume that $\overline{\{ x_k\mid k\!\in\!\omega\}}\cap\overline{\{ y_k\mid k\!\in\!\omega\}}\! =\!\emptyset$. We may also assume that either $x_k\! =\! x$ for each $k\!\in\!\omega$, or $(x_k)_{k\in\omega}$ is injective and $x_k\!\not=\! x$ for each $k\!\in\!\omega$. Moreover, we cannot have 
$(x_k,y_k)\! =\! (x,y)$ for each $k\!\in\!\omega$.\bigskip

\noindent\bf Case 1\rm\ $x_k\! =\! x$ and $y_k\!\not=\! y$ for each $k\!\in\!\omega$, and 
$(y_k)_{k\in\omega}$ is injective.\bigskip

 We define $f\! :\!\mathbb{K}\!\rightarrow\! X$ by setting $f(0)\! :=\! y$, $f(1)\! :=\! x$ and 
$f(2^{-k-1})\! :=\! y_k$. Note that $f$ is injective continuous and reduces 
$\mathbb{E}^{\bormone}_0$ to $E$.\bigskip

\noindent\bf Case 2\rm\ $y_k\! =\! y$ and $x_k\!\not=\! x$ for each $k\!\in\!\omega$, and 
$(x_k)_{k\in\omega}$ is injective.\bigskip

 We argue as in Case 1.\bigskip

\noindent\bf Case 3\rm\ $x_k\!\not=\! x$ and $y_k\!\not=\! y$ for each $k\!\in\!\omega$, and 
$(x_k)_{k\in\omega},(y_k)_{k\in\omega}$ are injective.\bigskip

 Note that we may assume that either $(x,x_k)\!\in\! E$ for each $k\!\in\!\omega$, or 
$(x,x_k)\!\notin\! E$ for each $k\!\in\!\omega$.\bigskip

\noindent\bf Case 3.1\rm\ $(x,x_k)\!\in\! E$ for each $k\!\in\!\omega$.\bigskip

 Note that $x,x_k,y_l$ are in the same $E$-class, which does not contain $y$. We define 
$f\! :\!\mathbb{K}\!\rightarrow\! X$ by setting $f(0)\! :=\! y$ and $f(2^{-k})\! :=\! y_k$. Note that $f$ is injective continuous and reduces $\mathbb{E}^{\bormone}_0$ to $E$.\bigskip

\noindent\bf Case 3.2\rm\ $(x,x_k)\!\notin\! E$ for each $k\!\in\!\omega$.\bigskip

 The previous discussion shows that we may assume that $(x,y_k),(y,x_k),(y,y_k)\!\notin\! E$ for each $k\!\in\!\omega$. By Ramsey's theorem (see 19.A in [K1]), we may assume that either 
$(x_k,x_l)\!\in\! E$ for each $k\!\not=\! l$, or $(x_k,x_l)\!\notin\! E$ for each $k\!\not=\! l$.\bigskip

\noindent\bf Case 3.2.1\rm\ $(x_k,x_l)\!\in\! E$ for each $k\!\not=\! l$.\bigskip

 We argue as in Case 3.1.\bigskip
  
\noindent\bf Case 3.2.2\rm\ $(x_k,x_l)\!\notin\! E$ for each $k\!\not=\! l$.\bigskip 

 The previous discussion shows that we may assume that $(x_k,y_l),(y_k,y_l)\!\notin\! E$ for each $k\!\not=\! l$. We define $f\! :\!\mathbb{H}\!\rightarrow\! X$ by setting $f(0,0)\! :=\! x$, $f(1,0)\! :=\! y$, $f(0,2^{-k})\! :=\! x_k$ and $f(1,2^{-k})\! :=\! y_k$. Note that $f$ is injective continuous and reduces $\mathbb{E}^{\bormone}_3$ to $E$.
\hfill{$\square$}
 
\section{$\!\!\!\!\!\!$ Some facts about the rank two} 
  
\begin{lem} \label{subrat} Let $D$ be a non-nowhere dense subset of $2^\omega$ contained in  
$\mathbb{P}_f$. Then there is $f\! :\! 2^\omega\!\rightarrow\! 2^\omega$ injective continuous such that $f[\mathbb{P}_f]\!\subseteq\! D$ and $f[\mathbb{P}_\infty]\!\subseteq\!\mathbb{P}_\infty$.
\end{lem}
 
\noindent\bf Proof.\rm\ Let $s\!\in\! 2^{<\omega}$ such that $N_s\!\subseteq\!\overline{D}$. Note that $N_s\!\subseteq\!\overline{N_s\cap D}$, so that $N_s\cap D$ is dense and co-dense in 
$N_s$. In particular, by Baire's theorem, $N_s\cap D$ is not separable from 
$N_s\cap\mathbb{P}_\infty$ by a $\bormtwo$ set. Theorem \ref{Lo-SR} gives 
$f\! :\! 2^\omega\!\rightarrow\! 2^\omega$ injective continuous such that 
$f[\mathbb{P}_f]\!\subseteq\! N_s\cap D$ and 
$f[\mathbb{P}_\infty]\!\subseteq\! N_s\cap\mathbb{P}_\infty$.\hfill{$\square$}

\begin{lem} \label{subirrat} Let $G$ be a non-meager subset of $2^\omega$ having the Baire property and contained in $\mathbb{P}_\infty$. Then there is 
$f\! :\! 2^\omega\!\rightarrow\! 2^\omega$ injective continuous such that 
$f[\mathbb{P}_\infty ]\!\subseteq\! G$ and $f[\mathbb{P}_f]\!\subseteq\!\mathbb{P}_f$.\end{lem}
 
\noindent\bf Proof.\rm\ As $G$ has the Baire property and is not meager, there is 
$s\!\in\! 2^{<\omega}$ such that $N_s\cap G$ is comeager in $N_s$. By Baire's theorem, 
$N_s\cap\mathbb{P}_f$ is not separable from $N_s\cap G$ by a $\bormtwo$ set. Theorem 
\ref{Lo-SR} gives $f\! :\! 2^\omega\!\rightarrow\! 2^\omega$ injective continuous such that 
$f[\mathbb{P}_f]\!\subseteq\! N_s\cap\mathbb{P}_f$ and 
$f[\mathbb{P}_\infty]\!\subseteq\! N_s\cap G$.\hfill{$\square$}\bigskip

\noindent\bf Convention.\rm\ In the rest of Sections 7 to 9, we will perform a number of Cantor-like constructions. The following will always hold. We fix $s\!\in\! 2^{<\omega}$, and inductively construct a sequence $(n_t)_{t\in 2^{<\omega}}$ of natural numbers, and a sequence 
$(U_t)_{t\in 2^{<\omega}}$ of clopen subsets of $2^\omega$ satisfying the following conditions:
$$\begin{array}{ll}
& (1)~U_{t\varepsilon}\!\subseteq\! U_t\!\subseteq\! N_s\cr
& (2)~\alpha_{n_t}\!\in\! U_t\cr
& (3)~\mbox{diam}(U_t)\!\leq\! 2^{-\vert t\vert}\cr
& (4)~U_{t0}\cap U_{t1}\! =\!\emptyset\cr
& (5)~n_{t0}\! =\! n_t
\end{array}$$
Assume that this is done. Using (1)-(3), we define $f\! :\! 2^\omega\!\rightarrow\! N_s$ by 
$\{ f(\beta )\}\! :=\!\bigcap_{n\in\omega}~U_{\beta\vert n}$, and $f$ is injective continuous by (4). If 
$t\!\in\! O$ and $\alpha\! =\! t0^\infty$, then $f(\alpha )\! =\!\alpha_{n_t}$ by (5). For the first step of the induction, we choose $n_\emptyset$ in such a way that $\alpha_{n_\emptyset}\!\in\! N_s$ and set $U_\emptyset\! :=\! N_s$. Condition (5) defines $n_{t0}$.

\begin{lem} \label{lemGdelta} Let $b\! :\!\mathbb{P}_\infty\!\rightarrow\! 2^\omega$ be a nowhere dense-to-one continuous map. Then there is ${f\! :\! 2^\omega\!\rightarrow\! 2^\omega}$  injective continuous such that $\mathbb{P}_f\! =\! f^{-1}(\mathbb{P}_f)$ and 
$b\big( f(\alpha )\big)\!\not=\! b\big( f(\beta )\big)$ if $\alpha\!\not=\!\beta\!\in\!\mathbb{P}_\infty$.
\end{lem}

\noindent\bf Proof.\rm\ We first prove the following.\bigskip

\noindent\bf Claim\it\ Let $\beta\!\in\!\mathbb{P}_f$. Then there is a sequence 
$(s^\beta_q)_{q\in\omega}$ of finite binary sequences such that\smallskip

(a) $\vert s^\beta_q\vert\! >\! q$,\smallskip

(b) $s^\beta_q\vert q\! =\!\beta\vert q$,\smallskip

(c) $s^\beta_q\!\not\subseteq\!\beta$,\smallskip

(d) $\forall p\!\not=\! q~~b[N_{s^\beta_p}\cap\mathbb{P}_\infty ]\cap 
b[N_{s^\beta_q}\cap\mathbb{P}_\infty ]\! =\!\emptyset$.\rm\bigskip

 Indeed, we first construct a sequence $(\beta_n)_{n\in\omega}$ of elements of 
$\mathbb{P}_\infty$ converging to $\beta$ and such that $\big( b(\beta_n)\big)_{n\in\omega}$ is injective. Assume that $(\beta_n)_{n\leq l}$ have been constructed. As $b$ is nowhere dense-to-one, we can find $\beta_{l+1}\!\in\! N_{\beta\vert (l+1)}\cap
\mathbb{P}_\infty\!\setminus\!\Big(\bigcup_{n\leq l}~b^{-1}\big(\{ b(\beta_n)\}\big)\Big)$.\bigskip

 We can extract a subsequence if necessary to ensure that $\big( b(\beta_n)\big)_{n\in\omega}$ converges to some $\gamma\!\in\! 2^\omega$, which is compact. Extracting again if necessary, we may assume that $b(\beta_n)\!\not=\!\gamma$ for each $n\!\in\!\omega$. As 
$b(\beta_0)\!\not=\!\gamma$ and $b$ is continuous, we can find $n_0\!\in\!\omega$ and 
$l_0\! >\! 0$ such that $b(\beta_0)(n_0)\!\not=\!\gamma (n_0)$, 
$\beta_0\vert l_0\!\not=\!\beta\vert l_0$ and 
$b(\alpha )\vert (n_0\! +\! 1)\! =\! b(\beta_0)\vert (n_0\! +\! 1)$ if 
$\alpha\!\in\! N_{\beta_0\vert l_0}\cap\mathbb{P}_\infty$. We set $s^\beta_0\! :=\!\beta_0\vert l_0$.\bigskip 

 Extracting again if necessary, we may assume that $\beta_n\vert 1\! =\!\beta\vert 1$ and 
$b(\beta_n)\vert (n_0\! +\! 1)\! =\!\gamma\vert (n_0\! +\! 1)$ for each $n\! >\! 0$. As 
$b(\beta_1)\!\not=\!\gamma$ and $b$ is continuous, we can find $n_1\! >\! n_0$ and 
$l_1\! >\! l_0$ such that $b(\beta_1)(n_1)\!\not=\!\gamma (n_1)$, 
$\beta_1\vert l_1\!\not=\!\beta\vert l_1$ and 
$b(\alpha )\vert (n_1\! +\! 1)\! =\! b(\beta_1)\vert (n_1\! +\! 1)$ if 
$\alpha\!\in\! N_{\beta_1\vert l_1}\cap\mathbb{P}_\infty$. We set $s^\beta_1\! :=\!\beta_1\vert l_1$. Note that $b(\alpha )(n_0)\!\not=\!\gamma (n_0)$ if 
$\alpha\!\in\! N_{s^\beta_0}\cap\mathbb{P}_\infty$, and $b(\alpha )(n_0)\! =\!\gamma (n_0)$ if 
$\alpha\!\in\! N_{s^\beta_1}\cap\mathbb{P}_\infty$, so that 
$b[N_{s^\beta_0}\cap\mathbb{P}_\infty ]\cap b[N_{s^\beta_1}\cap\mathbb{P}_\infty ]\! =\!
\emptyset$. We just have to continue like this to finish the construction of the desired $s^\beta_q$'s.
\hfill{$\diamond$}

\vfill\eject

 We set $s\! :=\!\emptyset$, and construct $(n_t)_{t\in 2^{<\omega}}$, 
$(U_t)_{t\in 2^{<\omega}}$, $(q^u_m)_{u\in O,m\in\omega}$ satisfying (1)-(5) and 
$$\begin{array}{ll}
& (6)~q^u_m\! <\! q^u_{m+1}\cr
& (7)~U_{t1}\cap\{\alpha_n\mid n\!\leq\!\vert t\vert\}\! =\!\emptyset\cr
& (8)~U_{u0^m1}\!\subseteq\! N_{s^{\alpha_{n_u}}_{q^u_m}}
\end{array}$$
Assume that this is done. If $\beta\!\in\!\mathbb{P}_\infty$, then there is an infinite strictly increasing sequence $(l_k)_{k\in\omega}$ of natural numbers with $\beta\vert l_k\!\in\! O$. Condition (7) implies that $f(\beta )\!\in\!\mathbb{P}_\infty$. Let 
$\alpha\!\not=\!\beta\!\in\!\mathbb{P}_\infty$, which gives $u\!\in\! O$ and $m\!\not=\! p$ such that 
$\alpha\!\in\! N_{u0^m1}$ and $\beta\!\in\! N_{u0^p1}$. Conditions (8) and (d) in the claim imply that $b\big( f(\alpha )\big)\!\not=\! b\big( f(\beta )\big)$. So we are done.\bigskip

 Let us prove that the construction is possible. Assume that 
$(n_t)_{\vert t\vert\leq l}$, $(U_t)_{\vert t\vert\leq l}$, $(q^u_m)_{u\in O,\vert u\vert +m+1\leq l}$ satisfying (1)-(8) have been constructed, which is the case for $l\! =\! 0$.\bigskip

 Let $t\! :=\! u0^m\!\in\! 2^l$, with $u\!\in\! O$. As $\alpha_{n_t}\!\in\! U_t$, 
$N_{\alpha_{n_t}\vert q}\!\subseteq\! U_t$ if $q$ is big enough, say $q\!\geq\! q_t$. We choose 
$q^u_m\! >\!\mbox{max}(\mbox{max}_{j<m}~q^u_j,q_t)$, and $n_{t1}$ in such a way that
$\alpha_{n_{t1}}\!\in\! N_{s^{\alpha_{n_t}}_{q^u_m}}\!\setminus\!\{\alpha_n\mid n\!\leq\! l\}$. Note that 
$s^{\alpha_{n_t}}_{q^u_m}\vert q^u_m\! =\!\alpha_{n_t}\vert q^u_m\!\supseteq\!\alpha_{n_t}\vert q_t$, so that 
$\alpha_{n_{t1}}\!\in\! U_t$ and $\alpha_{n_{t1}}\!\not=\!\alpha_{n_{t0}}$. We choose disjoint clopen sets $U_{t0},U_{t1}$ with diameter at most $2^{-l-1}$ such that 
$\alpha_{n_{t\varepsilon}}\!\in\! U_{t\varepsilon}\!\subseteq\! U_t$ and satisfying (7)-(8).
\hfill{$\square$}
 
\section{$\!\!\!\!\!\!$ Non-$\boratwo$ equivalence relations}

\bf Notation.\rm\ We set $\mathbb{C}\! :=\!\mathbb{P}_\infty$.\bigskip

\noindent\bf Proof of Theorem \ref{main1} when ${\bf\Gamma}\! =\!\boratwo$.\rm\ By Lemma 
\ref{basic}.(b), (a) and (b) cannot hold simultaneously. So assume that (a) does not hold. By Proposition \ref{injcontA}, we may assume that $X\! =\! 2^\omega$ and $\mathbb{C}$ is an equivalence class of $E$.\bigskip

\noindent\bf Case 1\rm\ $[\alpha ]_E$ is nowhere dense for each $\alpha\!\notin\!\mathbb{C}$.\bigskip

 We inductively construct a sequence $(n_k)_{k\in\omega}$ of natural numbers as follows. Let 
$(O_k)_{k\in\omega}$ be a basis for the topology of $2^\omega$ made of nonempty sets. Pick $n_0\!\in\!\omega$ such that $\alpha_{n_0}\!\in\! O_0$. As $[\alpha_{n_0}]_E$ is nowhere dense, we can find $n_1\!\in\!\omega$ such that $\alpha_{n_1}\!\in\! O_1\!\setminus\!\overline{[\alpha_{n_0}]_E}$. As $[\alpha_{n_1}]_E$ is nowhere dense, we can find $n_2\!\in\!\omega$ such that 
$\alpha_{n_2}\!\in\! O_2\!\setminus\! (\overline{[\alpha_{n_0}]_E}\cup\overline{[\alpha_{n_1}]_E})$. And so on. Note that $(\alpha_{n_k})_{k\in\omega}$ is dense and co-dense in $\mathbb{C}\cup\{\alpha_{n_k}\mid k\!\in\!\omega\}$ (which is co-countable in $2^\omega$), so that $\{\alpha_{n_k}\mid k\!\in\!\omega\}$ is not $\bormtwo$, by Baire's theorem. By Hurewicz's theorem, there is 
$f\! :\! 2^\omega\!\rightarrow\!\mathbb{C}\cup\{\alpha_{n_k}\mid k\!\in\!\omega\}$ injective continuous such that 
$\mathbb{C}\! =\! f^{-1}(\mathbb{C})$. Note that $f$ reduces $\mathbb{E}^{\boratwo}_0$ to $E$.\bigskip

\noindent\bf Case 2\rm\ there is $\alpha\!\notin\!\mathbb{C}$ such that $[\alpha ]_E$ is not nowhere dense.\bigskip

 Let $s\!\in\! 2^{<\omega}$ such that $N_s\!\subseteq\!\overline{[\alpha ]_E}$. Note that the countable and thus 
$\boratwo$ set $N_s\cap [\alpha ]_E$ is dense and co-dense in $N_s\cap (\mathbb{C}\cup [\alpha ]_E)$ (which is co-countable in $N_s$). As in the Case 1, we get 
${f\! :\! 2^\omega\!\rightarrow\! N_s\cap (\mathbb{C}\cup [\alpha ]_E)}$ injective continuous such that 
$\mathbb{C}\! =\! f^{-1}(\mathbb{C})$. Note that $f$ reduces $\mathbb{E}^{\boratwo}_1$ to $E$. This finishes the proof.\hfill{$\square$}\bigskip

\noindent\bf Proof of Theorem \ref{main2} when ${\bf\Gamma}\! =\!\boratwo$.\rm\ If 
$(\mathbb{X},\mathbb{E})\!\in\! {\mathcal B}^{\boratwo}$, then 
$\mathbb{E}\!\notin\!\boratwo$ by Theorem \ref{prope}.(a), so that (a) and (b) cannot hold simultaneously. Assume that (a) does not hold. By Theorem \ref{main1}, we may assume that the equivalence classes of $E$ are 
$\boratwo$. By Theorem \ref{contB}, we may assume that $X\! =\! 2\!\times\! 2^\omega$, 
$$\Delta (2\!\times\! 2^\omega )\cup\big\{\big( (0,\alpha ),(1,\alpha )\big)\mid\alpha\!\in\!\mathbb{C}\big\}
\!\subseteq\! E$$ 
and $\big\{\big( (0,\alpha ),(1,\alpha )\big)\mid\alpha\!\notin\!\mathbb{C}\big\}\!\subseteq\!\neg E$.\bigskip

 We will now prove that we may assume that 
$\big( (\varepsilon ,\alpha ),(1\! -\!\varepsilon ,\beta )\big)\!\notin\! E$ if $\varepsilon\!\in\! 2$ and 
$\alpha ,\beta\!\notin\!\mathbb{C}$. Indeed, assume first that 
$(E_{(\varepsilon ,\alpha )})_{1-\varepsilon}\!\setminus\!\mathbb{C}$ is not nowhere dense in $2^\omega$ for some 
$\varepsilon\!\in\! 2$ and some $\alpha\!\notin\!\mathbb{C}$. Then 
$\big( (1\! -\!\varepsilon ,\beta ),(1\! -\!\varepsilon ,\gamma )\big)\!\in\! E$, $\big( (0,\beta ),(1,\gamma )\big)\!\notin\! E$ and $\big( (1,\beta ),(0,\gamma )\big)\!\notin\! E$ if 
$\beta ,\gamma\!\in\! (E_{(\varepsilon ,\alpha )})_{1-\varepsilon}\!\setminus\!\mathbb{C}$. Lemma \ref{subrat} gives 
$f\! :\! 2^\omega\!\rightarrow\! 2^\omega$ injective continuous such that 
$f[2^\omega\!\setminus\!\mathbb{C}]\!\subseteq\! (E_{(\varepsilon,\alpha )})_{1-\varepsilon}\!\setminus\!\mathbb{C}$ and $f[\mathbb{C}]\!\subseteq\!\mathbb{C}$, so we are done. Assume now that 
$(E_{(\varepsilon ,\alpha )})_{1-\varepsilon}\!\setminus\!\mathbb{C}$ is nowhere dense in 
$2^\omega$ for each $\varepsilon\!\in\! 2$ and each $\alpha\!\notin\!\mathbb{C}$.  We set 
$s\! :=\!\emptyset$, and construct $(n_t)_{t\in 2^{<\omega}}$, $(U_t)_{t\in 2^{<\omega}}$ satisfying (1)-(5) and the following:
$$(6)~\alpha_{n_{t1}}\!\notin\!\{\alpha_{n_t}\}\cup\{\alpha_n\mid n\!\leq\!\vert t\vert\}\cup
\bigcup_{\varepsilon\in 2,s\in 2^{\vert t\vert}}~(E_{(\varepsilon ,\alpha_{n_s})})_{1-\varepsilon}\cup
\bigcup_{\varepsilon\in 2,s\in 2^{\vert t\vert},s<_{\mbox{lex}}t}~
(E_{(\varepsilon ,\alpha_{n_{s1}})})_{1-\varepsilon}$$
Assume that this is done. If $\beta\!\in\!\mathbb{C}$, then there is an infinite strictly increasing sequence 
$(l_k)_{k\in\omega}$ of natural numbers with $\beta\vert l_k\!\in\! O$. Condition (6) implies that 
$f(\beta )\!\in\!\mathbb{C}$. Now let $\beta\!\not=\!\beta'\!\notin\!\mathbb{C}$. Condition (6) implies that 
$\Big(\big(\varepsilon ,f(\beta )\big) ,\big( 1\! -\!\varepsilon ,f(\beta')\big)\Big)\!\notin\! E$ for each 
$\varepsilon\!\in\! 2$. So we are done. Let us prove that the construction is possible. Assume that 
$(n_t)_{\vert t\vert\leq l}$ and $(U_t)_{\vert t\vert\leq l}$ satisfying (1)-(6) have been constructed, which is the case for $l\! =\! 0$. Let $t\!\in\! 2^l$. We define $n_{t1}$ by induction on $t$ with respect to the lexicographical ordering. We choose it in such a way that 
$$\alpha_{n_{t1}}\!\in\! U_t\!\setminus\! 
\Big(\{\alpha_{n_t}\}\cup\{\alpha_n\mid n\!\leq\! l\}\cup
\bigcup_{\varepsilon\in 2,s\in 2^l}~
\overline{(E_{(\varepsilon,\alpha_{n_s})})_{1-\varepsilon}\!\setminus\!\mathbb{C}}\cup
\bigcup_{\varepsilon\in 2,s\in 2^l,s<_{\mbox{lex}}t}~
\overline{(E_{(\varepsilon ,\alpha_{n_{s1}})})_{1-\varepsilon}\!\setminus\!\mathbb{C}}\Big) .$$
We do this for each $t\!\in\! 2^l$, in the lexicographical ordering. Then we choose disjoint clopen subsets $U_{t0},U_{t1}$ of $U_t$ with diameter at most $2^{-l-1}$ with 
$\alpha_{n_{t\varepsilon}}\!\in\! U_{t\varepsilon}$ for each $\varepsilon\!\in\! 2$.\bigskip

 Similarly, we may assume that either 
$\big( (\varepsilon ,\alpha ),(\varepsilon ,\beta )\big)\!\in\! E$ for each 
$\alpha\!\not=\!\beta\!\notin\!\mathbb{C}$, or $\big( (\varepsilon ,\alpha ),(\varepsilon ,\beta )\big)$ is not in $E$ for each $\alpha\!\not=\!\beta\!\notin\!\mathbb{C}$, for each $\varepsilon\!\in\! 2$.\bigskip

 Let us prove that $E$ has meager classes. We argue by contradiction, which gives 
$(\varepsilon ,\alpha )\!\in\! 2\!\times\! 2^\omega$ such that $[(\varepsilon ,\alpha )]_E$ is not meager. As 
$[(\varepsilon ,\alpha )]_E$ is in ${\bf\Gamma}\! =\!\boratwo$, we get $\varepsilon'\!\in\! 2$ and $s\!\in\! 2^{<\omega}$ such that $\{\varepsilon'\}\!\times\! N_s\!\subseteq\! [(\varepsilon ,\alpha )]_E$. Assume, for example, that 
$\varepsilon'\! =\! 0$, so that $\{ 1\}\!\times\! (N_s\cap\mathbb{C})\!\subseteq\! [(\varepsilon ,\alpha )]_E$. Thus $(\{ 1\}\!\times\! N_s)\cap [(\varepsilon ,\alpha )]_E$ is comeager in $\{ 1\}\!\times\! N_s$ and 
$\boratwo$, which gives $t\!\in\! 2^{<\omega}$ such that 
$$\{ 1\}\!\times\! N_{st}\!\subseteq\! [(\varepsilon ,\alpha )]_E.$$ 
Thus $(0,st0^\infty ),(1,st0^\infty )\!\in\! [(\varepsilon ,\alpha )]_E$ and 
$\big( (0,st0^\infty ),(1,st0^\infty )\big)\!\in\! E$, which is absurd.\bigskip

 The Sarbadhikari theorem gives an increasing sequence $(F_l)_{l\in\omega}$ of Borel relations on $2\!\times\! 2^\omega$ with closed nowhere dense vertical sections whose union contains $E$ (see 5.12.11 in [Sr]).\bigskip

 We will now prove that we may assume that 
$\big( (\varepsilon ,\alpha ),(\varepsilon',\beta )\big)\!\notin\! E$ if $\varepsilon ,\varepsilon'\!\in\! 2$, and either 
$\alpha\!\in\!\mathbb{C}$ and $\beta\!\notin\!\mathbb{C}$, or $\alpha\!\notin\!\mathbb{C}$ and 
$\beta\!\in\!\mathbb{C}$. We set $s\! :=\!\emptyset$, and construct $(n_t)_{t\in 2^{<\omega}}$, 
$(U_t)_{t\in 2^{<\omega}}$ satisfying (1)-(5) and the following:
$$\begin{array}{ll}
& (6)~U_{t1}\cap\Big(\{\alpha_{n_t}\}\cup\{\alpha_n\mid n\!\leq\!\vert t\vert\}\cup
\bigcup_{\varepsilon ,\varepsilon'\in 2,s\in 2^{\vert t\vert}}~
\big( (F_{\vert t\vert})_{(\varepsilon ,\alpha_{n_s})}\big)_{\varepsilon'}\Big)\! =\!\emptyset
\end{array}$$
Assume that this is done. If $\beta\!\in\!\mathbb{C}$, then there is an infinite strictly increasing sequence 
$(l_k)_{k\in\omega}$ of natural numbers with $\beta\vert l_k\!\in\! O$. Condition (6) implies that 
$f(\beta )\!\in\!\mathbb{C}$. Condition (6) also implies that 
$\Big(\big(\varepsilon ,f(\gamma )\big) ,\big(\varepsilon',f(\delta )\big)\Big)\!\notin\! E$ if 
$\varepsilon ,\varepsilon'\!\in\! 2$, and either $\gamma\!\in\!\mathbb{C}$ and $\delta\!\notin\!\mathbb{C}$, or 
$\gamma\!\notin\!\mathbb{C}$ and $\delta\!\in\!\mathbb{C}$. So we are done. Let us prove that the construction is possible. Assume that $(n_t)_{\vert t\vert\leq l}$ and $(U_t)_{\vert t\vert\leq l}$ satisfying (1)-(6) have been constructed, which is the case for $l\! =\! 0$. Let $t\!\in\! 2^l$. As $(F_l)_{(\varepsilon ,\alpha )}$ is closed nowhere dense for each $(\varepsilon ,\alpha )\!\in\! 2\!\times\! 2^\omega$, 
$\big( (F_l)_{(\varepsilon ,\alpha )}\big)_{\varepsilon'}$ is a nowhere dense closed subset of 
$2^\omega$ for each $(\varepsilon ,\varepsilon',\alpha )$ in $2^2\!\times\! 2^\omega$. We choose 
$n_{t1}$ in such a way that
$$\alpha_{n_{t1}}\!\in\! U_t\!\setminus\!
\Big(\{\alpha_{n_t}\}\cup\{\alpha_n\mid n\!\leq\! l\}\cup
\bigcup_{\varepsilon ,\varepsilon'\in 2,s\in 2^l}~
\big( (F_l)_{(\varepsilon ,\alpha_{n_s})}\big)_{\varepsilon'}\Big) .$$
Then we choose disjoint clopen sets $U_{t0},U_{t1}$ with diameter at most $2^{-l-1}$ such that 
$\alpha_{n_{t\varepsilon}}\!\in\! U_{t\varepsilon}\!\subseteq\! U_t$ and satisfying (6).\bigskip

 We will now prove that we may assume that $E\cap (2\!\times\!\mathbb{C})^2\! =\!
\big\{\big( (\varepsilon ,\alpha ),(\varepsilon',\alpha )\big)\mid
\varepsilon ,\varepsilon'\!\in\! 2\wedge\alpha\!\in\!\mathbb{C}\big\}$. By Proposition \ref{foursame}, we just have to prove that we may assume that 
$$E\cap (\{ 0\}\!\times\!\mathbb{C})^2\! =\!\big\{\big( (0,\alpha ),(0,\alpha )\big)\mid\alpha\!\in\!\mathbb{C}\big\} .$$ 
We set 
$E'\! :=\!\big\{ (\alpha ,\beta )\!\in\! 2^\omega\!\times\! 2^\omega\mid\big( (0,\alpha ),(0,\beta )\big)\!\in\! E\big\}$, so that we must see that we may assume that $E'\cap\mathbb{C}^2\! =\!\Delta (\mathbb{C})$. As $E$ is a Borel equivalence relation on $2\!\times\! 2^\omega$ with $\boratwo$ classes, we can write $E'\! =\!\bigcup_{q\in\omega}~K_q$, where $K_q$ is a Borel relation on 
$2^\omega$ with nonempty closed vertical sections, by the Saint Raymond theorem (see 35.45 in [K1]). By Theorem 3.6 in [Lo1], there is for each $n\!\in\!\omega$ a finer Polish topology $\tau_q$ on $2^\omega$ such that $K_q\!\in\!\bormone\big( (2^\omega ,\tau_q)\!\times\! 2^\omega\big)$. By 8.38 in [K1], there is a dense $G_\delta$ subset $G_q$ of $2^\omega$ on which $\tau_q$ coincides with the usual topology on $2^\omega$, so that 
$K_q\cap (G_q\!\times\! 2^\omega )\!\in\!\bormone (G_q\!\times\! 2^\omega )$. We equip the hyperspace 
$F(2^\omega )$ of closed subsets of $2^\omega$ with the Effros Borel structure (see 12.C in [K1]). The following maps are Borel.\smallskip

\noindent (i) $\psi_q\! :\! 2^\omega\!\rightarrow\! F(2^\omega )\!\setminus\!\{\emptyset\}$ defined by 
$\psi_q(\gamma )\! :=\! (K_q)_\gamma$.\smallskip

 Indeed, 
$(K_q)_\gamma\cap N_s\!\not=\!\emptyset\ \Leftrightarrow\ \exists\beta\!\in\! N_s~~(\gamma ,\beta )\!\in\! K_q$, so that $\{\gamma\!\in\! 2^\omega\mid (K_q)_\gamma\cap N_s\!\not=\!\emptyset\}$ is analytic. Assume, for simplicity of the notation, that $K_q$ is $\Borel$. If $\gamma\!\in\! 2^\omega$, then $(K_q)_\gamma\cap N_s$ is $\Borel (\gamma )$ and compact. By 4F.11 in [Mos], 
$(K_q)_\gamma\cap N_s$ is not empty if and only if it contains a $\Borel (\gamma )$ point. This shows that $\{\gamma\!\in\! 2^\omega\mid (K_q)_\gamma\cap N_s\!\not=\!\emptyset\}$ is also co-analytic, and thus Borel. Thus $\psi_q$ is Borel.\smallskip

\noindent (ii) $\phi_\alpha\! :\! F(2^\omega )\!\setminus\!\{\emptyset\}\!\rightarrow\!\mathbb{R}$ defined by $\phi_\alpha (K)\! :=\! d(\alpha ,K)$.\smallskip

 By 12.13 in [K1], there is a sequence $(d_k)_{k\in\omega}$ of Borel functions from $F(2^\omega )$ into $2^\omega$ such that $\big( d_k(K)\big)_{k\in\omega}$ is dense in $K$ if $K\!\in\! F(2^\omega )$ is not empty. We get the following, for $a,b\!\in\!\mathbb{R}$:
$$d(\alpha ,K)\! >\! a\Leftrightarrow\exists p\!\in\!\omega ~~\forall k\!\in\!\omega ~~
d\big(\alpha ,d_k(K)\big)\! >\! a\! +\! 2^{-p}\mbox{,}$$
$$d(\alpha ,K)\! <\! b\Leftrightarrow\exists k\!\in\!\omega ~~
d\big(\alpha ,d_k(K)\big)\! <\! b\mbox{,}$$
showing that $\phi_\alpha$ is Borel.

\vfill\eject

\noindent (iii) $\varphi_{q,\alpha}\! :\! 2^\omega\!\rightarrow\!\mathbb{R}$ defined by 
$\varphi_{q,\alpha}(\gamma )\! :=\! d\big(\alpha ,(K_q)_\gamma\big)$.\smallskip

 Indeed, $\varphi_{q,\alpha}\! =\!\phi_\alpha\circ\psi_q$. Consequently, 8.38 in [K1] gives a dense 
$G_\delta$ subset $H_{q,\alpha}$ of $2^\omega$ on which $\varphi_{q,\alpha}$ is continuous. We set 
$H\! :=\!\mathbb{C}\cap\bigcap_{q\in\omega}~G_q\cap\bigcap_{q,n\in\omega}~H_{q,\alpha_n}$, so that $H$ is also a dense $G_\delta$ subset of $2^\omega$. In particular, 
$H\! =\!\bigcap_{l\in\omega}~O_l$, where $(O_l)_{l\in\omega}$ is a decreasing sequence of dense open subsets of $2^\omega$.\smallskip

 We set $s\! :=\!\emptyset$, and construct $(n_t)_{t\in 2^{<\omega}}$, 
$(U_t)_{t\in 2^{<\omega}}$ satisfying (1)-(5) and the following:
$$\begin{array}{ll}
& (6)~U_{t1}\!\subseteq\! O_{\vert t\vert}\!\setminus\!\{\alpha_n\mid n\!\leq\!\vert t\vert\}\cr
& (7)~\big(\bigcup_{s\in 2^{\vert t\vert},\eta\in 2,s\eta\not= t1}~
(U_{s\eta}\cap H)\!\times\! (U_{t1}\cap H)\big)\cap (\bigcup_{q\leq\vert t\vert +1}~K_q)\! =\!\emptyset
\end{array}$$
Assume that this is done. If $\beta\!\in\!\mathbb{C}$, then there is an infinite strictly increasing sequence 
$(l_k)_{k\in\omega}$ of natural numbers with $\beta\vert l_k\!\in\! O$. Condition (6) implies that 
$f(\beta )\!\in\!\mathbb{C}$. Conditions (6)-(7) imply that $\big( f(\gamma ),f(\delta )\big)\!\notin\! E'$ if $\gamma\!\not=\!\delta\!\in\!\mathbb{C}$. So we are done. Let us prove that the construction is possible. Assume that 
$(n_t)_{\vert t\vert\leq l}$ and $(U_t)_{\vert t\vert\leq l}$ satisfying (1)-(7) have been constructed, which is the case for $l\! =\! 0$. Note first that $E'$ is a meager relation on $2^\omega$ since $E$ has meager classes and is Borel. In particular, $E'\cap H^2$ is meager in $H^2$ since $H$ is a dense $G_\delta$ subset of $2^\omega$. Moreover, 
$\bigcup_{q\leq l+1}~K_q\cap H^2$ is a closed relation on $H$ contained in $E'$, so that 
$\bigcup_{q\leq l+1}~K_q\cap H^2$ is nowhere dense in $H^2$. Now let $s\!\not=\! t\!\in\! 2^l$ (we have $s1$ and $t1$ in mind). Note that $(U_s\cap H)\!\times\! (U_t\cap H)\!\not\subseteq\!\bigcup_{q\leq l+1}~K_q$. So we can find a nonempty clopen subset $V_s$ of $U_s\cap O_l\!\setminus\! (\{\alpha_{n_s}\}\cup\{\alpha_n\mid n\!\leq\! l\} )$ such that $\big( (V_s\cap H)\!\times\! (V_t\cap H)\big)\cap (\bigcup_{q\leq l+1}~K_q)\! =\!\emptyset$ if 
$s\!\not=\! t\!\in\! 2^l$. Now let $s,t\!\in\! 2^l$ (we have $s0$ and $t1$ in mind). We choose 
$\gamma_{s,t}\!\in\! V_t\cap H$. If $q\!\leq\! l\! +\! 1$, then 
${\alpha_{n_s}\!\notin\! (K_q)_{\gamma_{s,t}}}$ since $\alpha_{n_s}\!\notin\!\mathbb{C}$ and 
$(K_q)_{\gamma_{s,t}}\!\subseteq\! (E')_{\gamma_{s,t}}\!\subseteq\!\mathbb{C}$. As 
$(K_q)_{\gamma_{s,t}}$ is closed in $2^\omega$, it is compact. This gives $p_{s,t,q}\!\in\!\omega$ such that $d(\alpha_{n_s},(K_q)_{\gamma_{s,t}})\! >\! 2^{-p_{s,t,q}}$. The continuity of 
$\varphi_{q,\alpha_{n_s}}$ on $H$ gives $l_{s,t,q}\!\in\!\omega$ such that 
$d(\alpha ,(K_q)_\gamma )\! >\! 2^{-p_{s,t,q}}$ if $\alpha\!\in\! N_{\alpha_{n_s}\vert l_{s,t,q}}$ and 
$\gamma\!\in\! H\cap N_{\gamma_{s,t}\vert l_{s,t,q}}$, 
$N_{\gamma_{s,t}\vert l_{s,t,q}}\!\subseteq\! V_t$, 
$N_{\alpha_{n_s}\vert l_{s,t,q}}\!\subseteq\! U_s$, 
$N_{\gamma_{s,t}\vert l_{s,t,q}}\cap N_{\alpha_{n_s}\vert l_{s,t,q}}\! =\!\emptyset$. We replace 
$V_t$ with $N_{\gamma_{s,t}\vert l_{s,t,q}}$ and $U_s$ with $N_{\alpha_{n_s}\vert l_{s,t,q}}$ for the biggest $l_{s,t,q}$ with $q\!\leq\! l\! +\! 1$, which gives $V'_t$ and $U'_s$. We do this for all the possible $s,t\!\in\! 2^l$, which lead to $\tilde V_t$'s and $\tilde U_s$'s. We now choose $\alpha_{n_{t1}}\!\in\!\tilde V_t$. Then we choose disjoint clopen sets $U_{t0},U_{t1}$ with diameter at most $2^{-l-1}$ such that 
$\alpha_{n_{t\varepsilon}}\!\in\! U_{t\varepsilon}\!\subseteq\!\tilde U_t$ and satisfying (6)-(7).\bigskip

 There are three cases to conclude.\bigskip

\noindent\bf Case 1\rm\ $\big( (\varepsilon ,\alpha ),(\varepsilon ,\beta )\big)\!\notin\! E$ for each 
$\alpha\!\not=\!\beta\!\notin\!\mathbb{C}$ and each $\varepsilon\!\in\! 2$.\bigskip

 Then $E\! =\!\mathbb{E}^{\boratwo}_3$.\bigskip
 
\noindent\bf Case 2\rm\ $\big( (\varepsilon ,\alpha ),(\varepsilon ,\beta )\big)\!\in\! E$ and 
$\big( (1\! -\!\varepsilon ,\alpha ),(1\! -\!\varepsilon ,\beta )\big)\!\notin\! E$ for each 
$\alpha\!\not=\!\beta\!\notin\!\mathbb{C}$ and some $\varepsilon\!\in\! 2$.\bigskip

 Exchanging the first coordinate if necessary, we may assume that $\varepsilon\! =\! 1$. Then 
$E\! =\!\mathbb{E}^{\boratwo}_4$.\bigskip
 
\noindent\bf Case 3\rm\ $\big( (\varepsilon ,\alpha ),(\varepsilon ,\beta )\big)\!\in\! E$ for each 
$\alpha\!\not=\!\beta\!\notin\!\mathbb{C}$ and each $\varepsilon\!\in\! 2$.\bigskip

 Then $E\! =\!\mathbb{E}^{\boratwo}_5$.\hfill{$\square$}
 
\section{$\!\!\!\!\!\!$ Non-$\bormtwo$ equivalence relations}

\bf Notation.\rm\ We set $\mathbb{C}\! :=\!\mathbb{P}_f\! :=\!\{\alpha_n\mid n\!\in\!\omega\}$.

\vfill\eject

\noindent\bf Proof of Theorem \ref{main1} when ${\bf\Gamma}\! =\!\bormtwo$.\rm\ By Lemma \ref{basic}.(b), (a) and (b) cannot hold simultaneously. So assume that (a) does not hold. By the proof of Proposition \ref{injcontA}, we may assume that $X\! =\! 2^\omega$ and $\mathbb{C}$ is an equivalence class of $E$. By the Harrington-Kechris-Louveau theorem (see Theorem 13 in [M]), either there is $b\! :\! 2^\omega\!\rightarrow\! 2^\omega$ Borel with 
$E\! =\! (b\!\times\! b)^{-1}\big(\Delta (2^\omega )\big)$, or $(2^\omega ,\mathbb{E}_0)\sqsubseteq_c (2^\omega,E)$. In the latter case, the map $\phi\! :\! 2^{<\omega}\!\rightarrow\! 2^{<\omega}$ defined inductively by 
$\phi (\emptyset )\!:=\!\emptyset$, $\phi (s1)\! :=\!\phi (s)1\phi (s)$ and 
$$\phi (s0)\! :=\!\phi (s)0^{1+\vert\phi (s)\vert}$$ 
induces $f\! :\! 2^\omega\!\rightarrow\! 2^\omega$ injective continuous reducing $\mathbb{E}^{\bormtwo}_0$ to 
$\mathbb{E}_0$, showing that $\mathbb{E}^{\bormtwo}_0$ is below $E$. We apply Lemma \ref{subirrat} to a dense $G_\delta$ subset $G$ of $2^\omega$ contained in $\neg\mathbb{C}$ on which 
$b$ is continuous, so that we may assume that there is $b\! :\!\neg\mathbb{C}\!\rightarrow\! 2^\omega$ continuous such that $(\alpha ,\beta)\!\in\! E\Leftrightarrow b(\alpha )\! =\! b(\beta )$ if $\alpha ,\beta\!\notin\!\mathbb{C}$.\bigskip

\noindent\bf Case 1\rm\ $[\alpha ]_E$ is meager for each $\alpha\!\in\! 2^\omega$ (i.e., $b$ is nowhere dense-to-one).\bigskip

 Lemma \ref{lemGdelta} gives ${f\! :\! 2^\omega\!\rightarrow\! 2^\omega}$  injective continuous such that 
$\mathbb{C}\! =\! f^{-1}(\mathbb{C})$ and $b\big( f(\alpha )\big)\!\not=\! b\big( f(\beta )\big)$ if 
$\alpha\!\not=\!\beta\!\notin\!\mathbb{C}$. It remains to note that $f$ reduces $\mathbb{E}^{\bormtwo}_0$ to $E$.\bigskip

\noindent\bf Case 2\rm\ there is $\alpha\!\in\! 2^\omega$ such that $[\alpha ]_E$ is not meager.\bigskip

 We apply Lemma \ref{subirrat} to $G\! :=\! [\alpha ]_E$, which gives 
$f\! :\! 2^\omega\!\rightarrow\! 2^\omega$ injective continuous reducing 
$\mathbb{E}^{\bormtwo}_1$ to $E$.\hfill{$\square$}\bigskip

\noindent\bf Proof of Theorem \ref{main2} when ${\bf\Gamma}\! =\!\bormtwo$.\rm\ If 
$(\mathbb{X},\mathbb{E})\!\in\! {\mathcal B}^{\bormtwo}$, then 
$\mathbb{E}\!\notin\!\bormtwo$, so that (a) and (b) cannot hold simultaneously. Assume that (a) does not hold. By Theorem \ref{main1}, we may assume that the equivalence classes of $E$ are 
$\bormtwo$. By Theorem \ref{contB}, we may assume that $X\! =\! 2\!\times\! 2^\omega$, 
$$\Delta (2\!\times\! 2^\omega )\cup\big\{\big( (0,\alpha ),(1,\alpha )\big)\mid\alpha\!\in\!\mathbb{C}\big\}
\!\subseteq\! E$$ 
and $\big\{\big( (0,\alpha ),(1,\alpha )\big)\mid\alpha\!\notin\!\mathbb{C}\big\}\!\subseteq\!\neg E$.\bigskip

\noindent\bf Claim\it\ $(E_{(\varepsilon ,\gamma )})_{\varepsilon'}\cap\mathbb{C}$ is nowhere dense in 
$2^\omega$ for each $\gamma\!\in\! 2^\omega$ and each $\varepsilon ,\varepsilon'\!\in\! 2$.\rm\bigskip

 Indeed, we argue by contradiction, which gives $s\!\in\! 2^{<\omega}$ such that 
$N_s\!\subseteq\!\overline{(E_{(\varepsilon ,\gamma )})_{\varepsilon'}\cap\mathbb{C}}$. As 
$E_{(\varepsilon ,\gamma )}$ is $\bormtwo$, 
$(E_{(\varepsilon ,\gamma )})_{\varepsilon'}\cap N_s\!\setminus\!\mathbb{C}$ is comeager in $N_s$. Moreover, $(E_{(\varepsilon ,\gamma )})_{\varepsilon'}\cap\mathbb{C}\!\subseteq\! 
(E_{(\varepsilon ,\gamma )})_{1-\varepsilon'}$ since 
$$\big\{\big( (0,\alpha ),(1,\alpha )\big)\mid\alpha\!\in\!\mathbb{C}\big\}\!\subseteq\! E.$$
This implies that $(E_{(\varepsilon ,\gamma )})_{1-\varepsilon'}\cap N_s$ is a dense $\bormtwo$ subset of $N_s$, and is therefore comeager in $N_s$. Thus 
$(E_{(\varepsilon ,\gamma )})_0\cap (E_{(\varepsilon ,\gamma )})_1\cap N_s\!\setminus\!\mathbb{C}$ is comeager in $N_s$ and contains some $\beta$. Therefore $\big( (0,\beta ),(1,\beta )\big)$ is in $E$, which is absurd.
\hfill{$\diamond$}\bigskip

 The Sarbadhikari theorem gives an increasing sequence $(F_l)_{l\in\omega}$ of Borel relations on $2^\omega$ with closed nowhere dense vertical sections whose union contains 
$E\cap\big( (2\!\times\! 2^\omega )\!\times\! (2\!\times\!\mathbb{C})\big)$.\bigskip

 We will now prove that we may assume that 
$$E\!\subseteq\!\Delta (2\!\times\! 2^\omega )\cup
\big\{\big( (\varepsilon ,\alpha ),(1\! -\!\varepsilon ,\alpha )\big)\mid\varepsilon\!\in\! 2\wedge
\alpha\!\in\!\mathbb{C}\big\}\cup
\big\{\big( (\varepsilon ,\alpha ),(\varepsilon',\beta )\big)\mid\varepsilon ,\varepsilon'\!\in\! 2\wedge
\alpha\!\not=\!\beta\!\notin\!\mathbb{C}\big\} .$$

 We set $s\! :=\!\emptyset$, and construct $(n_t)_{t\in 2^{<\omega}}$, 
$(U_t)_{t\in 2^{<\omega}}$ satisfying (1)-(5) and the following:\bigskip

\leftline{$(6)~U_{t1}\cap\big(\{\alpha_{n_t}\}\cup\{\alpha_n\mid n\!\leq\!\vert t\vert\}\cup
\bigcup_{\varepsilon ,\varepsilon',\eta\in 2,s\in 2^l,\eta =0\vee s<_{\mbox{lex}}t}~
\overline{(E_{(\varepsilon ,\alpha_{n_{s\eta}})})_{\varepsilon'}\cap\mathbb{C}}$}\smallskip

\rightline{$\cup\bigcup_{\varepsilon ,\varepsilon'\in 2,s\in 2^l}~
\big( (F_l)_{(\varepsilon ,\alpha_{n_s})}\big)_{\varepsilon'}\big)\! =\!\emptyset$}\bigskip

\noindent Assume that this is done. If $\beta\!\notin\!\mathbb{C}$, then there is an infinite strictly increasing sequence $(l_k)_{k\in\omega}$ of natural numbers with $\beta\vert l_k\!\in\! O$, so that $f(\beta )\!\notin\!\mathbb{C}$ by Condition (6). Note that $\Big(\big(\varepsilon ,f(\alpha )\big) ,\big(\varepsilon',f(\alpha')\big)\Big)\!\notin\! E$ if 
$\alpha\!\not=\!\alpha'\!\in\!\mathbb{C}$, by Condition (6). Moreover, 
$\Big(\big(\varepsilon ,f(\alpha )\big) ,\big(\varepsilon',f(\beta )\big)\Big)\!\notin\! E$ if 
$\alpha\!\in\!\mathbb{C}$ and $\beta\!\notin\!\mathbb{C}$, by Condition (6). Thus we are done. Let us prove that the construction is possible. Assume that $(n_t)_{\vert t\vert\leq l}$ and $(U_t)_{\vert t\vert\leq l}$ satisfying (1)-(6) have been constructed, which is the case for $l\! =\! 0$. Let $t\!\in\! 2^l$. We define $n_{t1}$ by induction on $t$ with respect to the lexicographical ordering. As $(F_l)_{(\varepsilon ,\alpha )}$ is closed nowhere dense for each 
$(\varepsilon ,\alpha )\!\in\! 2\!\times\! 2^\omega$, $\big( (F_l)_{(\varepsilon ,\alpha )}\big)_{\varepsilon'}$ is a closed nowhere dense subset of $2^\omega$ for each $(\varepsilon ,\varepsilon',\alpha )$ in $2^2\!\times\! 2^\omega$. We choose $n_{t1}$ in such a way that\bigskip

\leftline{$\alpha_{n_{t1}}\!\in\! U_t\!\setminus\!\Big(\{\alpha_{n_t}\}\cup\{\alpha_n\mid n\!\leq\! l\}\ \cup
\bigcup_{\varepsilon ,\varepsilon'\in 2,s\in 2^l}~
\overline{(E_{(\varepsilon ,\alpha_{n_s})})_{\varepsilon'}\cap\mathbb{C}}\cup
\big( (F_l)_{(\varepsilon ,\alpha_{n_s})}\big)_{\varepsilon'}\cup$}\smallskip

\rightline{$\bigcup_{\varepsilon ,\varepsilon'\in 2,s\in 2^l,s<_{\mbox{lex}}t}~
\overline{(E_{(\varepsilon ,\alpha_{n_{s1}})})_{\varepsilon'}\cap\mathbb{C}}\Big)\mbox{,}$}\bigskip
 
\noindent We do this for each $t\!\in\! 2^l$, in the lexicographical ordering. Then we choose disjoint clopen sets $U_{t0},U_{t1}$ with diameter at most $2^{-l-1}$ such that $\alpha_{n_{t\varepsilon}}\!\in\! U_{t\varepsilon}\!\subseteq\! U_t$ and satisfying (6).\bigskip

 We now prove that we may assume that 
$$E\!\subseteq\!\Delta (2\!\times\! 2^\omega )\cup
\big\{\big( (\varepsilon ,\alpha ),(1\! -\!\varepsilon ,\alpha )\big)\mid\varepsilon\!\in\! 2\wedge
\alpha\!\in\!\mathbb{C}\big\}\cup
\big\{\big( (\varepsilon ,\alpha ),(\varepsilon ,\beta )\big)\mid\varepsilon\!\in\! 2\wedge
\alpha\!\not=\!\beta\!\notin\!\mathbb{C}\big\} .$$
Theorem 3.6 in [Lo1] gives a finer Polish topology $\sigma$ on $2\!\times\! 2^\omega$ such that 
$E\!\in\!\bormtwo\big( (2\!\times\! 2^\omega ,\sigma )^2\big)$ since the equivalence classes of $E$ are $\bormtwo$. Corollary 1.2 in [Ha-K-Lo] gives another Polish topology $\tau$ on 
$2\!\times\! 2^\omega$, finer than $\sigma$, such that 
$E\!\in\!\bormone\big( (2\!\times\! 2^\omega ,\tau )^2\big)$. By 8.38 in [K1], there is a dense 
$G_\delta$ subset of $\neg\mathbb{C}$ on which $\tau$ and the usual topology coincide. This shows that we may assume that $E\cap\big(2\!\times\! (\neg\mathbb{C})\big)^2$ is closed in 
$\big(2\!\times\! (\neg\mathbb{C})\big)^2$, which gives a closed relation $F$ on $2\!\times\! 2^\omega$ with 
$$E\cap\big(2\!\times\! (\neg\mathbb{C})\big)^2\! =\! F\cap\big(2\!\times\! (\neg\mathbb{C})\big)^2.$$
Fix $\varepsilon\!\in\! 2$. Note that $F\cap\big\{\big( (\varepsilon ,\alpha ),(1\! -\!\varepsilon ,\alpha )\big)\mid\alpha\!\in\!\mathbb{C}\big\}$ is nowhere dense in 
$$\big\{\big( (\varepsilon ,\alpha ),(1\! -\!\varepsilon ,\alpha )\big)\mid\alpha\!\in\! 2^\omega\big\} .$$ 
Indeed, we argue by contradiction, which gives $s\!\in\! 2^{<\omega}$ such that 
$$\big\{\big( (\varepsilon ,\alpha ),(1\! -\!\varepsilon ,\alpha )\big)\mid\alpha\!\in\! N_s\big\}\!\subseteq\!
\overline{F\cap\big\{\big( (\varepsilon ,\alpha ),(1\! -\!\varepsilon ,\alpha )\big)\mid\alpha\!\in\!\mathbb{C}\big\}}
\!\subseteq\! F\mbox{,}$$ 
and $\alpha\!\notin\!\mathbb{C}$ such that $\big( (\varepsilon ,\alpha ),(1\! -\!\varepsilon ,\alpha )\big)$ is in $F$, and thus in $E$, which cannot be. This gives $t\!\in\! O$ such that 
$\big( (\varepsilon ,t0^\infty ),(1\! -\!\varepsilon ,t0^\infty )\big)\!\notin\! F$, and $l\!\in\!\omega$ with 
$\big( (\{\varepsilon\}\!\times\! N_{t0^l})\!\times\! (\{1\! -\!\varepsilon\}\!\times\! N_{t0^l})\big)\cap F
\! =\!\emptyset$. So we are done.\bigskip

 The previous point shows that we may assume that 
$E\cap\big(\{\varepsilon\}\!\times\! (\neg\mathbb{C})\big)^2$ is a closed equivalence relation on 
$\{\varepsilon\}\!\times\! (\neg\mathbb{C})$ for each $\varepsilon\!\in\! 2$. By 18.D in [K1], there is a 
$\sigma (\ana )$-measurable map 
$$S\! :\!\{\varepsilon\}\!\times\! (\neg\mathbb{C})\!\rightarrow\!\{\varepsilon\}\!\times\! (\neg\mathbb{C})$$ 
such that $S(\varepsilon ,\alpha )\! =\! S(\varepsilon ,\beta )~E~(\varepsilon ,\alpha )$ if 
$\big( (\varepsilon ,\alpha ),(\varepsilon ,\beta )\big)\!\in\! E\cap\big(\{\varepsilon\}\!\times\! (\neg\mathbb{C})\big)^2$. By 8.38 and 29.D in [K1], there is a dense $G_\delta$ subset $G$ of $\neg\mathbb{C}$ such that the restriction of 
$S$ to $\{\varepsilon\}\!\times\! G$ is continuous. So we may assume that there is 
$b_\varepsilon\! :\!\neg\mathbb{C}\!\rightarrow\! 2^\omega$ continuous such that 
$\big( (\varepsilon ,\alpha ),(\varepsilon ,\beta )\big)\!\in\! E\Leftrightarrow b_\varepsilon (\alpha )\! =\! b_\varepsilon (\beta )$ if $\alpha ,\beta\!\notin\!\mathbb{C}$. Assume first that $b_\varepsilon$ is nowhere dense-to-one. By Lemma \ref{lemGdelta}, there is ${f\! :\! 2^\omega\!\rightarrow\! 2^\omega}$  injective continuous such that 
$\mathbb{C}\! =\! f^{-1}(\mathbb{C})$ and $b_\varepsilon\big( f(\alpha )\big)\!\not=\! b_\varepsilon\big( f(\beta )\big)$ if $\alpha\!\not=\!\beta\!\notin\!\mathbb{C}$. This implies that $f$ reduces $\mathbb{E}^{\bormtwo}_3$ to $E$ if both $b_0$ and $b_1$ are nowhere dense-to-one. If $b_0$ is not nowhere dense-to-one and $b_1$ is nowhere dense-to-one, then using Lemma \ref{subirrat} we see that $\mathbb{E}^{\bormtwo}_4$ is reducible to $E$. This is also the case if $b_1$ is not nowhere dense-to-one and $b_0$ is nowhere dense-to-one, since we can exchange the first coordinate. If neither $b_0$, nor $b_1$ is nowhere dense-to-one, then $\mathbb{E}^{\bormtwo}_5$ is reducible to 
$E$, similarly.\hfill{$\square$}

\section{$\!\!\!\!\!\!$ Equivalence relations with countably many classes}

\subsection{$\!\!\!\!\!\!$ Non-$\boraxi$ equivalence relations with countably many classes}\indent

 If $\xi\!\geq\! 2$ is a countable ordinal, then Lemma \ref{true} provides 
$\mathbb{C}\!\in\!\bormxi (2^\omega)\!\setminus\!\boraxi$. Subsection \ref{SectionHur+} provides a partition 
$(\mathbb{C}_n)_{n\in\omega}$ of $\neg\mathbb{C}$ into $\borxi$ subsets of $2^\omega$, which allows to define an equivalence relation on $2^\omega$ by 
$\mathbb{E}^{\boraxi}_2\! :=\!\mathbb{C}^2\cup\bigcup_{n\in\omega}~\mathbb{C}_n^2$, as in the introduction.\bigskip

\noindent\bf Proof of Theorem \ref{c1dencls}.\rm\ By Lemma \ref{basic}.(a), the equivalence classes of $E$ are $\boraxi$ if $E$ is a $\boraxi$ subset of $X^2$. The converse comes from the fact that $E$ is the countable union of the square of its equivalence classes. By Lemma \ref{basic}, (a) and (b) cannot hold simultaneously. By Theorem \ref{main1}, we may assume that $\xi\!\geq\! 3$.\bigskip

 By Proposition \ref{injcontA}, we may assume that $X\! =\! 2^\omega$ and $\mathbb{C}$ is an equivalence class of the Borel relation $E$. As $E$ has countably many classes, we can write 
$\neg\mathbb{C}\! =\!\bigcup_{n\in I}~D_n$, where the $D_n$ are distinct $E$-classes and $I$ is countable and nonempty.\bigskip

 If there is $n$ such that the Borel set $D_n$ is not separable from the Borel set $\mathbb{C}$ by a 
$\bormxi$ set, then Theorem \ref{Lo-SR} gives $j\! :\! 2^\omega\!\rightarrow\! X$  injective continuous such that $\mathbb{C}\!\subseteq\! j^{-1}(\mathbb{C})$ and $\neg\mathbb{C}\!\subseteq\! j^{-1}(D_n)$. This implies that $(2^\omega ,\mathbb{E}^{\boraxi}_1)\sqsubseteq_c(X,E)$.\bigskip

 If the $D_n$'s are separable from $\mathbb{C}$ by a $\bormxi$ set, then they are separable from $\mathbb{C}$ by a $\borxi$ set. In particular, $I$ is infinite and we may assume that $I\! =\!\omega$. Theorem 
\ref{Hur+} provides $\phi\! :\!\omega\!\rightarrow\!\omega$ injective and 
$f\! :\! 2^\omega\!\rightarrow\! X$ injective continuous such that $\mathbb{C}\!\subseteq\! f^{-1}(\mathbb{C})$ and 
$\mathbb{C}_n\!\subseteq\! f^{-1}(D_{\phi (n)})$ for each $n\!\in\!\omega$. Note that $f$ reduces 
$\mathbb{E}^{\boraxi}_2$ to $E$ as desired.\hfill{$\square$}

\subsection{$\!\!\!\!\!\!$ Non-$\bormxi$ equivalence relations with countably many classes}
 
\bf Proof of Theorem \ref{C1denclPi}.\rm\ By Theorem \ref{prope}, (a) and (b) cannot hold simultaneously. By Proposition \ref{injcontA}, we may assume that $X\! =\!\mathbb{K}$ and $\mathbb{C}$ is an equivalence class of the Borel relation $E$. As $E$ has countably many classes, we can write 
$\neg\mathbb{C}\! =\!\bigcup_{n\in I}~C_n$, where the $C_n$ are distinct $E$-classes and $I$ is countable and nonempty. As $\mathbb{C}\!\notin\!\bormxi$, there is $n$ such that $C_n$ is not separable from 
$\mathbb{C}$ by a $\boraxi$ set. As $\mathbb{C}$ and $C_n$ are Borel, Theorem \ref{Lo-SR} gives 
$j\! :\!\mathbb{K}\!\rightarrow\! X$  injective continuous such that $\mathbb{C}\!\subseteq\! j^{-1}(\mathbb{C})$ and $\neg\mathbb{C}\!\subseteq\! j^{-1}(C_n)$. This implies that 
$(\mathbb{K},\mathbb{E}^{\bormxi}_1)\sqsubseteq_c(X,E)$ as desired.\hfill{$\square$}\bigskip

 In order to finish the study of Borel equivalence relations with countably many classes, it remains to characterize those which are not $\bormxi$ if $\xi\!\geq\! 3$. Lemma \ref{true} provides 
$\mathbb{C}\!\in\!\boraxi (2^\omega)\!\setminus\!\bormxi$. Subsection \ref{SectionHur+} provides a partition 
$(\mathbb{C}_n)_{n\in\omega}$ of $\mathbb{C}$ into $\borxi$ subsets of $2^\omega$, which allows to define an equivalence relation $\mathbb{E}^{\bormxi}_8$ on $2\!\times\! 2^\omega$ as in the introduction.\bigskip

\noindent\bf Notation.\rm\ Let $E$ be an equivalence relation on $2\!\times\! 2^\omega$. We set, for 
$\varepsilon ,\eta\!\in\! 2$, 
$$E_{\varepsilon ,\eta}\! :=\!\big\{ (\alpha ,\beta )\!\in\! 2^\omega\!\times\! 2^\omega\mid
\big( (\varepsilon ,\alpha ),(\eta ,\beta )\big)\!\in\! E\big\} .$$ 
Note that $E_{\varepsilon ,\varepsilon}$ is an equivalence relation on $2^\omega$.\bigskip

\noindent\bf Proof of Theorem \ref{pixidencl}.\rm\ By Theorem \ref{main2}, we may assume that 
$\xi\!\geq\! 3$. If $n\!\in\!\{ 1,8\}$, then $\mathbb{E}^{\bormxi}_n\!\notin\!\bormxi$, so that (a) and (b) cannot hold simultaneously. Assume that (a) does not hold. By Theorem \ref{C1denclPi}, we may assume that $E$ has $\bormxi$ classes. By Theorem \ref{contB}, in order to prove that 
$(2\!\times\! 2^\omega ,\mathbb{E}^{\bormxi}_8)\sqsubseteq_c(X,E)$, we may assume that 
$X\! =\! 2\!\times\! 2^\omega$, 
$\big\{\big( (0,\alpha ),(1,\alpha )\big)\mid\alpha\!\in\!\mathbb{C}\big\}\!\subseteq\! E$ 
and $\big\{\big( (0,\alpha ),(1,\alpha )\big)\mid\alpha\!\notin\!\mathbb{C}\big\}\!\subseteq\!\neg E$.\bigskip

 Note that $E_{\varepsilon ,\varepsilon}$ has countably many $\bormxi$ classes, for each 
$\varepsilon\!\in\! 2$, since the map $\alpha\!\mapsto\! (\varepsilon ,\alpha )$ reduces 
$E_{\varepsilon ,\varepsilon}$ to $E$. Consequently, we can write 
$\neg\mathbb{C}\! =\!\bigcup_{n\in\omega}~D^\varepsilon_n$, where the $D^\varepsilon_n$'s are 
$\bormxi$ and contained in distinct $E_{\varepsilon ,\varepsilon}$-classes. Note that there is 
$n\!\in\!\omega$ such that $D^\varepsilon_n$ is not separable from $\mathbb{C}$ by a $\boraxi$ set. Theorem \ref{Lo-SR} gives $g\! :\! 2^\omega\!\rightarrow\! 2^\omega$ injective continuous such that $\neg\mathbb{C}\!\subseteq\! g^{-1}(D^\varepsilon_n)$ and 
$\mathbb{C}\!\subseteq\! g^{-1}(\mathbb{C})$. So, replacing $E$ with 
$\big( (\mbox{Id}_2\!\times\! g)\!\times\! (\mbox{Id}_2\!\times\! g)\big)^{^-1}(E)$ if necessary, we may assume that $\neg\mathbb{C}$ is contained in a single $E_{\varepsilon ,\varepsilon}$-class 
$K_\varepsilon$, for each $\varepsilon\!\in\! 2$.\bigskip 

 Let us prove that $\neg\mathbb{C}$ is separable from $K_0\cap\mathbb{C}$ by a $\boraxi$ set, say $S$. We argue by contradiction. Theorem \ref{Lo-SR} gives $h\! :\! 2^\omega\!\rightarrow\! 2^\omega$ injective continuous such that $\neg\mathbb{C}\!\subseteq\! h^{-1}(\neg\mathbb{C})$ and 
 $\mathbb{C}\!\subseteq\! h^{-1}(K_0\cap\mathbb{C})$. We set 
$E'\! :=\!\big( (\mbox{Id}_2\!\times\! h)\!\times\! (\mbox{Id}_2\!\times\! h)\big)^{^-1}(E)$, so that $E'$ is a Borel equivalence relation on $2\!\times\! 2^\omega$ with countably many $\bormxi$ classes. Moreover, 
$\mathbb{C}^2\!\subseteq\! E'_{0,0}\cap\mathbb{C}^2\! =\! E'_{1,1}\cap\mathbb{C}^2$, by Proposition \ref{foursame}. So $\mathbb{C}$ is contained in an $E'_{1,1}$-class $C'$, which has to be $\bormxi$ as above. So let $\beta\!\in\! C'\!\setminus\!\mathbb{C}$, and $\alpha\!\in\!\mathbb{C}$. Then $(\alpha ,\beta )\!\in\! E'_{1,1}$, $(\beta ,\alpha )\!\in\! E'_{0,0}$, and 
$(\alpha ,\alpha )\!\in\! E'_{0,1}$, so that $(\beta ,\beta )\!\in\! E'_{0,1}$, which is absurd.\bigskip 

 Let us prove that $\neg\mathbb{C}$ is not separable from $\neg K_0$ by a $\boraxi$ set. We argue by contradiction, which gives $S'\!\in\!\boraxi$. Note that $2^\omega\! =\!\mathbb{C}\cup S'$ is a covering into 
$\boraxi$ sets. The reduction property of $\boraxi$ gives $\Delta\!\in\!\borxi$ with $\Delta\!\subseteq\! S'$ and 
$\neg\Delta\!\subseteq\!\mathbb{C}$ (see 22.16 in [K1]). Then $\neg\mathbb{C}\!\subseteq\!\Delta\!\subseteq\! K_0$, so that $\neg\mathbb{C}\! =\!\Delta\cap S\!\in\!\boraxi$, which is absurd.

\vfill\eject

 Theorem \ref{Lo-SR} gives $k\! :\! 2^\omega\!\rightarrow\! 2^\omega$ injective continuous such that 
$\mathbb{C}\!\subseteq\! k^{-1}(\neg K_0)$ and $\neg\mathbb{C}\!\subseteq\! k^{-1}(\neg\mathbb{C})$. So, replacing $E$ with $\big( (\mbox{Id}_2\!\times\! k)\!\times\! (\mbox{Id}_2\!\times\! k)\big)^{^-1}(E)$ if necessary, we may assume that $\neg\mathbb{C}$ is an $E_{0,0}$-class.\bigskip

 As $E_{0,0}$ has countably many $\bormxi$ classes, we can write $\mathbb{C}\! =\!\bigcup_{n\in\omega}~D_n$, where the $D_n$'s are distinct $\bormxi$ classes for $E_{0,0}$. Theorem \ref{Hur+} provides $\phi\! :\!\omega\!\rightarrow\!\omega$ injective and $k\! :\! 2^\omega\!\rightarrow\! 2^\omega$ injective continuous such that 
$\neg\mathbb{C}\!\subseteq\! k^{-1}(\neg\mathbb{C})$ and $\mathbb{C}_n\!\subseteq\! k^{-1}(D_{\phi (n)})$ for each 
$n\!\in\!\omega$. Replacing $E$ with 
$$\big( (\mbox{Id}_2\!\times\! k)\!\times\! (\mbox{Id}_2\!\times\! k)\big)^{^-1}(E)$$ 
if necessary, we consequently may assume that\bigskip
 
\noindent - $E_{0,0}\! =\! (\neg\mathbb{C})^2\cup\bigcup_{n\in\omega}~\mathbb{C}_n^2$,\smallskip
 
\noindent - $\big\{\big( (0,\alpha ),(1,\alpha )\big)\mid\alpha\!\in\!\mathbb{C}\big\}\!\subseteq\! E$,\smallskip
 
\noindent - $\big\{\big( (0,\alpha ),(1,\alpha )\big)\mid\alpha\!\notin\!\mathbb{C}\big\}\!\subseteq\!\neg E$,\smallskip
 
\noindent - $\neg\mathbb{C}$ is contained in an $E_{1,1}$-class.\bigskip

 Proposition \ref{foursame} shows that if $\varepsilon ,\eta\!\in\! 2$ and $\alpha ,\beta\!\in\!\mathbb{C}$, then 
$(\alpha ,\beta )\!\in\! E_{\varepsilon ,\eta}$ is equivalent to $(\alpha ,\beta )\!\in\! E_{0,0}$ (and 
$(\alpha ,\beta )\!\in\!\bigcup_{n\in\omega}~\mathbb{C}_n^2$).\bigskip

 Note that $E_{\varepsilon ,1-\varepsilon}\cap (\neg\mathbb{C})^2\! =\!\emptyset$. Indeed, we argue by contradiction and we may assume that $\varepsilon\! =\! 0$, which gives $\alpha ,\beta\!\notin\!\mathbb{C}$ such that 
$(\alpha ,\beta )\!\in\! E_{0,1}$. As $(\alpha ,\beta )\!\in\! E_{0,0}$, $(\beta ,\beta )\!\in\! E_{0,1}$, which is absurd.\bigskip
  
 We set, for $p\!\in\!\omega$, 
$B_{p+1}\! :=\!\{\beta\!\notin\!\mathbb{C}\mid\exists\alpha\!\in\!\mathbb{C}_p~~(\alpha ,\beta )\!\in\! E_{0,1}\}$. Note that $B_{p+1}$ is analytic. In fact, if $\beta\!\in\! B_{p+1}$ with witness $\alpha$ and $\gamma\!\in\!\mathbb{C}_p$, then $(\gamma ,\alpha )\!\in\! E_{0,0}$, so that 
$(\gamma ,\beta )\!\in\! E_{0,1}$ and 
$$B_{p+1}\! :=\!\{\beta\!\notin\!\mathbb{C}\mid\forall\gamma\!\in\!\mathbb{C}_p~~
(\gamma ,\beta )\!\in\! E_{0,1}\}$$ 
is also co-analytic and thus Borel. Moreover, the $B_{p+1}$'s are pairwise disjoint since two different 
$\mathbb{C}_p$'s are not $E_{0,0}$-related. We set 
$B_0\! :=\! (\neg\mathbb{C})\!\setminus\! (\bigcup_{p\in\omega}~B_{p+1})$. Then $(B_p)_{p\in\omega}$ is a partition of $\neg\mathbb{C}$ into Borel sets. Note that there is $p$ such that $B_p$ is not separable from $\mathbb{C}$ by a $\boraxi$ set. Theorem \ref{Hur+} provides $\psi\! :\!\omega\!\rightarrow\!\omega$ injective and 
$l\! :\! 2^\omega\!\rightarrow\! 2^\omega$ injective continuous such that $\neg\mathbb{C}\!\subseteq\! l^{-1}(B_p)$ and $\mathbb{C}_n\!\subseteq\! l^{-1}(\mathbb{C}_{\psi (n)})$ for each $n\!\in\!\omega$. So, replacing $E$ with 
$\big( (\mbox{Id}_2\!\times\! l)\!\times\! (\mbox{Id}_2\!\times\! l)\big)^{^-1}(E)$ if necessary, we may assume that 
$(\alpha ,\beta )\!\notin\! E_{0,1}$ if $\beta\!\notin\!\mathbb{C}$, $\alpha\!\in\!\mathbb{C}_n$ and $n\!\not=\! p$. As 
$\neg\mathbb{C}$ is not separable from $\bigcup_{n\not= p}~\mathbb{C}_n$ by a $\boraxi$ set, we can apply again Theorem \ref{Hur+}  to see that we may assume that $(\alpha ,\beta )\!\notin\! E_{0,1}$ if 
$\beta\!\notin\!\mathbb{C}$, $\alpha\!\in\!\mathbb{C}_n$ and $n\!\in\!\omega$. By symmetry, 
$(\alpha ,\beta )\!\notin\! E_{1,0}$ if $\alpha\!\notin\!\mathbb{C}$, $\beta\!\in\!\mathbb{C}_n$ and $n\!\in\!\omega$. Similarly, we may assume that $(\alpha ,\beta )\!\notin\! E_{1,0}$ if $\beta\!\notin\!\mathbb{C}$, 
$\alpha\!\in\! \mathbb{C}_n$ and $n\!\in\!\omega$. By symmetry, $(\alpha ,\beta )\!\notin\! E_{0,1}$ if 
$\alpha\!\notin\!\mathbb{C}$, $\beta\!\in\!\mathbb{C}_n$ and $n\!\in\!\omega$. Similarly again, we may assume that $(\alpha ,\beta )\!\notin\! E_{1,1}$ if $\alpha\!\notin\!\mathbb{C}$, $\beta\!\in\!\mathbb{C}_n$ and $n\!\in\!\omega$. By symmetry, $(\alpha ,\beta )\!\notin\! E_{1,1}$ if $\beta\!\notin\!\mathbb{C}$, $\alpha\!\in\!\mathbb{C}_n$ and 
$n\!\in\!\omega$. So we proved that we may assume that $E\! =\!\mathbb{E}^{\bormxi}_8$, i.e., 
$(\mathbb{H},\mathbb{E}^{\bormxi}_8)\sqsubseteq_c(X,E)$.\hfill{$\square$}
  
\section{$\!\!\!\!\!\!$ Borel equivalence relations with $F_\sigma$ classes}

\bf Proof of Theorem \ref{ctble}.\rm\ By Theorem \ref{prope}, (a) and (b) cannot hold simultaneously. So assume that (a) does not hold. As $E$ has $F_\sigma$ classes, its sections are in $\bf\Gamma$. By Theorem \ref{contB}, we may assume that $X\! =\!\mathbb{H}$, 
$\big\{\big( (0,\alpha ),(1,\alpha )\big)\mid\alpha\!\in\!\mathbb{C}\big\}\!\subseteq\! E$, and 
$\big\{\big( (0,\alpha ),(1,\alpha )\big)\mid\alpha\!\notin\!\mathbb{C}\big\}\!\subseteq\!\neg E$.\bigskip

 Recall that $E_{\varepsilon ,\varepsilon}$ is a Borel equivalence relation on $2^\omega$ with $F_\sigma$ classes. In order to simplify the notation, we may assume by relativization that 
$\xi\! :=\!\mbox{rk}({\bf\Gamma})\! <\!\omega_1^{\mbox{CK}}$ and $\mathbb{C},E\!\in\!\Borel$. We partly follow the proof of Silver's theorem (see [S]) given in [G]. So we set 
$$W\! :=\!\{\alpha\!\in\! 2^\omega\mid\exists U\!\in\!\Borel (2^\omega )~~
\alpha\!\in\! U\!\subseteq\! [\alpha ]_{E_{0,0}}\}\mbox{,}$$ 
and $V\! :=\! 2^\omega\!\setminus\! W$. The proof of Theorem 5.3.5 in [G] shows that $V\!\in\!\Ana$, and that $E_{0,0}\cap V^2$ is ${\it\Sigma}_{2^\omega}^2$-meager in $V^2$. Note also that $W$ contains $\Borel\cap 2^\omega$. As $\Borel (2^\omega )$ is countable, we can find a countable set $I$ and a sequence $(U_i)_{i\in I}$ of nonempty $\Borel$ sets each contained in a single $E_{0,0}$-class such that $W$ is contained in the $F_\sigma$ set $S\! :=\!\bigcup_{i\in I}~[U_i]_{E_{0,0}}$, where $[U_i]_{E_{0,0}}\! :=\!
\{\alpha\!\in\! 2^\omega\mid\exists\beta\!\in\! U_i~~(\alpha ,\beta )\!\in\! E_{0,0}\}$ is $\Ana$. Pick 
$\alpha_i\!\in\! U_i$ for each $i\!\in\! I$, so that $[U_i]_{E_{0,0}}\! =\! [\alpha_i]_{E_{0,0}}$ and $S$ is the disjoint union of the $[\alpha_i]_{E_{0,0}}$'s.\bigskip

 Let us prove that $V\cap\mathbb{C}$ is not separable from $V\!\setminus\!\mathbb{C}$ by a set in $\bf\Gamma$. We argue by contradiction, so that $\mathbb{C}\!\setminus\! S$ is also separable from $\neg (S\cup\mathbb{C})$ by a set in $\bf\Gamma$. As $\mathbb{C}\!\notin\! {\bf\Gamma}$ and $S\!\in\! F_\sigma\!\subseteq\! {\bf\Gamma}$, $S\cap\mathbb{C}$ is not separable from 
$S\!\setminus\!\mathbb{C}$ by a set in $\bf\Gamma$. This gives $i\!\in\! I$ such that 
$[\alpha_i]_{E_{0,0}}\cap\mathbb{C}$ is not separable from 
$[\alpha_i]_{E_{0,0}}\!\setminus\!\mathbb{C}$ by a set in $\bf\Gamma$. In particular, there is 
$\alpha\!\in\! [\alpha_i]_{E_{0,0}}\cap\mathbb{C}$. If 
$\beta\!\in\! [\alpha ]_{E_{0,0}}\cap\mathbb{C}$, then $\big( (0,\beta ),(1,\beta )\big)\!\in\! E$. Thus 
$\{ 1\}\!\times\! ([\alpha ]_{E_{0,0}}\cap\mathbb{C})$ is contained in the $F_\sigma$ set 
$[(0,\alpha )]_E\cap (\{ 1\}\!\times\! [\alpha ]_{E_{0,0}})$. This gives 
$\gamma\!\in\! [\alpha ]_{E_{0,0}}\!\setminus\!\mathbb{C}$ such that 
$\big( (0,\alpha ),(1,\gamma )\big)\!\in\! E$. As $\big( (0,\alpha ),(0,\gamma )\big)\!\in\! E$, 
$\big( (0,\gamma ),(1,\gamma )\big)\!\in\! E$, which is absurd.\bigskip

 Theorem \ref{presC-R} provides $f\! :\! 2^\omega\!\rightarrow\! 2^\omega$ injective continuous such that $\mathbb{C}\! =\! f^{-1}(\mathbb{C})$ and $\big( f(\alpha ),f(\beta )\big)$ is not in 
$E_{0,0}$ if $\alpha\!\not=\!\beta$. This shows that we may assume that $E$ coincides with 
$\mathbb{E}^{\bf\Gamma}_3$ on $(\{ 0\}\!\times\! 2^\omega )^2$.\bigskip

 Similarly, we may assume that $E$ coincides with $\mathbb{E}^{\bf\Gamma}_3$ on $(\{ 1\}\!\times\! 2^\omega )^2$. By Proposition \ref{foursame}, $E$ coincides with $\mathbb{E}^{\bf\Gamma}_3$ on $(\{\varepsilon\}\!\times\!\mathbb{C})\!\times\! (\{\eta\}\!\times\!\mathbb{C})$ for each 
$\varepsilon ,\eta\!\in\! 2$. Pick $\alpha ,\beta ,\gamma\!\in\! 2^\omega$. If both 
$\big( (0,\alpha ),(1,\beta )\big)$ and $\big( (0,\alpha ),(1,\gamma )\big)$ are in $E$, then 
$\beta\! =\!\gamma$. Similarly, if 
$\big( (0,\beta ),(1,\alpha )\big) ,\big( (0,\gamma ),(1,\alpha )\big)\!\in\! E$, then 
$\beta\! =\!\gamma$. This shows that $E$ coincides with $\mathbb{E}^{\bf\Gamma}_3$ on 
$(\{\varepsilon\}\!\times\!\mathbb{C})\!\times\!
\big(\{ 1\! -\!\varepsilon\}\!\times\! (\neg\mathbb{C})\big)$ and 
$\big(\{\varepsilon\}\!\times\! (\neg\mathbb{C})\big)\!\times\!
(\{ 1\! -\!\varepsilon\}\!\times\!\mathbb{C})$ for each $\varepsilon\!\in\! 2$, and also that 
$E\cap\Big(\big(\{\varepsilon\}\!\times\! (\neg\mathbb{C})\big)\!\times\!
\big(\{ 1\! -\!\varepsilon\}\!\times\! (\neg\mathbb{C})\big)\Big)$ is the graph of a Borel injection. In particular, $E$ is countable. We set $R\! :=\!\bigcup_{\varepsilon ,\eta\in 2}~E_{\varepsilon ,\eta}$. Note that $R'$ is a locally countable relation on $2^\omega$. Corollary \ref{avoidctble} provides 
$l\! :\! 2^\omega\!\rightarrow\! 2^\omega$ injective continuous such that 
$\mathbb{C}\! =\! l^{-1}(\mathbb{C})$ and $\big( l(\alpha ),l(\beta )\big)\!\notin\! R'$ if 
$\alpha\!\not=\!\beta$. So we may assume that $E$ coincides with 
$\mathbb{E}^{\bf\Gamma}_3$.\hfill{$\square$}

\vfill\eject
  
\section{$\!\!\!\!\!\!$ References}

\noindent [A-K]\ \ S. Adams and A. S. Kechris, Linear algebraic groups and countable Borel equivalence relations,~\it J. Amer. Math. Soc.\rm ~13 (2000), 909-943

\noindent [C-L-M]\ \ J. D. Clemens, D. Lecomte and B. D. Miller, Essential countability of treeable equivalence relations,~\it Adv. Math.\ \rm 265 (2014), 1-31

\noindent [D-SR]\ \ G. Debs and J. Saint Raymond, Borel liftings of Borel sets: 
some decidable and undecidable statements,~\it Mem. Amer. Math. Soc.\rm ~187, 876 (2007)

\noindent [G]\ \ S. Gao,~\it Invariant Descriptive Set Theory,~\rm Pure and Applied Mathematics, A Series of Monographs and Textbooks, 293, Taylor and Francis Group, 2009

\noindent [Ha-K-Lo]\ \ L. A. Harrington, A. S. Kechris and A. Louveau, A Glimm-Effros 
dichotomy for Borel equivalence relations,~\it J. Amer. Math. Soc.\rm ~3 (1990), 903-928

\noindent [H-K]\ \ G. Hjorth and A. S. Kechris, New dichotomies for Borel equivalence relations,~\it Bull. Symbolic Logic\rm ~3, 3 (1997), 329-346

\noindent [J-K-Lo]\ \ S. Jackson, A. S. Kechris and A. Louveau, Countable Borel equivalence relations,~\it J. Math. Log. \rm ~2, 1 (2002), 1-80

\noindent [K1]\ \ A. S. Kechris,~\it Classical Descriptive Set Theory,~\rm Springer-Verlag, 1995

\noindent [K2]\ \ A. S. Kechris,~\it The theory of countable Borel equivalence relations,~\rm preprint, 2018 (see the author's webpage at http://www.math.caltech.edu/)

\noindent [L1]\ \ D. Lecomte, A dichotomy characterizing analytic graphs of uncountable Borel chromatic number in any dimension,~\it Trans. Amer. Math. Soc.\rm~361 (2009), 4181-4193

\noindent [L2]\ \ D. Lecomte, How can we recognize potentially $\bormxi$ subsets of the plane?,~\it  J. Math. Log.\ \rm  9, 1 (2009), 39-62

\noindent [L3]\ \ D. Lecomte, Potential Wadge classes,~\it\ Mem. Amer. Math. Soc.,\rm ~221, 1038 (2013)

\noindent [L4]\ \ D. Lecomte, A separation result for countable unions of Borel rectangles,~\it\ arXiv:1708.00642\rm

\noindent [Lo1]\ \ A. Louveau, A separation theorem for $\Ana$ sets,\ \it Trans. 
Amer. Math. Soc.\ \rm 260 (1980), 363-378

\noindent [Lo2]\ \ A. Louveau, Ensembles analytiques et bor\'eliens dans les 
espaces produit,~\it Ast\'erisque (S. M. F.)\ \rm 78 (1980)

\noindent [Lo-SR]\ \ A. Louveau and J. Saint Raymond, Borel classes and closed games: 
Wadge-type and Hurewicz-type results,\ \it Trans. Amer. Math. Soc.\ \rm 304 (1987), 431-467

\noindent [M]\ \ B. D. Miller, Forceless, ineffective, powerless proofs of descriptive dichotomy theorems, Lecture 3, The Harrington-Kechris-Louveau theorem,~\it\ draft\rm ~(see the author's webpage at the address https://dl.dropboxusercontent.com/u/47430894/Web/otherwork/paristhree.pdf)

\noindent [Mos]\ \ Y. N. Moschovakis,~\it Descriptive set theory,~\rm North-Holland, 1980

\noindent [R]\ \ C. Rosendal, Cofinal families of Borel equivalence relations and quasiorders,~\it\ J. Symbolic Logic\rm ~70, 4 (2005), 1325-1340

\noindent [S]\ \ J. H. Silver, Counting the number of equivalence classes of Borel and coanalytic
equivalence relations,~\it Ann. Math. Logic\ \rm 18, 1 (1980), 1-28

\noindent [Sr]\ \ S. M. Srivastava,~\it A course on Borel sets,~\rm Graduate Texts in Mathematics, 180, Springer-Verlag, New-York, 1998
 
\end{document}